\documentclass[11pt]{article}
 \usepackage{graphicx,amsmath,amsfonts,amssymb,color,mathrsfs}

\usepackage{epsfig}

 \newcommand{\QED}{\hfill \thicklines \framebox(6.6,6.6)[l]{}}
 \newenvironment{proof}{\noindent {\bf \sc Proof.} \rm}{\QED}

\numberwithin{equation}{section}

 \newtheorem{assumption}{Assumption}
 \newtheorem{convention}{Convention}
 \newtheorem{theorem}{Theorem}[section]
 \newtheorem{lemma}{Lemma}[section]
 \newtheorem{proposition}{Proposition}[section]
 
 \newtheorem{remark}{Remark}[section]

 \newtheorem{corollary}{Corollary}[section]
 \newtheorem{definition}{Definition}[section]

\newcommand{\eqnb}{\begin{eqnarray*}}
\newcommand{\eqne}{\end{eqnarray*}}

\renewcommand\Im{\operatorname{Im}}

\title{\Large \bf A Kernel Method for Exact Tail Asymptotics \\--- Random Walks in the Quarter
Plane\footnotetext{$\,$ \hspace{-6.2ex} $^{*}$ Postal address:
Department of Mathematics, Mount Saint Vincent University,
Halifax, NS Canada B3M 2J6
 \newline
$^{**}$ Postal address: School of Mathematics and Statistics,
Carleton University, Ottawa, ON Canada K1S 5B6} \\
\rm{(In memory of Dr.~Philippe Flajolet)} }
\author{Hui Li$^{*}$, Mount Saint Vincent University \\
Yiqiang Q. Zhao$^{**}$,  Carleton University }
\date{Revised, July 2012}

\begin{document}
 \maketitle

\begin{abstract}
In this paper, we propose a kernel method for exact tail
asymptotics of a random walk to neighborhoods in the quarter
plane. This is a two-dimensional method, which does not require a
determination of the unknown generating function(s). Instead, in
terms of the asymptotic analysis and a Tauberian-like theorem, we
show that the information about the location of the dominant
singularity or singularities and the detailed asymptotic property
at a dominant singularity is sufficient for the exact tail
asymptotic behaviour for the marginal distributions and also for
joint probabilities along a coordinate direction. We provide all
details, not only for a ``typical'' case, the case with a single
dominant singularity for an unknown generating function, but also
for all non-typical cases which have not been studied before. A
total of four types of exact tail asymptotics are found for the
typical case, which have been reported in the literature. We also
show that on the circle of convergence, an unknown generating
function could have two dominant singularities instead of one,
which can lead to a new periodic phenomena. Examples are
illustrated by using this kernel method. This paper can be
considered as a systematic summary and extension of existing
ideas, which also contains new and interesting research results.
\medskip

\end{abstract}

\section{Introduction}

Two-dimensional discrete random walks in the quarter plane are
classical models, that could be either probabilistic or
combinatorial. Studying these models is important and often
fundamental for both theoretical and applied purposes. For a
stable probabilistic model, it is of significant interest to study
its stationary probabilities. However, only for very limited
special cases, a closed-form solution is available for the
stationary probability distribution. This adds value to studying
tail asymptotic properties in stationary probabilities, since
performance bounds and approximations can often be developed from
the tail asymptotic property. The focus of this paper is to
characterize exact tail asymptotics. Specifically, we propose a
kernel method to systematically study the exact tail behaviour for
the stationary probability distribution of the random walk in the
quarter plane.

The kernel method proposed here is an extension of the classical
kernel method, first introduced by Knuth~\cite{Knuth:69} and later
developed as the kernel method by Banderier \textit{et
al.}~\cite{B-BM-D-F-G-GB:02}. The standard kernel method deals
with the case of a functional equation of the fundamental form
$K(x, y)F(x, y) = A(x, y)G(x) + B(x, y)$, where $F(x, y)$ and
$G(x)$ are unknown functions. The key idea in the kernel method is
to find a branch $y=y_0(x)$, such that, at $(x,y_0(x))$, the
kernel function is zero, or $K(x, y_0(x))=0$. When analytically
substituting this branch into the right hand side of the
fundamental form, we then have $G(x)=-B(x, y_0(x))/A(x, y_0(x))$,
and hence,
\[
    F(x, y) = \frac{-A(x, y)B(x,y_0(x))/A(x, y_0(x)) + B(x, y)}{K(x, y)}.
\]
However, applying the above idea to the fundamental form of a
two-dimensional random walk does not immediately lead to a
determination of the generating function $P(x,y)$. Instead, it
provides a relationship between two unknown generating functions
$\pi_1(x)$ and $\pi_2(y)$, referred to as the generating functions
for the boundary probabilities. This is the key challenge in the
analysis of using a kernel method. Therefore, a good understanding
on the interlace of these two functions is crucial.

Following the early research by Malyshev~\cite{Malyshev:72,
Malyshev:73}, the algebraic method targeting on expressing the
unknwon generating functions was further systematically updated in
Fayolle, Iasnogorodski and Malyshev~\cite{FIM:1999} based on the
study of the kernel equation. The authors indicated in their book
that: ``Even if asymptotic problems were not mentioned in this
book, they have many applications and are mostly interesting for
higher dimensions.'' The proposed kernel method in this paper is a
continuation of the study in \cite{FIM:1999}. Research on tail
asymptotics for various models following the method (determination
of the unknown generating function(s) first) of \cite{FIM:1999} or
other closely related methods can be found in
 Flatto and McKean~\cite{Flatto-McKeqn:77},
 Fayolle and Iasnogorodski~\cite{FI:1979},
 Fayolle, King and Mitrani~\cite{FKM:82},
 Cohen and Boxma~\cite{cohen-boxma:83},
 Flatto and Hahn~\cite{Flatto-Hahn:84},
 Flatto~\cite{Flatto:85},
 Wright~\cite{Wright:92},
 Kurkova and Suhov~\cite{Kurkova-Suhov:03},
 Bousquet-Melou~\cite{Bousquet-Melou:05},
 Morrison~\cite{Morrison:07},
 Li and Zhao~\cite{Li-Zhao:09, Li-Zhao:10},
 Guillemin and Leeuwaarden~\cite{Guillemin-Leeuwaarden:09},
 and Li, Tavakoli and Zhao~\cite{Li-Tavakoli-Zhao:11}.

Different from the work mentioned above, which requires
characterizing or expressing the unknown generating function, such
as a closed-form solution or an integral expression through
boundary value problems, the proposed kernel method only requires
the information about the dominant singularities of the unknown
function, including the location and detailed asymptotic property
at the dominant singularities. Because of this, the method makes
it possible to systematically deal with all random walks instead
of a model based treatment. In a recent research, Li and
Zhao~\cite{Li-Zhao:10b} applied this method to a specific model,
and Li, Tavakoli and Zhao~\cite{Li-Tavakoli-Zhao:11} to the
singular random walks. For exact tail asymptotics without a
determination of the unknown generating function(s) or Laplace
transformation function(s), different methods were used in the
following studies:
 Abate and Whitt~\cite{Abate-Whitt:97},
 Lieshout and Mandjes~\cite{LM:2008},
 Miyazawa and Rolski~\cite{Miyazawa-Rolski:09},
 Dai and Miyazawa~\cite{Dai-Miyazawa:10}.

Other methods for studying two-dimensional problems, including
exact tail asymptotics, also exist, for example, based on large
deviations, on properties of the Markov additive process (including matrix-analytic methods), or on
asymptotic properties of the Green functions. References include
    Borovkov and Mogul'skii~\cite{Borovkov-M:01},
    McDonald~\cite{McDonald:99},
    Foley and McDonald~\cite{Foley-McDonald:01,Foley-McDonald:05a,Foley-McDonald:05b},
    Khanchi~\cite{Khanchi:08,Khanchi:09},
    Adan, Foley and McDonald~\cite{Adan-Foley-McDonald:09},
    Raschel~\cite{Raschel:10},
    Miyazawa~\cite{Miyazawa:07,Miyazawa:09,Miyazawa:08},
    Kobayashi and Miyazawa~\cite{Kobayashi-Miyazawa:2011},
    Takahashi, Fujimoto and Makimoto~\cite{TFM:01},
    Haque~\cite{Haque:03},
    Miyazawa~\cite{Miyazawa:04},
    Miyazawa and Zhao~\cite{Miyazawa-Zhao:04},
    Kroese, Scheinhardt and Taylor~\cite{KST:04},
    Haque, Liu and Zhao~\cite{Haque-Liu-Zhao:05},
    Li and Zhao~\cite{Li-Zhao:05},
    Motyer and Taylor~\cite{Motyer-Taylor:06},
    Li, Miyazawa and Zhao~\cite{Li-Miyazawa-Zhao:07},
    He, Li and Zhao~\cite{He-Li-Zhao:08},
    Liu, Miyazawa and Zhao~\cite{Liu-Miyazawa-Zhao:08},
    Tang and Zhao~\cite{Tang-Zhao:08},
    Kobayashi, Miyazawa and Zhao~\cite{Kobayashi-Miyazawa-Zhao:10}, among others.
For more references, people may refer to a recent survey on tail
asymptotics of multi-dimensional reflecting processes for queueing
networks by Miyazawa~\cite{Miyazawa:11}.

The main focus of this paper is to propose a kernel method for
exact tail asymptotics of random walks in the quarter plane
following the ideal in \cite{FIM:1999}, based on which a complete
description of the exact tail asymptotics for stationary
probabilities of a non-singular genus 1 random walk is obtained.
 We claim that the unknown generating function $\pi_1(x)$, or equivalently, $\pi_2(y)$, has either one or two dominant singularities.
For the case of either one dominant singularity, or two dominant singularities with different asymptotic properties,
 a total of four types of exact tail asymptotics
exists: (1) exact geometric decay; (2) a geometric decay
multiplied by a factor of $n^{-1/2}$; (3) a geometric decay
multiplied by a factor of $n^{-3/2}$; and (4) a geometric decay
multiplied by a factor of $n$. These results are essentially not
new (for examples see references \cite{Borovkov-M:01,
Foley-McDonald:01, Foley-McDonald:05b, Miyazawa:09, Khanchi:09})
except that the fourth type is missing from previous studies for
the discrete random walk, but was reported for the continuous
random walk in \cite{Dai-Miyazawa:10}. For the case of two
dominant singularities with the same asymptotic property, a new
periodic phenomena in the tail asymptotic property is discovered,
which has not been reported in previous literature. For the tail
asymptotic behaviour of the non-boundary joint probabilities along
a coordinate direction, a new method based on recursive
relationships of probability generating functions will be applied,
which is an extension of the idea used in \cite{Li-Zhao:10b}.

For an unknown generating function of probabilities, a
Tauberian-like theorem is used as a bridge to link the asymptotic
property of the function at its dominant singularities to the tail
asymptotic property of its coefficient, or in our case, stationary
probabilities. This theorem does not require the monitonicity in
the probabilities, which is required by a standard Tauberian
theorem and cannot be verified in general,  or Heaviside
operational calculus, which is usually very difficult to be
rigorous. However, the price paid for applying the Tauberian-like
theorem requires more in analyticity of the function and detailed
information about all dominant singularities, or singularities on
the circle of convergence. Therefore we need to provide
information about how many singularities exist on the circle of
convergence and their detailed properties, such as the nature of
the singularity and the multiplicity in the case of the pole, for
the random walk. It is not always true that only one singularity
exists on the circle of convergence. Technical details are needed
to address these issues.

The kernel method immediately leads to exact tail asymptotics in
the boundary probabilities, in both directions, based on which
exact tail asymptotics in a marginal distribution will become
clear. However, it does not directly lead to exact tail asymptotic
properties for the joint probabilities along a coordinate
direction, except for the boundary probabilities as mentioned
above. Therefore, further efforts are required. In this paper, we
propose a method, based on difference equations of the unknown
generating functions, to do the asymptotic analysis, which
successfully overcomes the hurdle for exact tail asymptotics for
joint probabilities.

The rest of the paper is organized into eight sections. In
section~2, after the model description, the so-called fundamental
form for the random walk in the quarter plane is provided,
together with a stability condition. Section~3 contains necessary
properties for the two branches (or an algebraic function) defined
by the kernel equation and for the branch points of the branches.
These properties are either directly from \cite{FIM:1999} or
further refinements. Section~4 consists of six subsections for the
purpose of characterizing the asymptotic properties of the unknown
generating functions $\pi_1(x)$ and $\pi_2(y)$ at their dominant
singularities. Specifically, two Tauberian-like theorems are
introduced in subsection~1; the interlace between the two unknown
generating functions is discussed in subsection~2, which plays a
key role in the proposed kernel method; detailed properties for
singularities of the unknown generating functions are obtained in
subsections~3--5, which finally lead to the main theorem
(Theorem~\ref{theorem3.1}) in this section provided in the last
subsection. In Section~5, asymptotic analysis for the boundary
generating functions is carried out, which directly leads to the
tail asymptotics for the boundary probabilities in terms of the
Tauberian-like theorem. In Section~6, based on the asymptotic
results obtained for the generating function of boundary
probabilities in the previous section, and the fundamental form,
exact tail asymptotic properties for the two marginal
distributions are provided.
    Exact tail asymptotic
properties for joint probabilities along a coordinate direction is
addressed in Section~7, which is not a direct result from the
kernel method. Instead, we propose a difference equation method to
carry out an asymptotic analysis of a sequence of unknown
generating functions. The last section contributes to some
concluding remarks and two examples by applying the kernel method.

\section{Description of the Random Walk}

The random walk in the quarter plane used in this paper to
demonstrate the kernel method is a reflected random walk or a
Markov chain with the state space $\mathbb{Z}_+^2=\{(m,n); m, n
\text{ are non-negative integers }\}$. To describe this process,
we divide the whole quadrant $\mathbb{Z}_+^2$ into four regions:
the interior $S_+ = \{(m,n); m, n =1, 2, \ldots \}$, horizontal
boundary $S_1 = \{(m,0); m =1, 2, \ldots \}$, vertical boundary
$S_2 = \{(0,n); n =1, 2, \ldots \}$, and the origin $S_0=
\{(0,0)\}$, or $\mathbb{Z}_+^2= S_+ \cup S_1 \cup S_2 \cup S_0$.
In each of these regions, the transition is homogeneous.
Specifically, let $X_+$, $X_1$, $X_2$ and $X_0$ be random
variables having the distributions, respectively, $p_{i,j}$ with
$i, j =0, \pm 1$; $p^{(1)}_{i,j}$ with $i =0, \pm 1$ and $j =0,
1$; $p^{(2)}_{i,j}$ with $i =0, 1$ and $j =0, \pm 1$; and
$p^{(0)}_{i,j}$ with $i, j =0, 1$. Then, the transition
probabilities of the random walk (Markov chain)
$L_t=(L_1(t),L_2(t))$ are given by
\begin{align*}
    P(L_{t+1}= (m_2,n_2) &| L_t  = (m_1,n_1) \} = \\
    & \left \{ \begin{array}{ll}
P(X_+ = (m_2-m_1,n_2-n_1)), & \text{if } (m_2,n_2) \in S, (m_1,n_1) \in S_+, \\
P(X_k = (m_2-m_1,n_2-n_1)), & \text{if } (m_2,n_2) \in S,(m_1,n_1)  \in S_k \text{ with } k=0, 1, 2.
\end{array} \right.
\end{align*}

\subsection{Ergodicity conditions}

A stability (ergodic) condition can be found in Theorem~3.3.1 of Fayolle,
Iasnogorodski and Malyshev~\cite{FIM:1999}, which has been amended
by Kobayashi and Miyazawa as Lemma~2.1 in
\cite{Kobayashi-Miyazawa:2011}.
This condition is stated in terms of the drift vectors defined by
\begin{eqnarray*}
    M &=&(M_{x,} M_{y})= \biggl ( \sum_i i \Bigl  ( \sum_j p_{i,j} \Bigr ), \sum_j j  \Bigl ( \sum_i p_{i,j} \Bigr ) \biggr ), \\
    M^{(1)} &=& (M_{x,}^{(1)} M_{y}^{(1)})=\biggl ( \sum_i i \Bigl ( \sum_j p^{(1)}_{i,j} \Bigr ), \sum_j j \Bigl ( \sum_i p^{(1)}_{i,j} \Bigr )\biggr ), \\
    M^{(2)} &=& (M_{x,}^{(2)} M_{y}^{(2)})=\biggl ( \sum_i i \Bigl ( \sum_j p^{(2)}_{i,j} \Bigr ), \sum_j j \Bigl ( \sum_i p^{(2)}_{i,j} \Bigr )\biggr ).
\end{eqnarray*}

%Theorem 1.2.1
\begin{theorem}[Theorem~3.3.1 in \cite{FIM:1999} and Lemma~2.1 in
\cite{Kobayashi-Miyazawa:2011}] \label{theoremergodicity}
 When $M\neq 0$, the random walk is
ergodic if and only if one of the following three conditions
holds:

\textbf{1.} $M_{x}<0$, $M_{y}<0$,
$M_{x}M_{y}^{(1)}-M_{y}M_{x}^{(1)}<0$ and
$M_{y}M_{x}^{(2)}-M_{x}M_{y}^{(2)}<0$;

\textbf{2.} $M_{x}<0$, $M_{y}\geq 0$,
$M_{y}M_{x}^{(2)}-M_{x}M_{y}^{(2)}<0$ and $M_{x}^{(1)}<0$ if
$M_{y}^{(1)}=0$;

\textbf{3.} $M_{x}\geq 0$, $M_{y}<0$,
$M_{x}M_{y}^{(1)}-M_{y}M_{x}^{(1)}<0$ and $M_{y}^{(2)}<0$ if
$M_{x}^{(2)}=0$.
\end{theorem}

Throughout the paper, we make the following assumption, unless
otherwise specified:

\begin{assumption} \label{assumption-1}
The random walk $L_t$ is irreducible, positive recurrent and aperiodic.
\end{assumption}
Under Assumption~\ref{assumption-1}, let $\pi_{m,n}$ be the unique stationary probability distribution
of the random walk.

\begin{remark}
It should be noted that for a stable random walk, the condition $M
\neq 0$ is equivalent to that both sequences $\{\pi_{m,0}\}$ and
$\{\pi_{0,n}\}$ are light-tailed (for example, see Lemma~3.3 of
\cite{Kobayashi-Miyazawa:2011}), which is not our focus of this
paper. Therefore, Theorem~\ref{theoremergodicity} provides a
necessary and sufficient stability condition for the light-tailed
case.
\end{remark}

\subsection{Fundamental Form}

Define the following generating functions of the probability sequences for the interior states, horizontal boundary states and
vertical boundary states, respectively,
\begin{eqnarray*}
    \pi(x,y) &=& \sum_{m=1}^{\infty} \sum_{n=1}^{\infty} \pi_{m,n}x^{m-1}y^{n-1}, \\
    \pi_1(x) &=&\sum_{m=1}^{\infty} \pi_{m,0}x^{m-1}, \\
    \pi_2(y) &=&\sum_{n=1}^{\infty} \pi_{0,n}y^{n-1}.
\end{eqnarray*}

The so-called fundamental form of the random walk provides a functional equation relating the three unknown generating functions $\pi(x,y)$, $\pi_1(x)$ and
$\pi_2(y)$. To state the fundamental form, we define

\begin{eqnarray*}
    h(x,y) &=&xy\left( \sum_{i=-1}^{1}\sum_{j=-1}^{1}p_{i,j}x^{i}y^{j}-1\right) \\
    &=&a(x)y^{2}+b(x)y+c(x)=\tilde{a}(y)x^{2}+\tilde{b}(y)x+\tilde{c}(y), \\
    h_{1}(x,y) &=&x\left(\sum_{i=-1}^{1}\sum_{j=0}^{1}p_{i,j}^{(1)}x^{i}y^{j}-1\right) \\
    &=&a_{1}(x)y+b_{1}(x)=\widetilde{a}_{1}(y)x^{2}+\widetilde{b}_{1}(y)x+\widetilde{c}_{1}(y), \\
    h_{2}(x,y) &=&y\left(\sum_{i=0}^{1}\sum_{j=-1}^{1}p_{i,j}^{(2)}x^{i}y^{j}-1\right) \\
    &=&\widetilde{a}_{2}(y)x+\widetilde{b}_{2}(y)=a_{2}(x)y^{2}+b_{2}(x)y+c_{2}(x), \\
    h_{0}(x,y) &=&\left(\sum_{i=0}^{1}\sum_{j=0}^{1}p_{i,j}^{(0)}x^{i}y^{j}-1\right) \\
    &=&a_{0}(x)y+b_{0}(x)=\widetilde{a}_{0}(y)x+\widetilde{b}_{0}(y),
\end{eqnarray*}
where
\begin{eqnarray*}
a(x) &=&p_{-1,1}+p_{0,1}x+p_{1,1}x^{2}, \\
b(x) &=&p_{-1,0}-(1-p_{0,0})x+p_{1,0}x^{2}, \\
c(x) &=&p_{-1,-1}+p_{0,-1}x+p_{1,-1}x^{2}, \\
\tilde{a}(y) &=&p_{1,-1}+p_{1,0}y+p_{1,1}y^{2}, \\
\tilde{b}(y) &=&p_{0,-1}-(1-p_{0,0})y+p_{0,1}y^{2}, \\
\tilde{c}(y) &=&p_{-1,-1}+p_{-1,0}y+p_{-1,1}y^{2},
\end{eqnarray*}
\begin{eqnarray*}
    a_{1}(x) &=&p_{-1,1}^{(1)}+p_{0,1}^{(1)}x+p_{1,1}^{(1)}x^{2},
b_{1}(x)=p_{-1,0}^{(1)}-\left( 1-p_{0,0}^{(1)}\right) x+p_{1,0}^{(1)}x^{2},\\
    \widetilde{a}_{1}(y) &=&p_{1,0}^{(1)}+p_{1,1}^{(1)}y,
\widetilde{b}_{1}(y)=p_{0,0}^{(1)}-1+p_{0,1}^{(1)}y, \widetilde{c}_{1}(y)=p_{-1,0}^{(1)}+p_{-1,1}^{(1)}y \\
    a_{2}(x) &=&p_{0,1}^{(2)} + p_{1,1}^{(2)}x,
b_{2}(x)=p_{0,0}^{(2)}-1+p_{1,0}^{(2)}x,
c_{2}(x)=p_{0,-1}^{(2)}+p_{1,-1}^{(2)}x \\
    \widetilde{a}_{2}(y) &=&p_{1,-1}^{(2)}+p_{1,0}^{(2)}y+p_{1,1}^{(2)}y^{2},
\widetilde{b}_{2}(y)=p_{0,-1}^{(2)}-\left( 1-p_{0,0}^{(2)}\right) y+p_{0,1}^{(2)}y^{2}, \\
    a_{0}(x) &=&p_{0,1}^{(0)}+p_{1,1}^{(0)}x,b_{0}(x)=p_{1,0}^{(0)}x-\left( 1-p_{0,0}^{(0)}\right) , \\
\widetilde{a}_{0}(y) &=&p_{1,0}^{(0)}+p_{1,1}^{(0)}y,\widetilde{b}_{0}(y)=p_{0,1}^{(0)}y-\left( 1-p_{0,0}^{(0)}\right) .
\end{eqnarray*}

The basic equation of the generating function of the joint
distribution, or the fundamental form of the random walk, is given
by
\begin{equation} \label{eqn:fundamental}
    -h(x,y)\pi(x,y)=h_{1}(x,y)\pi_1(x)+h_{2}(x,y)\pi_2(y)+h_{0}(x,y)\pi_{0,0}.
\end{equation}
The reason for the above functional equation to be called
fundamental is largely due to the fact that through analysis of
this equation, the unknown generating functions can be determined
or expressed, for example, through algebraic methods and boundary
value problems as illustrated in Fayolle, Iasnogorodski and
Malyshev~\cite{FIM:1999}. The kernel method presented here also
starts with the fundamental form, but without expressing
generating functions first.

\begin{remark}
The generating function $\pi(x,y)$ is defined for $m,n>0$, excluding the boundary probabilities.
(\ref{eqn:fundamental}) was proved in (1.3.6) in \cite{FIM:1999}. Based on (\ref{eqn:fundamental}), one can also obtain a similar fundamental form using generating functions including boundary
probabilities: $\Pi(x,y)=\sum_{m=0}^\infty \sum_{n=0}^\infty \pi_{m,n} x^m y^n$, $\Pi_1(x)=\sum_{m=0}^\infty \pi_{m,0} x^m$ and
$\Pi_2(y)=\sum_{n=0}^\infty \pi_{0,n} y^n$.
\end{remark}

For the conclusion of this section, we can easily check the
following expressions, some of which will be needed in later
sections:
\begin{align}
    M_{y} &= a(1)-c(1) = \widetilde{a}^{\prime }(1) + \widetilde{b}^{\prime}(1) + \widetilde{c}^{\prime}(1), &
        M_{x} & = \widetilde{a}(1) - \widetilde{c}(1) = a^{\prime}(1)+b^{\prime}(1)+c^{\prime}(1),  \label{eqn:MxMy} \\
    M_{y}^{(1)} &= a_{1}(1) = \widetilde{a}^{\prime}_1(1) + \widetilde{b}^{\prime}_1(1) + \widetilde{c}^{\prime}_1(1), &
        M_{x}^{(1)} & = \widetilde{a}_1(1) - \widetilde{c}_1(1) = a_{1}^{\prime}(1) + b_{1}^{\prime}(1), \\
    M_{y}^{(2)} &= a_{2}(1) - c_{2}(1) = \widetilde{a}^{\prime}_2(1) + \widetilde{b}^{\prime}_2(1), &
        M_{x}^{(2)} &= \widetilde{a}_2(1) = a_{2}^{\prime}(1) + b_{2}^{\prime}(1) + c_{2}^{\prime}(1).
\end{align}

\section{Branch Points And Functions Defined by the Kernel Equation}

The property of the random walk relies on the property of the
kernel function $h$ and functions $h_1$ and $h_2$. The kernel
function plays a key role in the kernel method.

\begin{definition}
A random variable is called non-singular if the kernel function
$h(x,y)$, as a polynomial in the two variables $x$ and $y$, is
irreducible (equivalently, if $h=fg$ then either $f$ or $g$ is a
constant) and quadratic in both variables.
\end{definition}

Throughout the paper unless otherwise specified, we make the
second assumption below.
\begin{assumption} \label{assumption-2}
The random considered is non-singular.
\end{assumption}

The non-singular condition for a random walk is closely related to
the irreducibility of the marginal processes $L_1(t)$ and
$L_2(t)$, but they are not the same concept. A necessary and
sufficient condition for a random walk to be singular is given, in
terms of $p_{i,j}$, in Lemma~2.3.2 in \cite{FIM:1999}. Study on
tail asymptotics for a singular random walk is either easier or
similar to the non-singular case, which can be found in Li,
Tavakoli and Zhao~\cite{Li-Tavakoli-Zhao:11}.

The starting point of our analysis is the set of all pairs $(x,y)$
satisfy the kernel equation, or
\[
        B =\{(x,y) \in \mathbb{C}^2: h(x,y)=0\},
\]
where $\mathbb{C}$ is the set of all complex numbers. The kernel
function can be considered as a quadratic form in either $x$ or
$y$ with the coefficients being functions of $y$ or $x$,
respectively. Therefore, the kernel equation can be written as
\begin{equation} \label{eqn:algebraic-functions}
    a(x)y^{2}+b(x)y+c(x)=\widetilde{a}(y)x^{2}+\widetilde{b}(y)x+\widetilde{c}(y)=0.
\end{equation}

 For a fixed $x$, the two solutions to the kernel equation as a quadratic form in $y$ are given by
\begin{equation*}
    Y_{\pm }(x)=\frac{-b(x)\pm \sqrt{D_{1}(x)}}{2a(x)}
\end{equation*}
if $a(x) \neq 0$, where $D_{1}(x)=b^{2}(x)-4a(x)c(x)$. Notice that
non-singularity implies that $a(x) \not\equiv 0$ and, therefore,
only up to two values of $x$ could lead to $a(x)=0$ since $a(x)$
is a polynomial of degree up to 2.

 Similarly, for a fixed $y$, the two solutions to the kernel equation as a quadratic form in $x$ are given by
\begin{equation*}
    X_{\pm }(y)=\frac{-\widetilde{b}(y)\pm \sqrt{D_{2}(y)}}{2\widetilde{a}(y)},
\end{equation*}
where $D_{2}(y)=\widetilde{b}^{2}(y)-4\widetilde{a}(y)\widetilde{c}(y)$.

It is important to study the set $B$, or equivalently $Y_{\pm}(x)$
or $X_{\pm}(y)$, since for all $(x,y) \in B$ with $|\pi(x,y)|
<\infty$, the right hand side of the fundamental form is also
zero, which provides a relationship between the two unknown
generating functions $\pi_1$ and $\pi_2$. In the above,
$\sqrt{D_1(x)}$ is well-defined if $D_1(x)\geq 0$ and similarly
$\sqrt{D_2(y)}$ is well-defined if $D_2(y)\geq 0$. As a function
of a complex variable, the square root is a two-valued function.
To specify a branch, when $z$ is complex, $\sqrt{z}$ is defined
such that $\sqrt{1}=1$.

Let $z = D_1(x)$. Then, both $Y_{-}(x)$ and $Y_{+}(x)$ are
analytic as long as $z \notin (-\infty, 0]$ and $a(x) \neq 0$. For
these two functions, we start from a region, in which they are
analytic, and consider an analytic continuation of these two
functions. In this consideration, the key is the continuation of
$\sqrt{D_1(x)}$.

\begin{definition} A branch point of $Y_{\pm}(x)$ ($X_{\pm}(y)$) is a value of $x$ ($y$) such that $D_1(x)=0$ ($D_2(y)=0$).
\end{definition}

To discuss the branch points, notice that
the discriminant $D_1$ ($D_2$) is a polynomial of degree up to four. Since the two cases are symmetric, we discuss $D_1(x)$ in detail only.
Rewrite $D_1(x)$ as
\begin{equation*}
    D_{1}(x)=d_{4}x^{4}+d_{3}x^{3}+d_{2}x^{2}+d_{1}x+d_{0},
\end{equation*}
where
\begin{align*}
    d_{0} &= p_{-1,0}^{2}-4p_{-1,1}p_{-1,-1}, \\
    d_{1} &= 2p_{-1,0}(p_{0,0}-1)-4(p_{-1,1}p_{0,-1}+p_{0,1}p_{-1,-1}), \\
    d_{2} &= (p_{0,0}-1)^{2}+2p_{1,0}p_{-1,0}-4(p_{1,1}p_{-1,-1}+p_{1,-1}p_{-1,1}+p_{0,1}p_{0,-1}), \\
    d_{3} &= 2p_{1,0}(p_{0,0}-1)-4(p_{1,1}p_{0,-1}+p_{0,1}p_{1,-1}), \\
    d_{4} &= p_{1,0}^{2}-4p_{1,1}p_{1,-1}.
\end{align*}
It can be easily checked that $d_1 \leq 0$ and $d_3 \leq 0$.

When $D_1$ is a polynomial of degree 4 (or $d_4 \neq 0$), there are four
branch points, denoted by $x_i$ ($y_i$), $i=1, 2, 3, 4$. Without
loss of generality, we assume that $|x_1| \leq |x_2| \leq |x_3|
\leq |x_4|$. When the
degree of $D_1(x)$ is $d<4$, for convenience, we let
$x_{d+k}=\infty$ for integer $k>0$ such that
$d+k \leq 4$. For example, if $d=3$, then $x_4=\infty$.
This can be justified by the following: consider
the polynomial $\tilde{D}_1(\widetilde{x})=D_1(x)/x^4$ in
$\widetilde{x}$, where $\widetilde{x}=1/x$. Then,
$\widetilde{x}=0$ is a $d$-tuple zero of
$\widetilde{D}_1(\tilde{x})$, and therefore $x=\infty$ can be
viewed as a $d$-tuple zero of $D_1(x)$.

The following lemma characterizes the branch points of
$Y_{\pm}(x)$ for all non-singular random walks, including the
heave-tailed case, or the case of $M =0$.

%Lemma 1.1
\begin{lemma} \label{lemma1.1}
\textbf{1.} For a non-singular random walk with $M_{y}\neq 0$,
$Y_{\pm}(x)$ has two branch points $x_{1}$ and $x_{2}$ inside the
unit circle and another two branch points $x_{3}$ and $x_{4}$
outside the unit circle. All these branch points lie on the real
line. More specifically,
\begin{description}
\item[(1)] if $p_{1,0}>2\sqrt{p_{1,1}p_{1,-1}}$, then
$1<x_{3}<x_{4}<\infty$;

\item[(2)] if $p_{1,0}=2\sqrt{p_{1,1}p_{1,-1}}$, then
$1<x_{3}<x_{4}=\infty$;

\item[(3)] if $p_{1,0}<2\sqrt{p_{1,1}p_{1,-1}}$, then $1<x_{3}
\leq - x_{4} < \infty$, where the equality holds if and only if
$d_1=d_3=0$.
\end{description}
Similarly,
\begin{description}
\item[(4)] if $p_{-1,0}>2\sqrt{p_{-1,1}p_{-1,-1}}$, then
$0<x_{1}<x_{2}<1$;

\item[(5)] if $p_{-1,0}=2\sqrt{p_{-1,1}p_{-1,-1}}$, then $x_{1}=0$
and $0<x_{2}<1$;

\item[(6)] if $p_{-1,0}<2\sqrt{p_{-1,1}p_{-1,-1}}$, then $0<-x_1
\leq x_{2}<1$, where the equality holds if and only if
$d_1=d_3=0$.
\end{description}

\textbf{2.} For a non-singular random walk with $M_{y}=0$ (in this
case $M_{x}\neq 0$ since we are only considering the genus 1 case
in this paper), either $x_{2}=1$ if $M_{x}<0$; or $x_{3}=1$ if
$M_{x}>0$. In the latter case, the system is unstable.
\end{lemma}

\proof We only need to prove \textbf{3.} and \textbf{6.} since all
other proofs can be found in Fayolle, Iasnogorodski and
Malyshev~\cite{FIM:1999} (Lemma~2.3.8 and Lemma~2.3.9). We provide details for \textbf{3.}
since \textbf{6.} can be proved similarly. Suppose otherwise
$x_{3}>-x_{4}$. From $d_{1}\leq 0$ and $d_{3}\leq 0$, we obtain
$D_{1}(-x_{3})=-d_{3}x_{3}^{3}-d_{1}x_{3}>0$. On the other hand,
$D_{1}(-\infty )=-\infty$ since $d_{4}<0$, which implies that
$D_{1}(x)=0$ has a fifth root in $(-\infty ,x_{3})$, but this is
impossible. The contradiction shows that $x_{3}\leq -x_{4}$. It is
clear that the equality holds if and only if $d_{1}=$ $d_{3}=0$.
 \QED

\begin{remark}
Similar results hold for the branch points $y_i$, $i=1, 2, 3, 4$, of $X_{\pm}(y)$.
\end{remark}

\begin{definition}
$p_{i,j}$ ($p_{i,j}^{(k)}$) is called X-shaped if $p_{i,j}=0$
($p_{i,j}^{(k)}=0$) for all $i$ and $j$ such that $|i+j|=1$. A
random walk is called X-shaped if $p_{i,j}$ and also
$p_{i,j}^{(k)}$ for $k=1, 2$ are all X-shaped.
\end{definition}

Based on Lemma~\ref{lemma1.1}, we can prove the following result.
\begin{corollary} \label{corollary-X}
$x_3=-x_4$ if and only if $p_{i,j}$ is X-shaped.
\end{corollary}

Throughout the rest of the paper, we define $[x_{3},x_{4}]=[-\infty ,x_{4}]\cup [x_{3},\infty ]$  when $x_{4}<-1$. Similarly,
$[y_{3},y_{4}]=[-\infty ,y_{4}]\cup [y_{3},\infty ]$] when $y_4<-1$. We define the following cut planes:
\begin{eqnarray*}
    \widetilde{\mathbb{C}}_{x} &=&\mathbb{C}_{x} \setminus [x_{3},x_{4}], \\
    \widetilde{\mathbb{C}}_{y} &=&\mathbb{C}_{y} \setminus [y_{3},y_{4}], \\
    \widetilde{\widetilde{\mathbb{C}}}_{x} &=&\mathbb{C}_{x} \setminus [x_{3},x_{4}] \cup [x_{1},x_{2}], \\
    \widetilde{\widetilde{\mathbb{C}}}_{y} &=&\mathbb{C}_{y} \setminus [y_{3},y_{4}]\cup [y_{1},y_{2}],
\end{eqnarray*}
where $\mathbb{C}_{x}$ and $\mathbb{C}_{y}$ are the complex planes for $x$ and $y$, respectively.

We now define two complex functions on the cut plane $\widetilde{\widetilde{\mathbb{C}}}_{x}$ based on $Y_{\pm}(x)$:
\begin{equation}
    Y_{0}(x) = \left \{ \begin{array}{ll}
    Y_{-}(x), &  \text{if } |Y_{-}(x)| \leq |Y_{+}(x)|, \\
    Y_{+}(x), &  \text{if } |Y_{-}(x)|>|Y_{+}(x)|; \end{array}
    \right.
\end{equation}
and
\begin{equation}
    Y_{1}(x) = \left \{ \begin{array}{ll}
    Y_{+}(x), &  \text{if } |Y_{-}(x)| \leq |Y_{+}(x)|, \\
    Y_{-}(x), &  \text{if } |Y_{-}(x)|>|Y_{+}(x)|. \end{array}
    \right.
\end{equation}
Obviously,  $Y_0$ is the function of $Y_{-}$ and $Y_{+}$ with the smaller modulus and $Y_{+}$ is the function with the larger modulus.

Functions $X_0(y)$ and $X_1(y)$ are defined on the cut plane $\widetilde{\mathbb{C}}_{y}$ in the same manner.

\begin{remark}
A branch point of $Y_{\pm}(x)$ ($X_{\pm}(y)$) is also referred to
as a branch point of $Y_0(x)$ and $Y_1(x)$ ($X_0(y)$ and
$X_1(y)$).
\end{remark}

\begin{remark} It is not always the case that $Y_0$ is a continuation of $Y_-$ and
$Y_1$ a continuation of $Y_+$. However, for $x \in \widetilde{\widetilde{\mathbb{C}}}_{x}$ with $a(x) \neq 0$, $Y_0(x)$ and $Y_1(x)$ are still the two
zeros of the kernel function $h(x,y)$. Parallel comments can be made on $X_0$ and $X_1$.
\end{remark}

A list of basic properties of $Y_0$ and $Y_1$ ($X_0$ and $X_1$) is provided in the following lemma.
\begin{lemma} \label{lemma1.1-b}
\textbf{1.} For $|x|=1$, $|Y_{0}(x)|\leq 1$ and $|Y_{1}(x)|\geq
1$, with equality only possibly for $x=\pm 1$. For $x=1$, we
have
\begin{align*}
    Y_{0}(1) &= \min \left( 1,\frac{\sum p_{i,-1}}{\sum p_{i,1}}=\frac{c(1)}{a(1)} \right ), \\
    Y_{1}(1) &= \max \left( 1,\frac{\sum p_{i,-1}}{\sum p_{i,1}}=\frac{c(1)}{a(1)} \right );
\end{align*}
for $x=-1$, the equality holds only if $p_{i,j}$ is X-shaped, for which we have
\begin{align*}
    Y_{0}(-1) &= -\min \left( 1,\frac{\sum p_{i,-1}}{\sum p_{i,1}}=\frac{c(1)}{a(1)} \right ), \\
    Y_{1}(-1) &= -\max \left( 1,\frac{\sum p_{i,-1}}{\sum p_{i,1}}=\frac{c(1)}{a(1)} \right ).
\end{align*}

\textbf{2.} The functions $Y_{i}(x)$, $i=0,1$, are meromorphic in
the cut plane $\widetilde{\widetilde{\mathbb{C}}}_{x}$. In
addition,
\begin{description}
\item[(a)] $Y_{0}(x)$ has two zeros and no poles. Hence $Y_{0}(x)$
is analytic in $\widetilde{\widetilde{\mathbb{C}}}_{x}$;

\item[(b)] $Y_{1}(x)$ has two poles and no zeros.

\item[(c)] $|Y_{0}(x)|\leq |Y_{1}(x)|$, in the whole cut complex
plane $\widetilde{\widetilde{\mathbb{C}}}_{x}$, and equality takes
place only on the cuts.
\end{description}

\textbf{3.} The function $Y_{0}(x)$ can become infinite at a point $x$ if
and only if,
\begin{description}
\item[(a)] $p_{11}=p_{10}=0$, in this case, $x=x_{4}=\infty$; or

\item[(b)] $p_{-11}=p_{-10}=0$, in this case, $x=x_{1}=0$.
\end{description}

Parallel conclusions can be made for functions $X_0(y)$ and $X_1(y)$.
\end{lemma}

\proof This lemma contains results in Lemma~2.3.4 and Theorem~5.3.3 in \cite{FIM:1999}. First, according to (ii) of Theorem~5.3.3 in
\cite{FIM:1999}, the functions $Y_0$ and $Y_1$ defined in this paper coincide the functions $Y_0$ and $Y_1$ in \cite{FIM:1999} due to
the uniqueness of the continuity. Then, all results in 1 come from Lemma~2.3.4 and Lemma~5.3.1 in \cite{FIM:1999} except for the expressions for $Y_0(-1)$ and $Y_1(-1)$, which
 can be obtained in the same fashion as for $Y_0(1)$ and $Y_1(1)$; results in
2 are given in (ii) of Theorem~5.3.3 in \cite{FIM:1999}; and the conclusion in 3 is the same as in (iii) Theorem~5.3.3 in \cite{FIM:1999}. \QED

\begin{remark} All the above properties can be directly obtained through elementary analysis of the square root function. \end{remark}

Throughout the rest of the paper, unless otherwise specified, we
make the following assumption:
\begin{assumption} \label{assumption-3}
All branch points $x_i$ and $y_i$, $i=1,2,3,4$, are distinct.
\end{assumption}
A random walk satisfying Assumption~\ref{assumption-3} is called a genus 1 random walk.

\begin{remark}
This assumption is equivalent to the assumption that the Riemann surface defined by the kernel equation has genus 1. The Riemann surface for the random walk is either genus 1 or genus 0.
A necessary and sufficient condition for the random walk in the quarter plane
to be genus 1 is given in Lemma~2.3.10 in \cite{FIM:1999}. Most of queueing application models are the case of genus 1. The genus 0 case can be analyzed similarly except for the heavy-tailed case,
the case where $M=0$. In general, analysis of the genus 0 case (except for the case of $M=0$) could be less challenging since expressions for the unknown generating functions $\pi_1(x)$ and $\pi_2(y)$ are either explicit or
less complex than for the genus 1 case, which can immediately lead to an analytic continuation of these unknown generating functions. Chapter~6 of \cite{FIM:1999} is devoted to the genus 0 case.
\end{remark}

\begin{corollary}
For a non-singular genus 1 random walk, if $p_{i,j}$ is X-shaped, then
all $p_{1,1}$, $p_{1,-1}$, $p_{-1,1}$ and $p_{-1,-1}$ are
positive.
\end{corollary}

 \proof If only one of $p_{1,1}$, $p_{1,-1}$, $p_{-1,1}$ and
$p_{-1,-1}$ is zero, then the random walk is non-singular having
genus 0 (Lemma~2.3.10 in \cite{FIM:1999}) and if at least two of
them are zero, then the random walk is singular (Lemma~2.3.2 in
\cite{FIM:1999}). \QED

\begin{corollary}
For a stable random walk with $M \neq 0$,

\textbf{1.} If $p^{(1)}_{i,j}$ is X-shaped, then
$p^{(1)}_{1,1}$ and $p^{(1)}_{-1,1}$ cannot be both zero; and

\textbf{2.} If $p^{(2)}_{i,j}$ is X-shaped, then
$p^{(2)}_{1,1}$ and $p^{(2)}_{-1,-1}$ cannot be both zero.
\end{corollary}

 \proof Otherwise, $p^{(k)}_{0,0}=1$, $k=1$ or 2, with which the random walk cannot be stable.
 \QED

For our purpose, more results about functions $Y_0$ and $Y_1$
($X_0$ and $X_1$) are needed. Once again, we consider $Y_0$ and
$Y_1$. $X_0$ and $X_1$ can be considered in the same way. Recall
that $Y_k$ ($k=1,2$) are defined on the cut plane $\mathbb{C}_{x}
\setminus [x_{3},x_{4}] \cup [x_{1},x_{2}]$, where two slits
$[x_1,x_2]$ and $[x_3,x_4]$ are removed from the complex plane
such that the functions $Y_k$ can stay always in one branch. Take
the slit $[x_1,x_2]$ as an example. For $x' \in [x_1,x_2]$, the
limit of $Y_k(x)$ when $x$ approaches to $x'$ from above the real
axis is different from the limit as $x$ approaches to $x'$ from
below the real axis. Let $x'' \in [x_1,x_2]$ be another point
satisfying $x'' >x'$. By
$Y_{0}[\underleftarrow{\overrightarrow{x'x''}}]$, we denote the
image contour, which is the limit of $Y_0(x)$ from above the real
axis when $x$ traverses from $x'$ to $x''$ and from below the real
axis when $x$ continues to traverse back from $x''$ to $x'$. For
convenience, we say that
$Y_{0}[\underleftarrow{\overrightarrow{x'x''}}]$ is the image of
the contour $\underleftarrow{\overrightarrow{x'x''}}$, traversed
from $x'$ to $x''$ along the upper edge of the slit $[x',x'']$ and
then back to $x'$ along the lower edge of the slit. In this way,
we can define the following image contours:
\begin{align}
    \mathscr{L} & = Y_{0}[\underleftarrow{\overrightarrow{x_{1}x_{2}}}], \;\;\;
    {\mathscr{L}}_{ext}  = Y_{0}[\underleftarrow{\overrightarrow{x_{3}x_{4}}}]; \\
    \mathscr{M} & = X_{0}[\underleftarrow{\overrightarrow{y_{1}y_{2}}}], \;\;\;
    {\mathscr{M}}_{ext} = X_{0}[\underleftarrow{\overrightarrow{y_{3}y_{4}}}],
\end{align}
respectively.
Furthermore, for an arbitrary simple closed curve $\mathscr{U}$,
by $G_{\mathscr{U}}$ we denote the interior domain bounded by
$\mathscr{U}$ and  by $G_{\mathscr{U}}^{c}$ the exterior domain.

The properties of the above image contours provided in the following lemma are important for the interlace between the two unknown functions $\pi_(x)$ and
$\pi_2(y)$ discussed in the next section. To state the lemma,
define the following determinant:
\begin{eqnarray}
    \Delta &=&\left\vert
\begin{tabular}{lll}
$p_{11}$ & $p_{10}$ & $p_{1,-1}$ \\
$p_{01}$ & $p_{00}$ & $p_{0,-1}$ \\
$p_{-1,1}$ & $p_{-1,0}$ & $p_{-1,-1}$%
\end{tabular}%
\right\vert  \notag .
\end{eqnarray}

%Lemma 1.4 (Theorem 5.3.3 (i) and Corollary 5.3.5 of Fayolle's book)

\begin{lemma}\label{lemma1.4}
 For non-singular genus 1
random walk without branch points on the unit circle, we have the
following properties:

\textbf{1.} The curve $\mathscr{M}$ and $\mathscr{M}_{ext }$
 are simple, closed and symmetrical about the
real axis in $\mathbb{C}_{x}$ plane. Moreover,

\textbf{(a)} If $\Delta >0$, then
\begin{equation*}
    \lbrack x_{1},x_{2}]\subset G_{\mathscr{M}}\subset G_{\mathscr{M}_{ext}}
\;\text{ and }\; [x_{3},x_{4}]\subset G_{\mathscr{M}_{ext}}^{c};
\end{equation*}

\textbf{(b)} If $\Delta <0$, then
\begin{equation*}
    \lbrack x_{1},x_{2}]\subset G_{\mathscr{M}_{ext}}\subset G_{\mathscr{M}}
\;\text{ and }\; [x_{3},x_{4}]\subset G_{\mathscr{M}}^{c};
\end{equation*}

\textbf{(c)} If $\Delta =0$, then
\begin{equation*}
    \lbrack x_{1},x_{2}]\subset G_{\mathscr{M}_{ext}}=G_{\mathscr{M}}
\;\text{ and }\; [x_{3},x_{4}]\subset G_{\mathscr{M}}^{c}.
\end{equation*}
Entirely symmetric results hold for $\mathscr{L}$ and
$\mathscr{L}_{ext}$.
\medskip

\textbf{2.} The branches $X_i$ and $Y_i$ have the following
properties:

\textbf{(a)} Both $X_0(y)$ and $Y_0(x)$ are conformal mappings:
$G_{\mathscr{M}}-[x_{1},x_{2}]\overset{Y_{0}(x)}{\underset{X_{0}(y)}{\rightleftarrows}}G_{\mathscr{L}}-[y_{1},y_{2}]$;

\textbf{(b)} $X_{0}(y)\in G_{\mathscr{M}} \cup
G_{\mathscr{M}_{ext}}$ and $X_{1}(y) \in G_{\mathscr{M}}^{c} \cup
G_{\mathscr{M}_{ext}}^{c}$. Symmetrically, $Y_{0}(x) \in
G_{\mathscr{L}} \cup G_{\mathscr{L}_{ext}}$ and  $Y_{1}(x)\in
G_{\mathscr{L}}^{c} \cup G_{\mathscr{L}_{ext}}^{c}$;

\textbf{(c)} If $G_{\mathscr{M}} \subset G_{\mathscr{M}_{ext}}$,
then
\begin{eqnarray*}
X_{0}\circ Y_{0}(t) &=&t, \text{ if } t \in G_{\mathscr{M}}, \\
X_{0}\circ Y_{0}(t) &\neq &t, \text{ if }t\in G_{\mathscr{M}}^{c}
    \text{ and } X_{0}\circ
    Y_{0}(G_{\mathscr{M}}^{c})=G_{\mathscr{M}}.
%X_{0}\circ Y_{1}(t) &=&t, \text{ if }t\in G_{\mathscr{M}_{ext}}, \\
%X_{0}\circ Y_{1}(t) &\neq &t, \text{ if } t \in
%    G_{\mathscr{M}_{ext}}^{c}\text{ and } X_{0}\circ Y_{1}(G_{\mathscr{M}_{ext}}^{c})=G_{\mathscr{M}_{ext}}, \\
%X_{1}\circ Y_{0}(t) &=&t, \text{ if }t\in G_{\mathscr{M}}^{c}, \\
%    X_{1}\circ Y_{0}(t) &\neq &t, \text{ if }t\in G_{\mathscr{M}}
%    \text{ and }X_{1}\circ Y_{0}(G_{\mathscr{M}})=G_{\mathscr{M}}^{c}, \\
%X_{1}\circ Y_{1}(t) &=&t, \text{ if }t\in G_{\mathscr{M}_{ext}}^{c}, \\
%    X_{1}\circ Y_{1}(t) &\neq &t, \text{ if }t\in G_{\mathscr{M}_{ext}}
%    \text{ and } X_{1}\circ Y_{1}(G_{\mathscr{M}_{ext}})=G_{\mathscr{M}_{ext}}^{c}.
\end{eqnarray*}

Symmetrically, if $G_{\mathscr{L}}\subset G_{\mathscr{L}_{ext}}$,
then
\begin{eqnarray*}
Y_{0}\circ X_{0}(t) &=&t\text{ if }t\in G_\mathscr{L}, \\
Y_{0}\circ X_{0}(t) &\neq &t\text{ if }t\in
G_\mathscr{L}^{c}\text{ and }Y_{0}\circ
X_{0}(G_\mathscr{L}^{c})=G_\mathscr{L}.
%Y_{0}\circ X_{1}(t) &=&t\text{ if }t\in G_\mathscr{L_{ext}}, \\
%Y_{0}\circ X_{1}(t) &\neq &t\text{ if }t\in
%G_\mathscr{L_{ext}}^{c}\text{ and
%} Y_{0}\circ X_{1}(G_\mathscr{L_{ext}}^{c})=G_\mathscr{L_{ext}}, \\
%Y_{1}\circ X_{0}(t) &=&t\text{ if }t\in G_\mathscr{L}^{c}, \\
%Y_{1}\circ X_{0}(t) &\neq &t\text{ if }t\in G_\mathscr{L}\text{
%and
%}Y_{1}\circ X_{0}(G_\mathscr{L})=G_\mathscr{L}^{c}, \\
%Y_{1}\circ X_{1}(t) &=&t\text{ if }t\in G_\mathscr{L_{ext}}^{c}, \\
%Y_{1}\circ X_{1}(t) &\neq &t\text{ if }t\in
%G_\mathscr{L_{ext}}\text{ and }Y_{1}\circ
%X_{1}(G_\mathscr{L_{ext}})=G_\mathscr{L_{ext}}^{c}.
\end{eqnarray*}
\end{lemma}

\proof
A proof of the lemma can be found in Theorem~5.3.3~(i) and
Corollary~5.3.5 in \cite{FIM:1999}. Parallel results when 1 is a
branch point (or both 1 and -1 are branch points) can be found in
Lemma ~2.3.6, Lemma~2.3.9 and Lemma~2.3.10 of \cite{FIM:1999}.
\QED

\begin{remark}
Results in this lemma can also be directly proved through elementary analysis without using advanced mathematical concepts used in
\cite{FIM:1999}.
\end{remark}

\section{Asymptotic Analysis of the Two Unknown Functions $\pi_1(x)$
And $\pi_2(y)$}

%\red{relationship; interlace in terms of conformal mappings; location and detailed property at the dominant singularity; Tauberian-like theorem
%(region); intuition of continuation;}

The key idea of the kernel method is to consider all $(x,y) \in B$
such that the right hand side of the fundamental form is also
zero, which provides a relationship between the two unknown
functions $\pi_1(x)$ and $\pi_2(y)$. Then, the interlace between
the unknown functions $\pi_1(x)$ and $\pi_2(y)$ plays the key role
in the asymptotic analysis of these two functions, from which
exact tail asymptotics of the stationary distribution can be
determined according to asymptotic analysis of the unknown
function at its singularities and the Tauberien-like theorem.

\subsection{Tauberian-like theorems}

Various approaches, say probabilistic or non-probabilistic,
including analytic or algebraic, are available for exact geometric
decay. However, asymptotic analysis seems unavoidable for exact
non-geometric decay.  A Tauberian, or Tauberian-like, theorem
provides a tool of connecting the asymptotic property at dominant
singularities of an analytic function at zero and the tail
property of the sequence of coefficients in the Taylor series of
the function. In our case, an unknown generating function of a
probability sequence is analytic at zero. Since these
probabilities are unknown, in general, it cannot be verified that
the probability sequence is (eventual) monotone, which is a
required condition for applying a standard Tauberian theorem. The
tool used in this paper is a Tauberian-like theorem, which does
not require this monitonicity. Instead, it imposes some extra
condition on analyticity of the unknown generating function.

Let $A(z)$ be analytic in $|z|<R$, where $R$ is the radius of convergence of the function $A(z)$. We first consider a special case in which $R$ is the only
singularity on the circle of convergence.

\begin{remark}
It should be noticed that for an analytic function at 0, if the coefficients of the Taylor expansion are all non-negative,
then the radius $R>0$ of convergence is a singularity of the function according to the well-known Pringsheim's Theorem.
\end{remark}

\begin{definition}[Definition~VI.1 in Flajolet and Sedgewick~\cite{Flajolet-Sedgewick:09}]
For given numbers $\varepsilon >0$ and $\phi$ with $0<\phi <\pi/2$, the open domain
$\Delta (\phi, \varepsilon)$ is defined by
\begin{equation}
    \Delta (\phi, \varepsilon) = \left \{z \in \mathbb{C}: |z| < 1+\varepsilon, z \neq 1, |z-1| > \phi \right \}.
\end{equation}
A domain is a $\Delta$-domain at 1 if it is a $\Delta(\phi, \varepsilon)$ for some $\varepsilon>0$ and $0< \phi< \pi/2$. For a complex number $\zeta \neq 0$,
a $\Delta$-domain at $\zeta$ is defined as the image $\zeta \cdot \Delta(\phi, \varepsilon)$ of a $\Delta$-domain $\Delta(\phi, \varepsilon)$ at 1
under the mapping $z \mapsto \zeta z$. A function is called $\Delta$-analytic if it is
analytic in some $\Delta$-domain.
\end{definition}

\begin{remark}
The region $\Delta (\phi, \varepsilon)$ is an intended disk with the radius of $1+\varepsilon$.
Readers may refer to Figure~VI.6 in \cite{Flajolet-Sedgewick:09} for a picture of the region.
Throughout the paper, without
otherwise stated, the limit of a $\Delta$-analytic function is always
taken in the $\Delta$-domain.
\end{remark}

\begin{theorem}[Tauberian-like theorem for single singularity] \label{tauberian-1}
Let $A(z)=\sum_{n\geq 0}a_{n}z^{n}$ be analytic at 0 with $R$ the
radius of convergence. Suppose that $R$ is a singularity of $A(z)$
on the circle of convergence such that $A(z)$ can be continued to
a $\Delta$-domain at $R$. If for a real number $\alpha \notin \{0,
-1, -2, \ldots\}$,
\begin{equation*}
    \lim_{z \rightarrow R}(1-z/R)^{\alpha}A(z)= g,
\end{equation*}
where $g$ is a non-zero constant, then,
\begin{equation*}
    a_{n} \sim  \frac{g}{\Gamma(\alpha)}  n^{\alpha-1} R^{-n},
\end{equation*}
where $\Gamma(\alpha)$ is the value of the gamma function at $\alpha$.
\end{theorem}

\proof
This is a immediate consequence of Corollary~VI.1 in \cite{Flajolet-Sedgewick:09} after the transform $z \mapsto R z$.
\QED

For the random walks studied in this paper, we will prove that the unknown generating function $\pi_1(x)$ ($\pi_2(y)$) has only one singularity on
the circle of its convergence, except the X-shaped random walk for which the convergent radius $R$ and $-R$ are the only singularities. To deal with the later case,
we introduce the following Tauberian-like theorem for the case of multiple singularities.

\begin{theorem}[Tauberian-like theorem for multiple singularities] \label{tauberian-2}
Let $A(z)=\sum_{n\geq 0}a_{n}z^{n}$ be analytic when $|z|<R$ and have a finite number of singularities $\zeta_k$, $k=1, 2, \ldots, m$ on the circle $|z|=R$ of convergence.
 Assume that there exists a $\Delta$-domain $\Delta_0$ at 1 such that $A$ can be continued to intersection of the $\Delta$-domains $\zeta_k$ at $\zeta_k$, $k=1, 2, \ldots, m$:
\[
    D =  \cap_{k=1}^m (\zeta_k \cdot \Delta_0).
\]
If for each $k$, there exists a real number $\alpha_k \notin \{0, -1, -2, \ldots\}$ such that
\begin{equation*}
    \lim_{z \rightarrow \zeta_k} (1-z/\zeta_k)^{\alpha_k}A(z)= g_k,
\end{equation*}
where $g_k$ is a non-zero constant, then,
\begin{equation*}
    a_{n} \sim  \sum_{k=1}^m \frac{g_k}{\Gamma(\alpha_k)}  n^{\alpha_k-1} \zeta_k^{-n}.
\end{equation*}
\end{theorem}

\proof This is an immediate corollary of Theorem~VI.5 in
\cite{Flajolet-Sedgewick:09} for the case where $\alpha _{k}$ is
real, $\beta _{k}=0$, $\sigma_{k}(z) =
\tau_{k}(z)=(1-z)^{-\alpha_{k}}$ and $\sigma_{k,n} =
\frac{g_{k}}{\Gamma (\alpha )}n^{\alpha _{k}-1}$. \QED

\subsection{Interlace of the two unknown functions $\pi_1(x)$ and $\pi_2(y)$}

The interlace of the unknown functions $\pi_1(x)$ and $\pi_2(y)$ is a key for asymptotic analysis of these functions.
Let
\begin{eqnarray*}
    \Gamma_{a} &=& \{x\in \mathbb{C:} |x|=a \}, \\
    D_{a} &=& \{x: |x|<a\}, \\
    \overline{D}_{a} &=&\{x:|x|\leq a\}.
\end{eqnarray*}
When $a=1$, we write $\Gamma=\Gamma_1$, $D = D_1$ and
$\overline{D}=\overline{D}_1$.

We fist state two literature results on the continuation of the functions $\pi_1(x)$ and $\pi_2(y)$.

%Lemma 1.3
\begin{lemma}[Theorem~3.2.3 in \cite{FIM:1999}] \label{lemma1.3}
For a stable non-singular random walk having genus 1, $\pi_1(x)$
is a meromorphic function in the complex cut plane
$\widetilde{\mathbb{C}}_{x}$. Similarly, $\pi_2(y)$ is a
meromorphic function in the complex cut plane
$\widetilde{\mathbb{C}}_{y}$.
\end{lemma}

This continuation result is crucial for tail asymptotic analysis.
The following intuition might be helpful to see why such a
continuation exist. When the right hand side of the fundamental
form is zero, the $x$ and $y$ are related, say through the
function $Y_0(x)$. Therefore, $x_3$ is the dominant singularity if
there are no other singularities exist between $(1,x_3)$. Based on
the expression for $\pi_1(x)$ obtained from the fundamental form,
all other singularities come from the zeros of $h_1(x,Y_0(x))$,
which are poles of $\pi_1(x)$, or the singularities of
$\pi_2(Y_0(x))$. A similar intuition holds for the function
$\pi_2(y)$. Based on the above intuition, it is reasonable to
expect Lemma~\ref{lemma1.3}.

\begin{remark}
An analytic continuation can be achieved through various methods.
In \cite{FIM:1999} and \cite{Flatto-Hahn:84}, it was proved in terms of
properties of Riemann surfaces. In \cite{Kobayashi-Miyazawa:2011} and \cite{Guillemin-Leeuwaarden:09}, direct methods were used
for a convergent region. For some cases, a simple proof exists by using the property of the conformal mapping $Y_0$ or $X_0$.
For example, for the case of $M_y>0$ and $M_x<0$,
we know, from Lemma~\ref{lemma1.1-b}-1, Lemma~\ref{lemma1.1} and Lemma~\ref{lemma1.1-b}-2 respectively, that
$|Y_0(x)|<1$ for $|x|=1$, $x_3>1$ and $Y_0$ is analytic in the cut
plan. Therefore, it is not difficult to see that we can find an
$\varepsilon>0$ such that for $|x|<1+\varepsilon$, the function
$\pi_2(Y_{0}(x))$ in (\ref{eqn:1.9}) is analytic, which leads to
the continuation of $\pi_1(x)$.
\end{remark}

%Lemma 1.2
\begin{lemma}[Lemma~2.2.1 in \cite{FIM:1999}] \label{lemma1.2}
Assume that the random walk is ergodic with $M \neq 0$ and the
polynomial $h(x,y)$ is irreducible. Then, exists an $\varepsilon
>0$ such that the functions $\pi_1(x)$ and $\pi_2(y)$ can be
analytically continued up to the circle $\Gamma_{1+\varepsilon }$
in their respective complex plane.  Moreover, they satisfy the
following equation in $D_{1+\varepsilon }^{2}\cap B$:
\begin{equation*}
    h_{1}(x,y)\pi_1(x)+h_{2}(x,y)\pi_2(y)+h_{0}(x,y)\pi_{0,0}=0.
\end{equation*}
\end{lemma}

\proof The analytic continuation is a direct consequence of Lemma~\ref{lemma1.3} and the equation is
directly from the fundamental form. \QED

%Theorem 1.1

\begin{theorem} \label{theorem1.1}
\textbf{1.}
Function $\pi_2(Y_{0}(x))$ is meromorphic in the cut complex plane
$\widetilde{\widetilde{\mathbb{C}}}_{x}$. Moreover, if
$Y_{0}(x_{3})$ is not a pole of $\pi_2(y)$, then $x_{3}$ is $x_{dom}$
of $\pi_2(Y_{0}(x))$ and there exist $\varepsilon >0$ and $0< \phi <\pi/2$ such that
\begin{equation*}
    \underset{ x\rightarrow x_{3}}{\lim }\pi_2(Y_{0}(x))= \pi_2(Y_{0}(x_{3}))
    \;\text{ and }\;
    \underset{ x\rightarrow x_{3}}{\lim }\pi_2^{\prime}(Y_{0}(x))=\pi_2^{\prime}(Y_{0}(x_{3})).
\end{equation*}

Similarly, $\pi_1(X_{0}(y))$ is
meromorphic in the cut complex plane $\widetilde{\widetilde{\mathbb{C}}}_{y}$. Moreover, if  $X_{0}(y_{3})$
is not a pole of $\pi_1(x)$, then $y_{3}$ is $y_{dom}$ of $\pi_1(X_{0}(y))$ and there exist $\varepsilon >0$ and $0< \phi <\pi/2$ such that
\[
    \underset{ y\rightarrow y_{3}}{\lim }\pi_1(X_{0}(y)) = \pi_1(X_{0}(y_{3}))
        \;\text{ and }\;
    \underset{y\rightarrow y_{3}}{\lim }\pi_1^{\prime}(X_{0}(y)) = \pi_1^{\prime}(X_{0}(y_{3})).
\]

\textbf{2.} In cut plane
$\widetilde{\widetilde{\mathbb{C}}}_{x}$, equation
\begin{equation} \label{eqn:1.8}
    h_{1}(x,Y_{0}(x))\pi_1(x)+h_{2}(x,Y_{0}(x))\pi_2(Y_{0}(x))+h_{0}(x,Y_{0}(x))\pi_{0,0}=0
\end{equation}
holds except at a pole (if there is any) of $\pi_1(x)$ or
$\pi_2(Y_{0}(x))$. Therefore,
\begin{equation} \label{eqn:1.9}
    \pi_1(x)=\frac{-h_{2}(x,Y_{0}(x))\pi_2(Y_{0}(x))-h_{0}(x,Y_{0}(x))\pi_{0,0}}{h_{1}(x,Y_{0}(x))},
\end{equation}
except at zero of $h_{1}(x,Y_{0}(x))$, or at a pole (if there is
any) of $\pi_1(x)$ or $\pi_2(Y_{0}(x))$.

Similarly, in the cut plane $\widetilde{\widetilde{\mathbb{C}}}_{y}$,
equation
\begin{equation}
    h_{1}(X_{0}(y),y)\pi_1(X_{0}(y)+h_{2}(X_{0}(y),y)\pi_2(y)+h_{0}(X_{0}(y),y)\pi_{0,0}=0
\end{equation}
holds except at a pole (if there is any) of $\pi_2(y)$ or
$\pi_1(X_{0}(y))$. Therefore,
\begin{equation} \label{eqn:1.11}
    \pi_2(y)=\frac{-h_{1}(X_{0}(y),y)\pi_1(X_{0}(y))-h_{0}(X_{0}(y),y)\pi_{0,0}}{h_{2}(X_{0}(y),y)},
\end{equation}
except at a zero of $h_{2}(X_{0}(y),y)$, or at a pole (if there is
any) of $\pi_2(y)$ or $\pi_1(X_{0}(y))$.
\end{theorem}

\proof We only prove the result for functions of $x$ and the result
for functions of $y$ can be proved in the same fashion.

\textbf{1.} From Lemma~\ref{lemma1.1} and Lemma~\ref{lemma1.3},
$Y_{0}(x)$ is analytic in the cut complex plane
$\widetilde{\widetilde{\mathbb{C}}}_{x}$ and $\pi_2(y)$ is
meromorphic in the cut complex plane $\widetilde{\mathbb{C}}_{y}$,
which implies $\pi_2(Y_{0}(x))$ is meromorphic in
$\widetilde{\widetilde{\mathbb{C}}}_{x}$ if $Y_{0}(x)\notin
[y_{3},y_{4}]$.  According to Lemma~\ref{lemma1.4}-2(b), for all
$x \in \mathbb{C}_x$, $Y_{0}(x)\in G_\mathscr{L} \cup
G_{\mathscr{L}_{ext}}$ and according to Lemma~\ref{lemma1.4}-1,
$[y_{3},y_{4}] \subset (G_\mathscr{L}\cup
G_{\mathscr{L}_{ext}})^{c}$, which confirms $Y_{0}(x)\notin
[y_{3},y_{4}]$. From the above, we have $\pi_2(y)$ is analytic at
$Y_0(x_3)$, then the limits in \textbf{1.} are immediate results of the
analytic properties of $\pi_2(Y_{0}(x))$.

\textbf{2.} Since both $\pi_1(x)$ and $ \pi_2(Y_{0}(x))$ are
meromorphic (proved in \textbf{1.}) and $Y_{0}(x)$ is analytic
(Lemma~\ref{lemma1.1}) in
$\widetilde{\widetilde{\mathbb{C}}}_{x}$, equation (\ref{eqn:1.8})
 in the cut
plane $\widetilde{ \widetilde{\mathbb{C}}_{x}}$ except at the
poles of $\pi_1(x)$ or $\pi_2(Y_{0}(x))$. \QED

\begin{remark} Let us
extend the definition of $\pi_1(x)$ to $x=x_3$ by $\pi_1(x_3) =
\lim_{x \to x_3} \pi_1(x)$ for $x$ in the cut plane. We say that
$x_3$ is a pole if the limit of $\pi_1(x)$ is infinite as $x \to
x_3$ in the cut plane.
\end{remark}

According to the above interlacing property and the Tauberian-like
theorem, for exact tail asymptotics of the boundary probabilities
$\pi_{n,0}$ and $\pi_{0,n}$, we only need to carry out an
asymptotic analysis at the dominant singularities of the functions
$\pi_1(x)$ and $\pi_2(y)$, respectively. There are only two
possible types of singularities, poles or branch points.  We need
to answer the following questions:

\textbf{Q1.} How many singularities on the circle of convergence (dominant singularities)?

\textbf{Q2.} What is the multiplicity of a pole?

\textbf{Q3.} Is the branch point also a pole?

For the random walk considered in this paper, we will answer all
these questions. We will see that on the convergent circle, there
is only one singularity or there are exactly two singularities.
For the former, Theorem~\ref{tauberian-1} will be applied, and for
the latter, Theorem~\ref{tauberian-2} will be applied.

\subsection{Poles of $\pi_1(x)$}

Parallel properties about poles of the function $\pi_2(y)$ can be obtained in the
same fashion, which will not be detailed here.

\begin{lemma} \label{lemma1.5}
\textbf{1.} Let $x \in G_{\mathscr{M}} \cap (\overline{D})^{c}$,
then the possible poles of $\pi_1(x)$ in $G_{\mathscr{M}}\cap
(\overline{D})^{c}$ are necessarily zeros of $h_{1}(x,Y_{0}(x))$,
and $|Y_{0}(x)|\leq 1$.

\textbf{2.} Let $ y \in G_\mathscr{L}\cap (\overline{D})^{c}$,
then the possible poles of $\pi_2(y)$ in $G_\mathscr{L}\cap
(\overline{D})^{c}$ are necessarily zeros of $ h_{2}(X_{0}(x),y)$,
and $|X_{0}(y)|\leq 1$.
\end{lemma}

\proof \textbf{1.} When $x\in
\mathscr{M}$, then $Y_0(x) =y \in [y_1,y_2]$.
From
Lemma~\ref{lemma1.1}, for $|x|=1$, $|Y_{0}(x)|\leq 1$. For $x \in
G_{\mathscr{M}} \cap (\overline{D})^{c}$, it follows from the
maximum modulus principle, we have $|Y_{0}(x)|\leq 1$. Hence,
$\pi_2(Y_{0}(x))$ is analytic in $G_{\mathscr{M}} \cap
(\overline{D})^{c}$. From Theorem~\ref{theorem1.1}, if
$h_{1}(x,Y_{0}(x)) \neq 0$,  equation (\ref{eqn:1.9}) holds, which
implies that the possible poles of $\pi_1(x)$ in $G_{\mathscr{M}}
\cap (\overline{D})^{c}$ are necessarily zeros of
$h_{1}(x,Y_{0}(x))$.

\textbf{2.} The proof is similar. \QED

\begin{theorem} \label{theorem1.2}
 Let $x_p$ be a pole of $\pi_1(x)$ with the
smallest modulus. Assume that $|x_p| \leq x_3$. Then, one of the
follow two cases must hold:

\textbf{1.} $x_p$ is a zero of $h_{1}(x,Y_{0}(x))$;

\textbf{2.} $\widetilde{y}_0=Y_{0}(x_p)$ is a zero of
$h_{2}(X_{0}(y),y)$ and $|\widetilde{y}_0|>1$.

Parallel results hold for a pole of $\pi_2(y)$.
\end{theorem}

\proof Suppose that $x_p$ is not a zero of $h_{1}(x,Y_{0}(x))$.
According to equation (\ref{eqn:1.9}) in Theorem~\ref{theorem1.1},
$x_p$ must be a pole of $\pi_2(Y_{0}(x))$ and
$|\widetilde{y}_0|>1$. Furthermore, by Lemma~\ref{lemma1.5},
$x_p \notin G_{\mathscr{M}}$. If $\widetilde{y}_0$ is not a
zero of $h_{2}(X_{0}(y),y)$, according to equation
(\ref{eqn:1.11}) in Theorem~\ref{theorem1.1},
$\widetilde{y}_0$ must be a pole of $\pi_1(X_{0}(y))$, that
is, $\widetilde{x}_0=X_0(\widetilde{y}_0)$ is a pole of $\pi_1(x)$.
It follows from Lemma~\ref{lemma1.5} that
$\widetilde{x}_0=X_{0}(\widetilde{y}_0)$ is a zero of
$h_{1}(x,Y_{0}(x))$ if $\widetilde{x}_0\in G_{\mathscr{M}}$.
 There are two possible cases:
$\Delta >0$ or $\Delta \leq 0$. If $\Delta >0$, by
Lemma~\ref{lemma1.4}-1(a) and 2(c), $\widetilde{x}_0 \in
G_{\mathscr{M}}$. In the case of $\Delta \leq 0$, according to
Lemma~\ref{lemma1.4}-1(b), 1(c) and 2(b), we also have
$\widetilde{x}_0\in G_{\mathscr{M}}$. However, this case is not
possible, since otherwise according to Lemma~\ref{lemma1.4}-1 we
would have $\widetilde{x}_0=x_p$ or $\widetilde{x}_0=-x_p$, both leading
to a contradiction. This completes the proof. \QED

%\textbf{3.} $x^{\ast\ast}=X_{0}(y^{\ast})$ is a zero of
%$h_{1}(x,Y_{0}(x))$ and $|x^{\ast\ast}|>1$.

\begin{remark} \label{remark1.3}
We will show in the next subsection that a pole
of $\pi_1(x)$ with the smallest modulus in the disk $|x|\leq x_{3}$ is real.
\end{remark}

\subsection{Zeros of $h_{1}(x,Y_{0}(x))$}
\label{section2}

In this subsection, we provide properties on the zeros of the
function $h_1(x,Y_0(x))$. The main result is stated in the
following theorem.
\begin{theorem} \label{theorem2.1}
For a non-singular random walk having genus 1, consider
the following two possible cases:

\textbf{1.} Either $p_{i,j}$ or $p_{i,j}^{(1)}$ is not X-shaped.
In this case, either $h_{1}(x,Y_{0}(x))$ has no zeros with modulus
in $(1,x_{3}]$, or it has only one simple zero, say $x^{\ast }$,
with modulus in $(1,x_{3}]$, and $x^{\ast }$ is positive.

\textbf{2.} Both $p_{i,j}$ and $p_{i,j}^{(1)}$ are X-shaped. In
this case, either $h_{1}(x,Y_{0}(x))$ has no zeros with modulus in
$(1,x_{3}]$, or it has exact two simple zeros, namely,
$x^{\ast}>0$ (with modulus in $(1,x_3]$) and $-x^{\ast}$, both are
zeros of $h_{1}(x,Y_{0}(x))$ or both are zeros of
$a(x)h_{1}(x,Y_{1}(x))$.
\end{theorem}

With this theorem and Theorem~\ref{theorem1.2}, we are able to
apply the Tauberian-like theorem to characterize the tail
asymptotic properties for the boundary probability sequence
$\pi_{n,0}$. To show the above Theorem, we need the following
several lemmas and two propositions. Instead of directly
considering the function $f_0(x)=h_{1}(x,Y_{0}(x))$, we consider a
polynomial $f(x)$, which is essentially the product of  $f_0(x)$
and $f_1(x)=h_{1}(x,Y_{1}(x))$:
\[
    f(x)= f_0(x)\widetilde{f_{1}}(x),
\]
where $\widetilde{f_{1}}(x)=a(x)f_{1}(x)$. It is easy to verify,
by noticing
\begin{equation*}
    Y_{0}(x)Y_{1}(x)=\frac{c(x)}{a(x)}\; \text{ and } \; Y_{0}(x)+Y_{1}(x)=-\frac{b(x)}{ a(x)},
\end{equation*}
that
\begin{eqnarray}
    f(x) &=& a(x)b_{1}^{2}(x)-b(x)b_{1}(x)a_{1}(x)+c(x)a_{1}^{2}(x)  \label{lemma2.1} \\
    &=&
    d_{6}x^{6}+d_{5}x^{5}+d_{4}x^{4}+d_{3}x^{3}+d_{2}x^{2}+d_{1}x+d_{0}.
\end{eqnarray}
Hence, a zero of $f_{i}(x)$, $i=0,1$, has to be a zero of $f(x)$,
and any zero of $f(x)$ is either a zero of $f_{0}(x)$ or a zero of
$\widetilde{f_{1}}(x)=a(x)f_{1}(x)$.

We can also write
\begin{equation}
    f(x) =  a(x)[a_{1}(x)]^{2}R_{-}(x)R_{+}(x), \label{eqn:2.1}
\end{equation}
where
\begin{equation}
    R_{\pm}(x)=F(x)\pm \frac{\sqrt{D_{1}(x)}}{2a(x)}
\end{equation}
with
\begin{equation}
    F(x)=\frac{b_{1}(x)}{a_{1}(x)}-\frac{b(x)}{2a(x)}.
\end{equation}

\begin{remark} \textbf{1.}
It can be easily seen that both $f_{0}(x)$ and
$\widetilde{f_{1}}(x)$ are analytic on the cut complex plan. In
fact, the analyticity of $f_0(x)$ is obvious and the analyticity
of $\widetilde{f}_1(x)$ is due to the cancellation of the zeros of
$a(x)$ and the pole of $f_1(x)$.

%\textbf{3.} We could not find a short proof. A different proof
%might be possible under the 1-arithmetic condition.
\end{remark}

All proofs for Lemmas~\ref{lemma2.2}--\ref{lemma2.6} and for
Proposition~\ref{theorem2.2} and Proposition~\ref{theorem2.3} are
organized into Appendix~\ref{appendix1}.

\begin{lemma} \label{lemma2.2}
\textbf{1.} \textbf{(a)} $Y_{0}^{\prime}(1)=\frac{M_{x}}{-M_{y}}$ if $M_{y}<0$;
\textbf{(b)} $Y_{1}^{\prime}(1)=\frac{M_{x}}{-M_{y}}$ if $M_{y}>0$; and
\textbf{(c)} $Y_{1}(1)=Y_{0}(1) $ and $x=1$ is a branch point of $Y_{1}(x)$ and $Y_{0}(x)$
if $M_{y}=0$. In this case, $Y_{1}^{\prime}(1)$ and $Y_{0}^{\prime}(1)$ do
not exist. Parallel results hold for functions $X_{k}(y)$.

\textbf{2.} If $M_{y}\neq 0$, then $f(x)$ has at least one
non-unit zero in $[x_{2},x_{3}]$ and 1 is a simple zero of $f(x)$.
Parallel results holds for the case of $M_{x}\neq 0$.

\end{lemma}

\begin{lemma} \label{lemma2.3}
\textbf{1.} Let $z$ be a branch point of $Y_{0}(x)$. If $f(z)=0$,
then $z$ cannot be a repeated root of $f(x)=0$.

\textbf{2.} $f(x)$ (therefore both $f_{0}(x)$ and $\widetilde{f}_{1}(x)$)
has (have) no zeros on the cuts, except possibly at a branch point. More
specifically, $f(x)<0$ if $a(x)<0$ and $f(x)>0$ if $a(x)>0$.

\textbf{3.} $f_{0}(x)$ and $\widetilde{f}_{1}(x)$ have no common
zeros except possibly at a branch point or at zero.

\textbf{4.} Consider the random walk in
Theorem~\ref{theorem2.1}-1. If $f_{0}(x)$ has a zero in
$[-x_{3},-1)$, then $f_{0}(x)$ has an additional (different) zero
in $[-x_{3},-1)$.

\textbf{5.} For the random walk in Theorem~\ref{theorem2.1}-1, if
$|x|\in (1,x_{3}]$, then $|Y_{0}(-|x|)|<Y_{0}(|x|)$.
\end{lemma}

\begin{lemma} \label{lemma2.4}
Consider the random walk in Theorem~\ref{theorem2.1}-1.
If $ M_{y}\leq 0$, then $x=1$ is the only zero of
$f_{0}(x)=h_{1}(x,Y_{0}(x))$ on the unit circle $|x|=1$. If
$M_{y}>0$, then $f_{0}(x)$ has no zero on unit circle $|x|=1$.
\end{lemma}

\begin{remark} \label{remark2.1}
From the proof of Lemma~\ref{lemma2.4}, we can see that for the
random walk considered in Theorem~\ref{theorem2.1}-2, $f_{0}(x)$
has no zeros with non-zero imaginary part on the unit circle.
\end{remark}

The proof of Theorem~\ref{theorem2.1} is based on detailed
properties of the function $f(x)$ and also the powerful continuity
argument to connect an arbitrary random walk to a simpler one. For
using this continuity argument, we consider the following special
random walk.

\textbf{Special Random Walk.} This is the random walk for which $p_{i,j}$ is cross-shaped (or $p_{i,j}=0$
whenever $|ij|=1$), and $p_{-1,1}^{(1)}=p_{-1,0}^{(1)}=0$.
We first prove the counterpart result to Theorem~\ref{theorem2.1} for the Special Random Walk.

\begin{proposition} \label{theorem2.2}
For the Special Random Walk, the following results hold:

\textbf{1.} $f(x)=0$ has six real roots with exact one non-unit root in
$[x_{2},x_{3}]$. More specifically, two roots are zero, two in $[x_{2,}x_{3}]$, one in $(-\infty ,x_{1}]$, and one
in $[x_{4},\infty )$.

\textbf{2.} If $f_{0}(x)$ has a zero, say $x^{\ast}$, in $(1$,
$x_{3}]$, then $x^{\ast}$ is the only zero of $f_{0}(x)$ with
modulus in $(1$, $x_{3}]$. Furthermore, $f_{0}(x)$ has no other
zeros with modulus greater than $1$ except possibly at $x=x_{4}$.
\end{proposition}

For the random walk considered in Theorem~\ref{theorem2.1}-2, we first prove the following results.
\begin{lemma} \label{lemma2.6}
For the random walk considered in Theorem~\ref{theorem2.1}-2 (or both $p_{i,j}$ and $p^{(1)}_{i,j}$ are X-shaped),
 $f(1)=f(-1)=0$, and $f(x)=0$ has two more real roots, say $0<x_0 \neq 1$ and $-x_0$, and two complex roots.
\end{lemma}

\begin{proposition} \label{theorem2.3}
For the random walk considered in Theorem~\ref{theorem2.1}-2 (or
both $p_{i,j}$ and $p^{(1)}_{i,j}$ are X-shaped), either the
two complex zeros of $f(x)$ are zeros of
$\widetilde{f_{1}}(x)=a(x)f_{1}(x)$ or they are inside the unit
circle.
\end{proposition}

\underline{\proof of Theorem~\ref{theorem2.1}.}
\textbf{1.}
For the random walk considered here (either $p_{i,j}$ or $p^{(1)}_{i,j}$ is not X-shaped), let
\begin{eqnarray*}
    \mathbf{p}&=&(p_{-1,-1}, p_{0,-1}, p_{1,-1}, p_{-1,0},p_{0,0},p_{0,1}, p_{-1,1},p_{0,1},p_{1,1}), \\
    \mathbf{p}^{(1)}&=&(p^{(1)}_{-1,0},p^{(1)}_{0,0},p^{(1)}_{0,1}, p^{(1)}_{-1,1},p^{(1)}_{0,1},p^{(1)}_{1,1}).
\end{eqnarray*}
Define
\begin{equation*}
    A=\Big \{ \big ( \mathbf{p},\mathbf{p}^{(1)} \big ): 0 \leq p_{i,j}, p^{(1)}_{i,j} \leq 1 \; \mbox{ and } \; \sum_{i,j} p_{i,j} = \sum_{i,j} p^{(1)}_{i,j} =1 \Big \}.
\end{equation*}

For an arbitrary random walk for which either $p_{i,j}$ or
$p^{(1)}_{i,j}$ is not X-shaped, let $\rho$ be the corresponding
point in $A$. We assume that $M_y \leq 0$ for the random walk
$\rho$ (and a similar proof can be found for the case of $M_y>0$).
Let $\rho_{0}$ be an arbitrarily chosen point in $A$ corresponding
the Special Random Walk. We prove the result by contradiction.
Suppose otherwise that the statement were not true. There would be
three possible cases: (i) $\Im (x^{\ast })\neq 0$; (ii)
$-x_{3}\leq x^{\ast}< -1$; and (iii) there exists $x_{0} \in
(1,x_{3}]$ with $x_{0}\neq x^{\ast}$ such that $f_{0}(x_{0})=0$.

\textbf{Case (i).} Clearly, $\overline{x^{\ast}}$ is also a root
of $f(x)=0$. Choose a simple connected path $\ell$ in $A$ to
connect $\rho$ to $\rho_{0}$ such that on $\ell$ (excluding
$\rho$, but including $\rho_{0}$) $M_{y}<0$. The zeros of $f(x)$
as a function of parameters in $A$ are continues on $\ell$. There
are two possible cases: (a) the zero function $x_0(\theta)$ (with
$x_0(\rho)=x^*$) never passes the unit circle when $\theta$
travels from $\rho$ to $\rho_0$; and (b) $x_0(\theta)$ passes the
unit circle at some point $\theta \in \ell$.

If (a) occurs, let $\theta_0$ be the first point at which
$x_0(\theta)=\overline{x}_0(\theta)$, where
$\overline{x}_0(\theta)$ is the zero function with
$\overline{x}_0(\rho)=\overline{x^*}$. If $\overline{x^*}$ is a
zero of $\widetilde{f}_{1}$, then $f_0$ and $\widetilde{f}_{1}$
would have a common zero  $x_0(\theta_0)=\overline{x}(\theta_0)$
at $\theta_0$, which contradicts Lemma~\ref{lemma2.3}-3. Hence,
the only possibility is that $\overline{x^*}$ is also a zero of
$f_0$. From $\theta_0$ on, both $x_0(\theta)$ and
$\overline{x}(\theta)$ should always be zeros of $f_0$, since
otherwise only at a branch point a zero of $f_0$ could be switched
to a zero of $\widetilde{f}_1$ and all branch points are real,
which means that $x_0(\theta)=\overline{x}(\theta)$ is a branch
point and a multiple roots, contradicting to
Lemma~\ref{lemma2.3}-1. As $\theta_0$ approaches $\rho_0$, it
leads to a contradiction that two zeros of $f_0$ are in $(1,x_3]$.

If (b) occurs, we can assume that when $x_0(\theta)$ pases the
unit circle it is a zero of $f_0$ based on the proof in (a). Then,
$f_0$ has two zeros since 1 is always a zero of $f_0$ independent
of the parameters (or $\theta$) when $M_y<0$, which is a different
zero from $x_0(\theta)$. This contradicts to the fact that $f_0$
has only one zero at the unit circle.

\textbf{Case (ii).} In this case, $f_0(x)$ would have another zero
in $[-x_3,-1)$ at $\rho$ according to Lemma~\ref{lemma2.3}-4.
Consider the same two cases (a) and (b) as in (i). We can then
follow a similar proof to show that case (ii) is impossible.

\textbf{Case (iii).} A similar proof will show that the case is
impossible.

\textbf{2.} This is a direct consequence of Lemma~\ref{lemma2.6}
and Proposition~\ref{theorem2.3}. \QED

The following Lemma gives a necessary and sufficient condition
under which $f_{0}(x)=h_1(x,Y_0(x))$ has a zero in $(1,x_{3}]$.
\begin{lemma} \label{lemma2.7}
Assume $M_{y}\neq 0$. We have following
results:

\textbf{1.} If $f_{0}(x_{3})\geq 0$, $f_{0}(x)$ has a zero in $(1,x_{3}]$;

\textbf{2.} If $f_{0}(x_{3})<0$, $f_{0}(x)$ has no zeros in $(1,x_{3}]$.
\end{lemma}

\proof
\textbf{1.} There are two cases: $M_{y}>0$ or $M_{y}<0$.
If $M_{y}>0$, then $f_{0}(1)<0$, which leads to the conclusion. If $M_{y}<0$, then
$f_{0}^{\prime}(1)<0$, which also leads to the conclusion since $f_{0}(1)=0$ and
$f_{0}(x_{3})\geq 0$.

\textbf{2.} Again there are two cases: $M_{y}>0$ or $M_{y}<0$. By simple
calculus, in either case, we obtain that if $f_{0}(x)=0$ had a root in $(1,x_{3}]$,
then it would have another root in $(1,x_{3}]$ since $f_{0}(x_{3})<0$. This
contradicts to Theorem~\ref{theorem2.1}.
\QED

\subsection{Zeros of $h_{2}(X_{0}(y),y)$}

Following the same argument in the previous subsection, we have
the following result:
\begin{theorem} \label{theorem2.1-b}
For a non-singular random walk having genus 1, consider the
following two possible cases:

\textbf{1.} Either $p_{i,j}$ or $p_{i,j}^{(2)}$ is not X-shaped.
In this case, either $h_{2}(X_{0}(y),y)$ has no zeros with modulus
in $(1,y_{3}]$, or it has only one simple zero, say $y^{\ast }$,
with modulus in $(1,y_{3}]$, and $y^{\ast }$ is positive.

\textbf{2.} Both $p_{i,j}$ and $p_{i,j}^{(2)}$ are X-shaped. In
this case, either $h_{2}(X_{0}(y),y)$ has no zeros with modulus in
$(1,y_{3}]$, or it has exact two simple zeros, namely,
$y^{\ast}>0$ (with modulus in $(1,y_3]$) and $-y^{\ast}$, both are
zeros of $g_{0}(y)$ or both are zeros of $g_{1}(y)$, where
\[
    g_0(y)=h_{2}(X_{0}(y),y) \qquad \text{and} \qquad g_1(y)=h_{2}(X_{1}(y),y).
\]
\end{theorem}

From the above analysis, we know that if $h_{1}(x,Y_{0}(x))$ has a
zero in $(1,x_{3}]$, then such a zero is unique. Similarly, if
$h_{2}(X_{0}(y),y)$ has a zero in $(1,y_{3}]$, then such a zero is
unique. For convenience, we make the following convention:

\begin{convention}  Let $x^{\ast}$ be the unique
zero in $(1,x_{3}]$ of the function $h_{1}(x,Y_{0}(x))$, if such a
zero exists, otherwise let $x^{\ast}=\infty$. Similarly Let
$y^{\ast}$ be the unique zero in $(1,y_{3}]$ of the function
$h_{2}(X_{0}(y),y)$ if such a zero exists, otherwise let
$y^{\ast}=\infty$.
\end{convention}

According to  Theorem~\ref{theorem1.2}, the unique pole in $(1,
x_3]$ of $\pi_1(x)$ is either $x^{\ast}$, or the image of the pole
under $Y_0$ is a zero of $h_2(X_0(y),y)$. Our focus in this
subsection is on this special case of $y^*$.

\begin{theorem} \label{theorem-h2}
If the pole in $(1, x_3]$ of $\pi_1(x)$ is not $x^{\ast}$, then,
it, denoted by $\widetilde{x}_1$, satisfies:

\textbf{1.} $\widetilde{x}_1=X_1(y^*)$, where $y^*$ is the unique
zero in $(1,y_{3}]$ of the function $h_{2}(X_{0}(y),y)$;

\textbf{2.} $\widetilde{x}_1$ is the only pole of $\pi_1(x)$ with
modulus in $(1,y_3]$, except for the case where both $p_{i,j}$ and
$p_{i,j}^{(2)}$ are X-shaped, for which $-\widetilde{x}_1$ is the
other pole of $\pi_1(x)$ with modulus in $(1,y_3]$.
\end{theorem}

\proof \textbf{1.} Let $\widetilde{x}$ be the solution of $y^*=Y_{0}(x)$.
Then,
$\widetilde{x}=\widetilde{x}_{0}\stackrel{\triangle}{=}X_{0}(y^*)$
or
$\widetilde{x}=\widetilde{x}_{1}\stackrel{\triangle}{=}X_{1}(y^*)$.
If $y^*\in G_\mathscr{L}$, then $\widetilde{x}=\widetilde{x}_{0}$
so that $y^*=Y_{0}(X_{0}(y^*))$. In this case, by
Lemma~\ref{lemma1.5}, $\widetilde{x}_{0}<1$. If $y^*\in
G_\mathscr{L}^{c}$, then $\widetilde{x}=\widetilde{x}_{1}$ so that
$y^*=Y_{0}(X_{1}(y^*)$ and $\widetilde{x}_{1}\in G_{\mathscr{M}}^{c}$.

\textbf{2.} It follows from the fact that the zero, $y^{\ast}$, of
$h_{2}(X_{0}(y),y)$ in $(1,y_{3}]$ is unique and the fact that
$y^{\ast}=Y_{0}(x)$ has only two possible solutions
$\widetilde{x}_{0}<1$ and $\widetilde{x}_{1}$. In the case where
both $p_{i,j}$ and $p_{i,j}^{(2)}$ are X-shaped, $-y^{\ast }$ is
the other zero of $h_{2}(X_{0}(y),y)$ with either $-y^{\ast
}=Y_{0}(-\widetilde{x}_{1})$ or $-y^{\ast
}=Y_{0}(-\widetilde{x}_{0})$.
 \QED

\begin{corollary} Let $\widetilde{x}$ be a solution of $y^*=Y_{0}(x)$.
In order for $\widetilde{x}$ to be in $(1,x_3]$ we need $y^* \in
G_\mathscr{L}^{c}$. Furthermore, we have
$y^*<y_{3}$.
\end{corollary}

\proof
The first conclusion is directly from the proof to Theorem~\ref{theorem-h2} and the second one
follows from that fact that by Lemma~\ref{lemma1.4}-1 and
Lemma~\ref{lemma1.4}-2(b), there exists no $x \in (1, x_3]$ such that
$y^*=y_3=Y_0(x)$. Therefore, we should have
$y^*<y_{3}$.
\QED

\begin{convention}  Let $\widetilde{x}_1=X_1(y^*)$ if the unique zero $y^*$ in
$(1,y_{3}]$ of the function $h_{2}(X_{0}(y),y)$ exists, otherwise
let $\widetilde{x}_1=\infty$.
\end{convention}

\subsection{Asymptotics behaviour of $\pi_1(x)$ and $\pi_2(y)$}
\label{section3}

In this subsection, we provide asymptotic behaviour of two unknown functions $\pi_1(x)$ and $\pi_2(y)$.
We only provide details for $\pi_1(x)$, since the behaviour for $\pi_2(y)$ can be characterized in the same fashion.

It follows from the discussion so far that:
\begin{description}
\item[(1)] If $p_{i,j}$ is not X-shaped, then, independent of the
properties of $p_{i,j}^{(1)}$ and $p_{i,j}^{(2)}$, there is only
one dominant singularity, which is the smallest one of $x^*$,
$\widetilde{x}_1$ and $x_3$. Here $x^*$, $\widetilde{x}_1$ and
$x_3$ are not necessarily all different.

\item[(2)] If $p_{i,j}$ is X-shaped, then both $x_3$ and $-x_3$
are branch points.
\begin{description}
 \item[(a)] If $p_{i,j}^{(1)}$ is not X-shaped, then $h_1(x,Y_0(x)$
has either no zero or one zero $x^*$ in $(1,x_3]$; and if
$p_{i,j}^{(1)}$ is X-shaped, then $h_1(x,Y_0(x)$ has either no
zero or two zeros $x^*  \in (1,x_3]$ and $-x^*$.
 \item[(b)] Similar to (a), $h_2(X_0(y),y)$ has either no zero in
 $(1,y_3]$ or one zero $y^*$ in it. For the latter, if $p_{i,j}^{(2)}$ is not X-shaped, then
$\widetilde{x}_1=X_1(y^*)$ is the only pole of $\pi_2(Y_0(x))$
with modulus in $(1,x_3]$; and if $p_{i,j}^{(2)}$ is X-shaped,
then $\widetilde{x}_1=X_1(y^*) \in (1,x_3]$ and
$-\widetilde{x}_1=X_1(-y^*)$ are the only two poles of
$\pi_2(Y_0(x))$ with modulus in $(1,x_3]$.
\end{description}
Therefore, in case (2), we either have only one dominant
singularity or exactly two dominant singularities depending on
which of $x^*$, $\widetilde{x}_1$ and $x_3$ is smallest and the
property of $p_{i,j}^{(k)}$, $k=1,2$.
\end{description}

The theorem in this subsection provides detailed asymptotic
properties at a dominant singularity for all possible cases. Let $x_{dom}$ be a dominant singularity of $\pi_1(x)$. Clearly,
$|x_{dom}|=$ $x^{\ast}$, $|x_{dom}|=\widetilde{x}_{1}$ or
$|x_{dom}|=x_{3}$.
To
state this theorem for the cases
where $x_{dom}=\pm x_{3}$, notice that through simple calculation we can
write
\begin{equation}
\label{eqn:3.1}
    h_{1}(x,Y_{0}(x))=p_{1}(x)+q_{1}(x)\sqrt{1-\frac{x}{x_{dom}}},
\end{equation}
\begin{equation} \label{eqn:3.2}
    Y_{0}(x)=p(x)+q(x)\sqrt{1-\frac{x}{x_{dom}}},
\end{equation}
\begin{equation} \label{eqn:3.3}
    Y_{0}(x_{dom})-Y_{0}(x) = \left(1-\frac{x}{x_{dom}}\right)p^{\ast}(x)-q(x)\sqrt{1-\frac{x}{x_{dom}}},
\end{equation}
\begin{equation} \label{eqn:3.4}
    h_{1}(x,Y_{0}(x))-h_{1}(x_{dom},Y_{0}(x_{dom})) = \left( 1-\frac{x}{x_{dom}}
\right) p_{1}^{\ast }(x)+q_{1}(x)\sqrt{1-\frac{x}{x_{dom}}},
\end{equation}
where
\begin{align*}
p(x) &=\frac{-b(x)}{2a(x)}, \quad p^{\ast }(x)=\frac{\frac{b(x)}{2a(x)}-
\frac{b(x_{dom})}{2a(x_{dom})}}{\frac{1}{x_{dom}}(x-x_{dom})}, \quad
p_{1}(x)=\frac{-b(x)a_{1}(x)}{2a(x)}+b_{1}(x), \\
    p_{1}^{\ast }(x) &= x_{dom}\left( \frac{a_{1}(x)-a_{1}(x_{dom})+b_{1}(x)-b_{1}(x_{dom})}{x_{dom}-x}\right),
\end{align*}
\[
q(x)=\left\{\begin{array}{ll}
-\frac{1}{2a(x)}\sqrt{\frac{D_{1}(x)}{1-\frac{x}{x_{dom}}}}, & \text{ if } x_{dom}=x_{3}, \\
\frac{1}{2a(x)}\sqrt{\frac{D_{1}(x)}{1-\frac{x}{x_{dom}}}}, & \text{ if } x_{dom}=-x_{3},
\end{array}
\right.
\]
and $q_{1}(x)=a_{1}(x)q(x)$.

Define
\begin{eqnarray*}
    L(x) &=&\frac{[h_{2}(x,Y_{0}(x))\pi_2(Y_{0}(x))+h_{0}(x,Y_{0}(x))\pi_{0,0}]h_{1}(x,Y_{1}(x))a(x)}{xf^{\prime}(x)},  \\
    \widetilde{L}(y)
    &=&\frac{[h_{1}(X_{0}(y),y)\pi_1(X_{0}(y))+h_{0}(X_{0}(y),y)\pi_{0,0}]h_{2}(X_{1}(y),y)\widetilde{a}(y)}{yg^{\prime}(y)},
\end{eqnarray*}
where $f(x)=a(x)h_{1}(x,Y_{0}(x))h_{1}(x,Y_{1}(x))$ is a
polynomial defined in Section~\ref{section2} and
$g(y)=\widetilde{a}(y)h_{2}(X_{1}(y),y)$ $h_{2}(X_{0}(y),y)$ is
the counterpart polynomial for function $h_2$.

The following Theorem shows the behaviour of $\pi_1(x)$ at
$x_{dom}$. Recall that $\widetilde{y}_0=Y_0(x^*)$.
\begin{theorem} \label{theorem3.1}
Assumed that both
$h_{2}(x^{\ast},Y_0(x^{\ast}))\pi_2(Y_0(x^{\ast}))+h_{0}(x^{\ast},Y_0(y^{\ast}))\pi_{0,0}
\neq 0$ and $h_{1}(X_{0}(\widetilde{y}_0),\widetilde{y}_0)$
$\pi(X_{0}(\widetilde{y}_0))+h_{0}(X_{0}(\widetilde{y}_0),\widetilde{y}_0)\pi_{0,0}
\neq 0$. For the function $\pi_1(x)$, a total of four types of
asymptotics exist as $x$ approaches to a dominant singularity of
$\pi_1(x)$, based on the detailed property of the dominant
singularity.

\textbf{Case 1}: If $|x_{dom}|=x^{\ast}<\min
\{\widetilde{x}_{1},x_{3}\}$, or $|x_{dom}|=\widetilde{x}_{1} <
\min\{x^{\ast},x_{3}\}$, or
$|x_{dom}|=x^{\ast}=\widetilde{x}_{1}=x_{3}$, then
\begin{equation*}
    \lim_{x\rightarrow x_{dom}}\left( 1-\frac{x}{x_{dom}}\right) \pi_1(x)=c_{0,1}(x_{dom}),
\end{equation*}
where
\[
    c_{0,1}(x_{dom}) = \left \{ \begin{array}{ll}
L(x_{dom}), & \text{if } x^{\ast}<\min \{\widetilde{x}_{1},x_{3}\}; \\
\displaystyle
\frac{-h_{2}(x_{dom},\widetilde{y}_0)\widetilde{y}_0\widetilde{L}(\widetilde{y}_0)}{h_{1}(x_{dom},\widetilde{y}_0)Y_{0}^{\prime}(x_{dom})x_{dom}},
&
    \text{if }\widetilde{x}_{1}<\min \{x^{\ast},x_{3}\}; \\
\displaystyle
\frac{h_{2}(x_{dom},\widetilde{y}_0)\widetilde{L}(\widetilde{y}_0)\widetilde{y}_0}{q_{1}(x_{dom})q(x_{dom})},
& \text{if } x^{\ast}=\widetilde{x}_{1}=x_{3}, \end{array} \right.
\]
with $\widetilde{y}_0=Y_{0}(x_{dom})$.

\textbf{Case 2}: If $|x_{dom}|=x^{\ast}=x_{3}<\widetilde{x}_{1}$
or $|x_{dom}|=\widetilde{x}_{1}=x_{3}<x^{\ast}$, then
\begin{equation*}
    \lim_{\frac{x}{x_{dom}}\rightarrow 1}\sqrt{1-x/x_{dom}}\pi_1(x)=c_{0,2}(x_{dom}),
\end{equation*}
where
\[
    c_{0,2}(x_{dom}) = \left \{ \begin{array}{ll}
\displaystyle
\frac{h_{2}(x_{dom},\widetilde{y}_0)\pi_2(\widetilde{y}_0)+h_{0}(x_{dom},\widetilde{y}_0)\pi_{0,0}}{-q_{1}(x_{dom})}, & \text{if }x^{\ast}=x_{3}<\widetilde{x}_{1}; \\
\displaystyle
    \frac{h_{2}(x_{dom},\widetilde{y}_0)\widetilde{y}_0\widetilde{L}(\widetilde{y}_0)}{h_{1}(x_{dom},\widetilde{y}_0)q(x_{dom})}, &   \text{if }
\widetilde{x}_{1}=x_{3}<x^{\ast}, \end{array} \right.
\]
with $\widetilde{y}_0=Y_{0}(x_{dom})$.

\textbf{Case 3:} If $|x_{dom}|=x_{3}<\min
\{\widetilde{x}_{1},x^{\ast}\}$, then
\begin{equation*}
    \lim_{x\rightarrow x_{dom}}\sqrt{1-x/x_{dom}}\pi_1^{\prime}(x)=c_{0,3}(x_{dom}),
\end{equation*}
where
\begin{equation*}
    c_{0,3}(x_{dom})=-\frac{q(x_{dom})x_{dom}^{2}}{2}\frac{d}{dy}\left [
\frac{h_{2}(x_{dom},y)\pi_2(y)+h_{0}(x_{dom},y)\pi_{0,0}}{-h_{1}(x_{dom},y)}
\right ] \bigg |_{y=Y_{0}(x_{dom})}.
\end{equation*}

\textbf{Case 4:} If $|x_{dom}|=x^{\ast}=\widetilde{x}_{1}<x_{3}$,
then
\begin{equation*}
    \lim_{x\rightarrow x_{dom}}\left( 1-\frac{x}{x_{dom}}\right)
    ^{2}\pi_1(x)=c_{0,4}(x_{dom}),
\end{equation*}
where
\begin{equation*}
    c_{0,4}(x_{dom})=\frac{h_{2}(x_{dom},\widetilde{y}_0)[h_{1}(\widetilde{x}_0,\widetilde{y}_0)
\pi_1(\widetilde{x}_0)+h_{0}(\widetilde{x}_0,\widetilde{y}_0)]\pi_{0,0}]}{x^{\ast2}h_{1}^{\prime}(x_{dom},\widetilde{y}_0)
Y_{0}^{\prime}(x_{dom})h_{2}^{\prime}(X_{0}(\widetilde{y}_0),\widetilde{y}_0)},
\end{equation*}
with $\widetilde{y}_0=Y_{0}(x_{dom})$ and
$\widetilde{x}_0=X_{0}(\widetilde{y}_0)$.

\end{theorem}

\proof
\textbf{Case~1.} If $x^{\ast}<\widetilde{x}_{1}$, then $x_{dom}$
is not a pole of $\pi_2(Y_{0}(x))$. According to
Theorem~\ref{theorem2.1}, $x_{dom}$ is a simple pole of
$\pi_1(x)$. From equation (\ref{eqn:1.9}) in
Theorem~\ref{theorem1.1} and Lemmas~\ref{lemma2.2} and
\ref{lemma2.3}, we have
\begin{eqnarray*}
    \pi_1(x) &=&\frac{-h_{1}(x,Y_0(x))\pi_2(Y_0(x))-h_{0}(x,Y_0(x))\pi_{0,0}}{h_{1}(x,Y_0(x))} \\
    &=&\frac{-[h_{1}(x,Y_0(x))\pi_2(Y_0(x))+h_{0}(x,Y_0(x))\pi_{0,0}]h_{1}(x,Y_{1})a(x)}{f(x)} \\
    &=&\frac{-[h_{1}(x,Y_0(x))\pi_2(Y_0(x))+h_{0}(x,Y_0(x))\pi_{0,0}]h_{1}(x,Y_{1})a(x)}{(x-x_{dom})f^{\ast}(x)},
\end{eqnarray*}
where $f^{\ast}(x_{dom})=f^{\prime}(x_{dom})\neq 0$. It follows
that
\begin{equation*}
    \lim_{x\rightarrow x_{dom}}\left( 1-\frac{x}{x_{dom}}\right)\pi_1(x)=L(x_{dom}).
\end{equation*}

Similarly, if $\widetilde{x}_{1}<x^{\ast}$, following the same
argument used in the above, we have
 $\lim_{y\rightarrow \widetilde{y}_0}\left( 1-\frac{y}{\widetilde{y}_0} \right)
\pi_2(y)=\widetilde{L}(\widetilde{y}_0)$ and
\begin{eqnarray*}
    &&\lim_{x\rightarrow x_{dom}}\left( 1-\frac{x}{x_{dom}}\right) \pi_1(x) \\
    &=&\lim_{x\rightarrow x_{dom}}\frac{-h_{2}(x,Y_0(x))
 \left(1-\frac{Y_0(x)}{ \widetilde{y}_0}\right) \pi_2(Y_0(x))-(1-\frac{Y_0(x)}{y_{dom}})
 h_{0}(x,Y_0(x))\pi_{0,0}}{\frac{1-\frac{Y_0(x)}{\widetilde{y}_0}}{1-\frac{x}{x_{dom}}}h_{1}(x,Y_0(x))} \\
    &=&\frac{-h_{2}(x_{dom},\widetilde{y}_0)\widetilde{L}(\widetilde{y}_0)\widetilde{y}_0}{h_{1}(x_{dom},\widetilde{y}_0)Y_0(x)^{\prime}(x_{dom})x_{dom}}.
\end{eqnarray*}

In the case of $x^{\ast}=\widetilde{x}_{1}=x_{3}=|x_{dom}|$, we
first have $\lim_{x\rightarrow x_{dom}} \left(
1-\frac{Y_0(x)}{\widetilde{y}_0} \right)
\pi_2(Y_{0}(x))=\widetilde{L}(\widetilde{y}_0)$. Then, using
equations (\ref{eqn:3.3}), (\ref{eqn:3.4}) and the expression for
$h_{1}(x_{3},\widetilde{y}_0)$, we obtain
\begin{eqnarray*}
    &&\lim_{x\rightarrow x_{dom}}\left( 1-\frac{x}{x_{dom}}\right) \pi_1(x) \\
    &=&\lim_{x\rightarrow x_{dom}}\frac{-h_{2}(x,Y_0(x))
\frac{\sqrt{1-\frac{x}{x_{dom}}}}{1-\frac{Y_0(x)}{\widetilde{y}_0}}
\left[\left( 1-\frac{Y_0(x)}{\widetilde{y}_0} \right)
\pi_2(Y_0(x))\right] -\sqrt{1-\frac{x}{x_{dom}}}
h_{0}(x,Y_0(x))\pi_{0,0}}{h_{1}(x,Y_0(x))/\sqrt{1-\frac{x}{x_{dom}}}} \\
    &=&\frac{\widetilde{L}(\widetilde{y}_0)h_{2}(x_{dom},\widetilde{y}_0)\widetilde{y}_0}{q_{1}(x_{dom})q(x_{dom})}.
\end{eqnarray*}

\textbf{Case 2.} If $x^{\ast}=x_{3,}$ then
$h_{1}(x_{dom},\widetilde{y}_0)=p_{1}(x_{dom})=0$, using equations
(\ref{eqn:1.9}), (\ref{eqn:3.1}), (\ref{eqn:3.3}) and
(\ref{eqn:3.4}), we can rewrite $\pi_1(x)$ as
\[
    \pi_1(x) =\frac{-h_{2}(x,Y_{0}(x))\pi_2(Y_{0}(x))-h_{0}(x,Y_{0}(x))\pi_{0,0}}{\sqrt{1-x/x_{dom}}\left[\sqrt{1-x/x_{dom}}p_{1}^{\ast}(x)+q_{1}(x) \right] }.
\]
It follows that
\begin{equation*}
    \lim_{x \to x_{dom}}\sqrt{1-x/x_{dom}}\pi_1(x)=\lim_{x \to x_{dom}}
\frac{-h_{2}(x,Y_{0}(x))
\pi_2(Y_{0}(x))-h_{0}(x,Y_{0}(x))\pi_{0,0}}{\left[
\sqrt{1-x/x_{dom}}p_{1}^{\ast}(x)+q_{1}(x)\right] }=c_{0,2}(x_{dom}).
\end{equation*}
Note that $q_{1}(x_{dom})\neq 0$.

Similarly, if $\widetilde{x}_{1}=x_{3}$, then $\widetilde{y}_0$ is
a pole of $\pi_2(y)$, which gives $\lim_{y\rightarrow
\widetilde{y}_0}\left( 1-\frac{y}{\widetilde{y}_0}\right)
\widetilde{ \pi}(y)=\widetilde{L}(\widetilde{y}_0)$.  Again, using
equations (\ref{eqn:1.9}), (\ref{eqn:3.1}), (\ref{eqn:3.3}) and
(\ref{eqn:3.4}), we obtain
\begin{eqnarray*}
    &&\lim_{x\rightarrow x_{dom}}\sqrt{1-\frac{x}{x_{dom}}}\pi_1(x) \\
    &=&\lim_{x\rightarrow
x_{dom}}\frac{-h_{2}(x,Y_0(x))\frac{\sqrt{1-\frac{x}{
x_{dom}}}}{1-\frac{Y_0(x)}{\widetilde{y}_0}}\left(
1-\frac{Y_0(x)}{\widetilde{y}_0}\right)
\pi_2(Y_0(x))-\sqrt{1-\frac{x}{x_{dom}}}h_{0}(x,Y_0(x))\pi_{0,0}}{
h_{1}(x,Y_{0}(x))} \\
    &=&-\frac{h_{2}(x_{dom},\widetilde{y}_0)\widetilde{L}(\widetilde{y}_0)}{
h_{1}(x_{dom},\widetilde{y}_0)}\lim_{x\rightarrow
x_{dom}}\frac{\widetilde{y}_0\sqrt{1-
\frac{x}{x_{dom}}}}{(1-x/x_{dom})p^{\ast}(x)-q(x)\sqrt{1-x/x_{dom}}} \\
    &=&\frac{h_{2}(x_{dom},\widetilde{y}_0)\widetilde{L}(\widetilde{y}_0)\widetilde{y}_0}{
h_{1}(x_{dom},\widetilde{y}_0)q(x_{dom})}.
\end{eqnarray*}

\textbf{Case 3.} Let
\begin{equation*}
    T(x,y)=\frac{h_{2}(x,y)\pi_2(y)+h_{0}(x,y)\pi_{0,0}}{-h_{1}(x,y)}.
\end{equation*}
Then,
\begin{equation*}
    \pi_1^{\prime}(x)=\frac{\partial T}{\partial x}+\frac{\partial T}{\partial y} \frac{dY_{0}(x)}{dx}
\end{equation*}
with
\begin{equation*}
    \frac{dY_{0}(x)}{dx}=p^{\prime}(x)+q^{\prime}(x)\sqrt{1-x/x_{dom}}-\frac{q(x)}{2x_{dom}\sqrt{1-x/x_{dom}}},
\end{equation*}
\begin{equation*}
    \frac{\partial T}{\partial x}=\frac{\widetilde{a}_{2}(y)\pi_2(y)
+\widetilde{a}_{0}(y)+[a_{1}^{\prime}(x)y+b_{1}^{\prime}(x)]T(x,y)}{-h_{1}(x,y)}
\end{equation*}
and
\begin{equation*}
    \frac{\partial T}{\partial y}=\frac{\frac{\partial h_{2}(x,y)}{\partial y}
\pi_2(y)+h_{2}(x,y)\pi_2^{\prime}(y)+\frac{\partial
h_{0}(x,y)\pi_{0,0}}{\partial y}+\frac{\partial
h_{1}(x,y)}{\partial y}T(x,y) }{-h_{1}(x,y)},
\end{equation*}
where $p(x)$ and $q(x)$ are defined by equation (\ref{eqn:3.2}).
Since $\lim_{x\rightarrow x_{dom}}\sqrt{1-x/x_{dom}}\frac{\partial
T}{\partial x}=0$, $\lim_{x\rightarrow x_{dom}}$
$\sqrt{1-x/x_{dom}}\frac{dY_{0}(x)}{dx}=-\frac{q(x_{dom})}{2x_{dom}}$
 and
$\frac{\partial T}{\partial y}$ is continuous at
$(x_{dom},Y_0(x_{dom}))$,
\begin{eqnarray}
    \lim_{x\rightarrow x_{dom}}\sqrt{1-x/x_{dom}}\pi_1^{\prime}(x)
    &=&-\frac{q(x_{dom})}{2x_{dom}}\frac{\partial T}{\partial y}|_{(x_{3},\widetilde{y}_0)} \\
    &=&-\frac{q(x_{dom})}{2x_{dom}}\frac{dT(x_{dom},y)}{dy}|_{y=\widetilde{y}_0}=c_{0,3}(x_{dom}).
\end{eqnarray}
It is easy to see $c_{3,0}(x_{dom}) \neq 0$, since otherwise
$\pi_1^{\prime}(x_{dom})<\infty $, which contradicts the fact that
$x_{3}$ is a branch point of $\pi_1(x)$.

%(For us to see: Let $U(x,y)=\frac{\partial h_{2}(x,y)}{\partial y}
%\pi_2(y)+h_{2}(x,y)\pi_2^{\prime}(y)+\frac{\partial
%h_{0}(x,y)\pi_{0,0}}{\partial y}+\frac{\partial
%h_{1}(x,y)}{\partial y} T(x,y) $. If $c_{3,0}=0$, then
%$U(x_{dom},Y_{0}(x_{dom}))=0$, which implies that, after some
%algebra manipulation, $\frac{\partial T}{\partial y}$ can written
%as
%\begin{equation*}
%\frac{\partial T}{\partial
%y}=\frac{U(x,y)-U(x_{dom},Y_{0}(x_{dom}))}{
%-h_{1}(x,y)}=\sqrt{1-x/x_{dom}}\Omega (x),
%\end{equation*}
%with $\lim_{x\rightarrow x_{3}}\Omega (x)<\infty $. It follows
%that $ \lim_{x\rightarrow x_{3}}\pi_1^{\prime}(x)<\infty $.)

\textbf{Case 4.} From equation (\ref{eqn:1.9}) and
(\ref{eqn:1.11}) in Theorem~\ref{theorem1.1}, we have
\[
\pi_1(x) =
\frac{h_{2}(x,Y_0(x))h_{1}(X_{0}(Y_0(x)),Y_0(x))\pi_1(X_{0}(Y_0(x)))
+ N(x)} {h_{1}(x,Y_0(x))h_{2}(X_{0}(Y_0(x)),Y_0(x))},
\]
where
\[
    N(x)=[h_{2}(x,Y_0(x))h_{0}(X_{0}(Y_0(x)),Y_0(x))-h_{2}(X_{0}(Y_0(x)),Y_0(x))h_{0}(x,Y_0(x))]\pi_{0,0}.
\]
 Since
\begin{equation*}
    \lim_{x\rightarrow x_{dom}}\frac{h_{1}(x,Y_0(x))}{x-x_{dom}} =\lim_{x\rightarrow x_{dom}}\frac{
h_{1}(x,Y_{0}(x))-h_{1}(x_{dom},Y_{0}(x_{dom}))}{x-x_{dom}}=h_{1}^{\prime}(x_{dom},\widetilde{y}_0)
\end{equation*}
and
\begin{eqnarray*}
    \lim_{x\rightarrow x_{dom}}\frac{h_{2}(X_{0}(Y_0(x)),Y_0(x))}{x-x_{dom}}
&=&\lim_{x\rightarrow x_{dom}}\frac{h_{2}(X_{0}(Y_0(x)),Y_0(x))-h_{2}(X_{0}(\widetilde{y}_0),\widetilde{y}_0)}{x-x_{dom}} \\
    &=&Y_{0}^{\prime}(x_{dom})h_{2}^{\prime}(X_{0}(\widetilde{y}_0),\widetilde{y}_0),
\end{eqnarray*}
we obtain
\begin{equation*}
    \lim_{x\rightarrow x_{dom}}\frac{\left( 1-\frac{x}{x_{dom}}\right )^{2}}{h_{1}(x,Y_0(x))h_{2}(X_{0}(Y_0(x)),Y_0(x))}
=\frac{1}{x_{dom}^{2}h_{1}^{\prime}(x_{dom},\widetilde{y}_0)Y_{0}^{\prime}(x_{dom})h_{2}^{\prime}(X_{0}(\widetilde{y}_0),\widetilde{y}_0)},
\end{equation*}
which yields
\begin{equation*}
    \lim_{x\rightarrow x_{dom}}\left( 1-\frac{x}{x_{dom}}\right )^{2}\pi_1(x)=c_{0,4}(x_{dom}).
\end{equation*}
\QED

\begin{remark}
It should be noted that the above theorem provides the asymptotic
behaviour at a dominant singularity, either positive or negative.
\end{remark}

\begin{corollary} \label{corollary3.1}
If
$h_{2}(x^{\ast},Y_0(x^{\ast}))\pi_2(Y_0(x^{\ast}))+h_{0}(x^{\ast},Y_0(y^{\ast}))\pi_{0,0}=
0$ or
$h_{1}(X_{0}(\widetilde{y}_0),\widetilde{y}_0)\pi(X_{0}(\widetilde{y}_0))+h_{0}(X_{0}(\widetilde{y}_0),$
$\widetilde{y}_0) \pi_{0,0}  = 0$, then the function $\pi_1(x)$,
as $x$ approaches to its dominant singularity, has one of the
three types of asymptotic properties shown in Case~1 to Case~3 of
Theorem~\ref{theorem3.1}.
\end{corollary}

\proof First suppose that $h_{2}(x^{\ast },Y_{0}(x^{\ast
}))\pi_{2}(Y_{0}(x^{\ast }))+h_{0}(x^{\ast },Y_{0}(x^{\ast
}))\pi_{0,0}=0$, but
$h_{1}(X_{0}(\widetilde{y}_{0}),\widetilde{y}_{0})\pi_{1}(X_{0}(\widetilde{y}_{0}))
+ h_{0}(X_{0}(\widetilde{y}_{0}),\widetilde{y}_{0})\pi_{0,0} \neq
0$. We then have the following four cases:

\textbf{1.} If $x^{\ast} < \widetilde{x}_{1}<x_{3}$, then
$\widetilde{x}_{1}$ is a pole and the dominant singular point of
$\pi_1(x)$ since $x^{\ast}$ is a removable singular point of
$\pi_1(x)$, which leads to $\lim_{x\rightarrow \widetilde{x}
}\left( 1-\frac{x}{\widetilde{x}_{1}}\right) \pi_1(x)=C$, the same
type in Case~1.

\textbf{2.} If $x^{\ast}<\widetilde{x}_{1}=x_3$, the same type of
asymptotic result as in Case~2 can be obtained.

\textbf{3.} If $x^{\ast }=x_{3}<\widetilde{x}_{1}$, then the
factor $\sqrt{1-x/x_{dom}}$ is cancelled out from both the
denominator and the numerator in the expression for $\pi _{1}(x)$.
By considering $\pi _{1}^{\prime }(x)$, we obtain the same type of
asymptotic result as that given in Case~3.

\textbf{4.} If $x^{\ast }=\widetilde{x}_{1}$, then $x^{\ast }$
would be a pole of $\pi _{2}(Y_{0}(x))$. This would imply
$\pi_{2}(Y_{0}(x^{\ast }))=\infty$, which contradict to
$h_{2}(x^{\ast },Y_{0}(x^{\ast }))
\pi_{2}(Y_{0}(x^{\ast}))+h_{0}(x^{\ast },Y_{0}(x^{\ast }))
\pi_{0,0}=0$ since $h_{2}(X_{0}(Y_{0}(x^{\ast }),Y_{0}(x^{\ast
}))=0$ implies $h_{2}(x^{\ast },Y_{0}(x^{\ast }))\neq 0$. Hence
this case is impossible.

Next, assume that both $h_{2}(x^{\ast },Y_{0}(x^{\ast }))
\pi_{2}(Y_{0}(x^{\ast})) + h_{0}(x^{\ast },Y_{0}(x^{\ast
}))\pi_{0,0}=0$ and
$h_{1}(X_{0}(\widetilde{y}_{0}),\widetilde{y}_{0})\pi_{1}(X_{0}(\widetilde{y}_{0}))
+ h_{0}(X_{0}(\widetilde{y}_{0}),\widetilde{y}_{0})\pi _{0,0}=0$.
We then have the following two cases:

\textbf{1.} If $\max \{\widetilde{x}_{1}$, $x^{\ast }\}<x_{3}$,
then both $\widetilde{x}_{1}$ and $x^{\ast }$ are removable poles
of $\pi _{1}(x)$. By considering $\pi_{1}^{\prime }(x)$, we obtain
the same type of asymptotic result as that given in Case~3.

\textbf{2.} If $\max \{\widetilde{x}_{1}$, $x^{\ast }\}=x_{3}$,
then the factor $\sqrt{1-x/x_{dom}}$ is cancelled out from both
the denominator and the numerator in the expression for
$\pi_{2}(Y_{0}(x))$ if $\widetilde{x}_{1}=x_{3}$ and for
$\pi_{1}(x)$ if $x^{\ast }=x_{3}$. By considering
$\pi_{1}^{\prime}(x)$, we obtain the same type of asymptotic
result as that given in Case~3.

Finally, the case in which
$h_{1}(X_{0}(\widetilde{y}_{0}),\widetilde{y}_{0})\pi_{1}(X_{0}(\widetilde{y}_{0}))
+ h_{0}(X_{0}(\widetilde{y}_{0}),\widetilde{y}_{0})\pi _{0,0}=0$,
but $h_{2}(x^{\ast },Y_{0}(x^{\ast })) \pi_{2}(Y_{0}(x^{\ast}))$
$+ h_{0}(x^{\ast },Y_{0}(x^{\ast }))\pi_{0,0} \neq 0$ can be
similarly considered. \QED

\begin{remark}
We believe that both
$h_{2}(x^{\ast},Y_0(x^{\ast}))\pi_2(Y_0(x^{\ast}))+h_{0}(x^{\ast},Y_0(x^{\ast}))\pi_{0,0}
\neq 0$ and $h_{1}(X_{0}(\widetilde{y}_0),\widetilde{y}_0)$
$\pi_1(X_{0}(\widetilde{y}_0))+h_{0}(X_{0}(\widetilde{y}_0),\widetilde{y}_0)\pi_{0,0}
\neq 0$ always hold, though at this moment we could not find a
proof. However, no new type of asymptotic property will appear
without this condition as shown in Corollary~\ref{corollary3.1}.
In the rest of the paper, the analysis will be carried out with
this condition, which is also valid without this condition.
\end{remark}

\begin{remark} \label{remark3.2}
When $x_{dom}=|x_{3}|<\min \{x^{\ast },\widetilde{x}_{1}\}$, the
numerator in the expression for $\pi_{1}(x)$ is not zero at
$x_{3}$.
\end{remark}

\section{Tail Asymptotics of Boundary Probabilities $\pi_{n,0}$ and $\pi_{0,n}$}

Since $\pi_1(x)$ and $\pi_2(y)$ are symmetric, properties for
$\pi_1(x)$ can be easily translated to the counterpart properties
for $\pi_2(y)$. Therefore, tail asymptotics for the boundary
probabilities $\pi_{0,n}$ can be directly obtained by symmetry.

 The exact tail
asymptotics of the boundary probabilities $\pi_{n,0}$ is a direct
consequence of Theorem~\ref{theorem3.1} and a Tauberian-like theorem applied to the function $\pi_1(x)$.
Specifically, if $\pi_1(x)$ has only one dominant singularity, then Theorem~\ref{tauberian-1}
is applied; and if $\pi_1(x)$ has two dominant singularity, then Theorem~\ref{tauberian-2}
is applied.

The following theorem shows that there are four types of exact
tail asymptotics, for large $n$, together with a possible periodic
property if $\pi_1(x)$ has two dominant singularities that have
the same asymptotic property.

In the theorem, let $x_{dom}$ be the positive dominant singularity of $\pi_1(x)$.
Consider the following four cases regarding which of $x^*$, $\widetilde{x}_1$ and $x_3$ will
be $x_{dom}$:
\begin{description}
\item[Case 1.] $x_{dom}=\min\{x^{\ast}, \widetilde{x}_{1}\}<x_{3}$
with $x^{\ast} \neq \widetilde{x}_{1}$, or
$x_{dom}=\widetilde{x}_{1}=x^{\ast}=x_{3}$;

\item[Case 2.] $x_{dom}=x_{3}=\min\{x^{\ast},\widetilde{x}_{1}\}$ with $x^{\ast}\neq
\widetilde{x}_{1}$;

\item[Case 3.] $x_{3}=x_{dom}<\min\{x^{\ast},\widetilde{x}_{1}\}$;

\item[Case 4.] $x_{dom}=x^{\ast}=\widetilde{x}_{1}<x_{3}$.
\end{description}

\begin{theorem}
\label{theorem4.1-a} Consider the stable non-singular genus 1
random walk. Corresponding to the above four cases, we have the
following tail asymptotic properties for the boundary
probabilities $\pi_{n,0}$ for large $n$. In all cases,
$c_{0,i}(x_{dom})$ ($1\leq i\leq 4$) are given in
Theorem~\ref{theorem3.1}.

\begin{description}
\item[1.] If $p_{i,j}$ is not X-shaped, then there are four types
of exact tail asymptotics:

\textbf{Case 1:} (Exact geometric decay)
\begin{equation} \label{eqn:exact}
    \pi_{n,0} \sim c_{0,1}(x_{dom})\left( \frac{1}{x_{dom}}\right)^{n-1};
\end{equation}

\textbf{Case 2:} (Geometric decay multiplied by a factor of
$n^{-1/2}$)
\begin{equation}
    \pi_{n,0} \sim \frac{c_{0,2}(x_{dom})}{\sqrt{\pi}}n^{-1/2}\left(\frac{1}{x_{dom}}\right)^{n-1};  \label{eqn:1/2}
\end{equation}

\textbf{Case 3:} (Geometric decay multiplied by a factor of
$n^{-3/2}$)
\begin{equation} \label{eqn:3/2}
    \pi_{n,0} \sim \frac{c_{0,3}(x_{dom})}{\sqrt{\pi}}n^{-3/2}\left(\frac{1}{x_{dom}}\right)^{n-1};
\end{equation}

\textbf{Case 4:} (Geometric decay multiplied by a factor of $n$)
\begin{equation}
    \pi_{n,0} \sim c_{0,4}(x_{dom})n\left( \frac{1}{x_{dom}}\right)^{n-1}; \label{eqn:n}
\end{equation}

\item[2.] If $p_{i,j}$ is X-shaped,
but both $p^{(1)}_{i,j}$ and $p^{(2)}_{i,j}$ are not X-shaped, we then have the following
exact tail asymptotic properties:

\textbf{Case 1:} (Exact geometric decay) It is given by (\ref{eqn:exact});

\textbf{Case 2:} (Geometric decay multiplied by a factor of
$n^{-1/2}$) It is given by (\ref{eqn:1/2});

\textbf{Case 3:} (Geometric decay multiplied by a factor of
$n^{-3/2}$)
\begin{equation} \label{eqn:p-3/2}
    \pi_{n,0} \sim \frac{\left[c_{0,3}(x_{dom})+(-1)^{n-1}c_{0,3}(-x_{dom})
\right] }{\sqrt{\pi}}n^{-3/2}\left( \frac{1}{x_{dom}}\right)^{n-1};
\end{equation}

\textbf{Case 4:} (Geometric decay multiplied by a factor of $n$)
It is given by (\ref{eqn:n});

\item[3.]  If $p_{i,j}$ and $p^{(1)}_{i,j}$ are X-shaped, but
$p^{(2)}_{i,j}$ is not, we then have the following
exact tail asymptotic properties:

\textbf{Case 1:} (Exact geometric decay) When $x^{\ast} \geq
\widetilde{x}_{1}$, it is given by (\ref{eqn:exact}); when
 $x_{dom}=x^{\ast} < \widetilde{x}_{1}$, it is given by
\begin{equation} \label{eqn:p-exact}
    \pi_{n,0} \sim \left[c_{0,1}(x_{dom})+(-1)^{n-1}c_{0,1}(-x_{dom})
\right] \left( \frac{1}{x_{dom}}\right)^{n-1};
\end{equation}

\textbf{Case 2:} (Geometric decay multiplied by a factor of
$n^{-1/2}$) When $x^{\ast}>\widetilde{x}_{1}$, it is given by (\ref{eqn:1/2}); when
 $x_{dom}=x^{\ast} < \widetilde{x}_{1}$, it is given by
 \begin{equation} \label{eqn:p-1/2}
    \pi_{n,0} \sim \frac{\left[c_{0,2}(x_{dom})+(-1)^{n-1}c_{0,2}(-x_{dom})
\right] }{\sqrt{\pi}}n^{-1/2}\left( \frac{1}{x_{dom}}\right)^{n-1};
\end{equation}

\textbf{Case 3:} (Geometric decay multiplied by a factor of
$n^{-3/2}$) It is given by (\ref{eqn:p-3/2}).

\textbf{Case 4:} (Geometric decay multiplied by a factor of $n$)
It is given by (\ref{eqn:n}).

\item[4.] If $p_{i,j}$ and $p^{(2)}_{i,j}$ are X-shaped, but
$p^{(1)}_{i,j}$ is not,  then it is the symmetric case to \textbf{3}. All
expression in \textbf{3} are valid after switching $x^{\ast}$ and
$\widetilde{x}_{1}$.

\item[5.] If all $p_{i,j}$, $p^{(1)}_{i,j}$ and $p^{(2)}_{i,j}$ are
X-shaped, we then have the following exact tail asymptotic properties:

\textbf{Case 1:} (Exact geometric decay) When $x^* \leq  \widetilde{x}_1$, it is given by (\ref{eqn:p-exact}); when $x^* > \widetilde{x}_1$, it is also given by (\ref{eqn:p-exact})
by replacing the dominant singularity $x^*$ by $\widetilde{x}_1$.

\textbf{Case 2:} (Geometric decay multiplied by a factor of $n^{-1/2}$)
 When $x^* < \widetilde{x}_1$, it is given by (\ref{eqn:p-1/2}); when $x^* > \widetilde{x}_1$, it is also given by (\ref{eqn:p-1/2})
by replacing the dominant singularity $x^*$ by $\widetilde{x}_1$.

\textbf{Case 3:} (Geometric decay multiplied by a factor of
$n^{-3/2}$) It is given by (\ref{eqn:p-3/2}).

\textbf{Case 4:} (Geometric decay multiplied by a factor of $n$) It is given by
\begin{equation} \label{eqn:p-n}
    \pi_{n,0} \sim \left[c_{0,4}(x_{dom})+(-1)^{n-1}c_{0,4}(-x_{dom})
\right] n\left( \frac{1}{x_{dom}}\right)^{n-1}.
\end{equation}

\end{description}
\end{theorem}

\proof \textbf{1.} Since $p_{i,j}$ is not X-shaped, all $-x_3$, $-x^*$ and $-\widetilde{x}_1$ are not
dominant singularities according to Corollary~\ref{corollary-X}, Theorem~\ref{theorem2.1} and Theorem~\ref{theorem2.1-b}.
Therefore, there is only one dominant singularity for $\pi_1(x)$. The tail asymptotic properties of $\pi_{n,0}$ follow from
Theorem~\ref{theorem3.1} and the direct application of the Tauberian-like theorem (Theorem~\ref{tauberian-1}).

 \textbf{2.} We only provide a proof to the cases, which are not identical to that in 1.
 \begin{description}
 \item[Case 1.] For the case that $x_{dom}=\widetilde{x}_{1}=x^{\ast}=x_{3}$, we notice that $-x_3$ is also a dominant singularity
 (Corollary~\ref{corollary-X}). In this case, the Tauberian-like theorem (Theorem~\ref{tauberian-2}) is used to have a tail asymptotic expression
 consisting of two terms, one, corresponding to the positive dominant singularity, with the exact geometric decay rate and the other, corresponding to
 the negative dominant singularity, with the geometric decay rate multiplied by a factor of $n^{-3/2}$. Therefore, the term with the geometric decay rate is the dominant (decay slower)
 term leading to the same tail asymptotic property given in (\ref{eqn:exact}).

 \item[Case 2.] Similar to Case~1, $-x_3$ is also a dominant singularity. The Tauberian-like theorem (Theorem~\ref{tauberian-2}) leads to a tail asymptotic expression
 consisting of two terms, one with the geometric rate  multiplied by a factor of $n^{-1/2}$ (dominant term) and the other by $n^{-3/2}$.

 \item[Case 3.] In this case, both $x_3$ and $-x_3$ are dominant singularities having the same asymptotic property according to Theorem~\ref{theorem3.1}.
 The tail asymptotic expression follows from the application of the Tauberian-like theorem (Theorem~\ref{tauberian-2}).
 \end{description}

 \textbf{3.} In this case, $-x_3$ and $-x^*$ are singularities, but $-\widetilde{x}_1$ is not. We only provide a proof to the cases, which are not identical to that in 1 or in 2.
 \begin{description}
 \item[Case 1.] For the case when $x^*=\widetilde{x}_1=x_3$, there are two dominant singularities. The Tauberian-like theorem (Theorem~\ref{tauberian-2}) leads
 to a tail asymptotic expression consisting of two terms, one (corresponding to the positive singularity) with a geometric decay rate, and the other (corresponding to
 the negative singularity) with the same geometric decay rate multiplied by a factor of $n^{-1/2}$ that is dominated by the geometric decay.

 When $x^* < \widetilde{x}_1$, both $x^*$ and $-x^*$ are dominant singularities with the same asymptotic property, which leads to the tail asymptotic expression by using
 Theorem~\ref{tauberian-2}.

 \item[Case 2.] For case when $x_3=x^*$, there are two dominant singularities having the same asymptotic property. The tail asymptotic expression follows from
 Theorem~\ref{tauberian-2}.

 \item[Case 4.] In this case, there are two dominant singularities, but the contribution from the positive dominant singularity dominates that from the negative
 dominant singularity. The tail asymptotic expression follows from Theorem~\ref{tauberian-2}.
 \end{description}

 \textbf{4.} The symmetric case to 3.

 \textbf{5.}  In this case, all $-x^*$, $-x^*$ and $-x_3$ are singularities. We only provide a proof to the cases, which are not considered in the above.
 \begin{description}
 \item[Case 1.] The only new situation here is the case when $x^*=x^*=x_3$. In this case, we have the same asymptotic property at both dominant singularities, which leads to
 (\ref{eqn:p-exact}).

 \item[Case 4.] In this case, we have the same asymptotic property at both dominant singularities, which leads to
 (\ref{eqn:p-n}).

 \end{description}

\QED

From the above theorem, it is clear that if there is only one dominant singularity,
then
the boundary probabilities $\pi_{n,0}$ have the following four
types of astymptotics: \textbf{1.} exact geometric; \textbf{2.}
geometric multiplied by a factor of $n^{-1/2}$; \textbf{3.}
geometric multiplied by a factor of $n^{-3/2}$; and \textbf{4.}
geometric multiplied by a factor of $n$.
If there are two dominant singularities, but with different asymptotic properties,
$\pi_{n,0}$ also has one of the above four types of tail asymptotic properties.
Finally, if we have the same asymptotic property at both dominant singularities, then $\pi_{n,0}$
reveals a periodic property with the above four types of tail asymptotics, which is a new discovery.

\section{Tail Asymptotics of the Marginal Distributions}

In the previous section, we have seen that the asymptotic
behaviour of the function $\pi_1(x)$ ($\pi_2(y)$) at its dominant
singularity or singularities determines the tail asymptotic property of the
boundary probabilities $\pi_{n,0}$ ($\pi_{0,n}$). According the
the fundamental form of the random walk, it, together with the
property of the kernel function $h(x,y)$, also determines the
tail asymptotic property of the marginal distribution
$\pi_{n}^{(1)} = \sum_j \pi_{n,j}$ (and $\pi_{n}^{(2)} = \sum_i
\pi_{i,n}$).

In this section, we provide details for the exact tail asymptotics of
the marginal distribution $\pi_{n}^{(1)}$. The exact tail asymptotics of $\pi_{n}^{(2)}$
can be easily obtained by symmetry. First, based on the
fundamental form, we have
\[
    \pi(x,y)=\frac{h_{1}(x,y)\pi_1(x)+h_{2}(x,y)\pi_2(y)+h_{0}(x,y)\pi_{0,0}}{-h(x,y)}
\]
and therefore,
\begin{eqnarray*}
    \pi(x,1) &=&\frac{h_{1}(x,1)\pi_1(x)+h_{2}(x,1)\pi_2(1)+h_{0}(x,1)\pi_{0,0}}{-h(x,1)} \\
&=&\frac{h_{1}(x,1)\pi_1(x)+h_{2}(x,1)\pi_2(1)+h_{0}(x,1)\pi_{0,0}}{-\widetilde{a}(1)[x-X_{0}(1)][x-X_{1}(1)]}.
\end{eqnarray*}

If $M_{x} \geq 0$, then $X_{1}(1)=1$, which implies that the
denominator of the expression for $\pi(x,1)$ does not have any
zero outside the unit circle. In this case, $\pi_{n}^{(1)}$ has
the same tail asymptotics as $\pi_{n,0}$. The only difference is
the expression for the coefficient, which can be obtained from straight
forward calculations.

If $M_{x}<0$, then $X_{0}(1)=1$ and $X_{1}(1)>1$. If $p_{i,j}$ is
not X-shaped, the analysis is so-called standard, details of which
will be provided here. If $p_{i,j}$ is X-shaped, then there are
four subcases based on if $p^{(k)}_{i,j}$ is X-shaped or not. For
these cases, detailed analysis varies, but similar. We provide
details here for the case where both $p^{(k)}_{i,j}$ for $k=1, 2$
are not X-shaped. Let
 $z=\min \{x^{\ast},\widetilde{x}_{1}\}$.
and consider the following four cases:

\textbf{1.}  $\min (X_{1}(1),z)<x_{3}$ and $X_{1}(1)\neq z$. In
this case, $\pi_{n}^{(1)}$ has an exact geometric decay with the
decay rate equal to $x_{dom}=\min (X_{1}(1),z)$:
\[
    \pi_{n}^{(1)}\sim c_{1}^{(x)} \left(
\frac{1}{x_{dom}}\right)^{n-1},
\]
where
\[
c_{1}^{(x)}=\left\{
\begin{array}{ll}
\frac{\lbrack h_{1}(X_{1}(1),1)\pi _{1}(X_{1}(1))+h_{2}(X_{1}(1),1)\pi
_{2}(1)+h_{0}(X_{1}(1),1)\pi _{0,0}]X_{1}(1)}{\widetilde{a}(1)(X_{1}(1)-1)},
& X_{1}(1)<z, \\
\frac{h_{1}(z,1)c_{0,1}(z)+h_{2}(z,1)\pi _{2}(1)+h_{0}(z,1)\pi _{0,0}}{-%
\widetilde{a}(1)[z-X_{0}(1)][z-X_{1}(1)]}, & X_{1}(1)>z,
\end{array}
\right.
\]
with $c_{0,1}(z)$ being given in Theorem~\ref{theorem3.1}.

\textbf{2.} $X_{1}(1)=z<x_{3}$. In this case, $X_{1}(1)=
\widetilde{x}_{1}$ is impossible, since otherwise
$h(X_{1}(1),1)=0$, which implies $1=Y_{0}(\widetilde{x}_{1})$ or
$1=Y_{1}(\widetilde{x}_{1})$. This is contradiction to
$Y_{0}(\widetilde{x}_{1})>1$. Hence, only $X_{1}(1)=x^{\ast}$ may
hold. There are two subcases:

\textbf{2(a):} $1=Y_{0}(x^{\ast})$. In this case,
$h_{1}(x^{\ast},1)=$ $h_{1}(x^{\ast},Y_{0}(x^{\ast}))=0$. We can
write $h_{1}(x,1)$ as
$h_{1}(x,1)=a_{1}(x)+b_{1}(x)=(x-X_{1}(1))h_{1}^{\ast}(x)$ with
$\frac{ h_{1}(x,1)}{x-X_{1}(1)}$ being a linear function of $x$,
which yields
\begin{eqnarray*}
\pi(x,1) &=&\frac{h_{1}(x,1)\pi_1(x)+h_{2}(x,1)\pi_2
(1)+h_{0}(x,1)\pi_{0,0}}{-\widetilde{a}(1)(x-1)[x-X_{1}(1)]} \\
&=&\frac{h_{1}^{\ast}(x)\pi_1(x)}{-\widetilde{a}(1)(x-1)}+\frac{h_{2}(x,1)
\pi_2(1)+h_{0}(x,1)\pi_{0,0}}{-\widetilde{a}(1)(x-1)[x-X_{1}(1)]}.
\end{eqnarray*}
Therefore, $\pi(x,1)$ has a single pole $X_{1}(1)$, which leads to
an exact geometric decay (recalling $\pi(1,1)\neq 1$):
\[
    \pi_{n}^{(1)}\sim c_{2,1}^{(x)}\left( \frac{1}{X_{1}(1)}\right)^{n-1}
\]
 with the
coefficient given by
\[
    c_{2,1}^{(x)}=\lim_{x\rightarrow X_{1}(1)}\left( 1-\frac{x}{X_{1}(1)}\right)
\pi(x,1)=\frac{[h_{2}(X_{1}(1),1)\pi_2(1)+h_{0}(X_{1}(1),1)\pi_{0,0}]X_{1}(1)}{\widetilde{a}(1)(X_{1}(1)-1)}.
\]

\textbf{2(b):} $1=Y_{1}(x^{\ast})$. In this case,
$h_{1}(x^{\ast},Y_{0}(x^{\ast}))=0$ and
$Y_{0}(x^{\ast})<Y_{1}(x^{\ast})$, we obtain $ h_{1}(x^{\ast},1)=$
$h_{1}(x^{\ast},Y_{1}(x^{\ast}))>0$, which implies that $x^{\ast}$
is a double pole of $\pi(x,1)$ (noting that
$h_{2}(x^{\ast},1)\pi_2(1)+h_{0}(x^{\ast},1)\pi_{0,0}>0$ since
$h_{2}(x^{\ast},1)>h_{2}(X_{0}(1),1)=0$ and
$h_{0}(x^{\ast},1)\pi_{0,0}>0$). The corresponding tail asymptotic
is given by
\[
\pi_{n}^{(1)}\sim c_{2,2}^{(x)}n\left(
\frac{1}{X_{1}(1)}\right)^{n-1},
\]
where
\begin{eqnarray*}
    c_{2,2}^{(x)}&=& \lim_{x\rightarrow X_{1}(1)}\left(1-\frac{x}{X_{1}(1)}\right
    )^{2}\pi(x,1)\\
    &=&\frac{X_{1}(1)[h_{1}(X_{1}(1),1)c_{0,1}(x^{\ast})+h_{2}(X_{1}(1),1)
\pi_2(1)+h_{0}(X_{1}(1),1)\pi_{0,0}]}{\widetilde{a}
(1)[X_{1}(1)-1]}
\end{eqnarray*}
with $c_{0,1}(x^{\ast})$ given in Theorem~\ref{theorem3.1}.
\medskip

\textbf{3.} $\min (X_{1}(1),z)=x_{3}$. In this case, there are
four possible subcases, for which proofs are omitted since they
are similar to that for the previous cases:

\textbf{3(a):} $X_{1}(1)=z=x_{3}$ leading to an exact geometric
decay:
\[
\pi_{n}^{(1)}\sim c_{3,1}^{(x)}\left( \frac{1}{x_{3}}\right
)^{n-1},
\]
where
\[
    c_{3,1}^{(x)} = \lim_{x\rightarrow x_{3}} \left (1-\frac{x}{x_{3}}\right )
\pi(x,1)=\frac{h_{2}(x_{3},1)\pi_2(1)+h_{0}(x_{3},1)\pi_{0,0}}{x_{3}\widetilde{a}(1)(x_{3}-1)}.
\]

\textbf{3(b):} $X_{1}(1)=x_{3}<z$ leading to an exact geometric
decay:
\[
    \pi_{n}^{(1)}\sim c_{3,2}^{(x)}\left(\frac{1}{x_{3}}\right)^{n-1},
\]
where
\[
    c_{3,2}^{(x)}= \lim_{x\rightarrow x_{3}}\left(
1-\frac{x}{x_{3}}\right) \pi(x,1) =
\frac{h_{1}(x_{3},1)\pi(x_{3})+h_{2}(x_{3},1)\pi_2
(1)+h_{0}(x_{3},1)\pi_{0,0}}{x_{3}\widetilde{a}(1)(x_{3}-1)}.
\]

\textbf{3(c):} $z=x_{3}<X_{1}(1)$ with $x^{\ast }\neq
\widetilde{x}_{1}$ leading to a geometric decay multiplied by the
factor $n^{-1/2}$:
\[
    \pi_{n}^{(1)}\sim c_{3,3}^{(x)}n^{-1/2}\left(
\frac{1}{x_{3}}\right)^{n-1},
\]
where
\[
    c_{3,3}^{(x)}=\lim_{x\rightarrow x_{3}}\left(
1-\frac{x}{x_{3}}\right)^{1/2}\pi(x,1)=\frac{
h_{1}(x_{3},1)c_{0,2}(x_3)}{\widetilde{a}(1)(X_{1}(1)-1)[X_{1}(1)-x_{3}]}
\]
with $c_{0,2}(x_3)$ given in Theorem~\ref{theorem3.1}.

\textbf{3(d):} $z=x^{\ast }=\widetilde{x}_{1}=x_{3}<X_{1}(1)$
leading to an exact geometric decay
\[
\pi _{n}^{(1)}\sim c_{3,4}^{(x)}\left( \frac{1}{x_{3}}\right) ^{n-1},
\]
where
\[
c_{3,4}^{(x)}=\lim_{x\rightarrow z}\left( 1-\frac{x}{z}\right) \pi (x,1)=
\frac{h_{1}(z,1)c_{0,1}(z)+h_{2}(z,1)\pi _{2}(1)+h_{0}(z,1)\pi _{0,0}}{-\widetilde{a}(1)[z-X_{0}(1)][z-X_{1}(1)]}.
\]

\textbf{4.} $x_{3}<\min (z,X_{1}(1))$ leading to a geometric decay
multiplied by the factor $ n^{-3/2}$:
\[
\pi_{n}^{(1)}\sim c_{4}^{(x)}n^{-3/2}\left(
\frac{1}{x_{3}}\right)^{n-1},
\]
where
\[
    c_{4}^{(x)}=\lim_{x\rightarrow x_{3}}\left(1-\frac{x}{x_{3}} \right)^{1/2}\pi^{\prime}(x,1)=
\frac{h_{1}(x_{3},1)c_{0,3}(x_{3})}{\widetilde{a}(1)(x_{3}-1)[X_{1}(1)-x_{3}]
}
\]
with $c_{0,3}(x_{3})$ given in Theorem~\ref{theorem3.1}.

For the completeness, we provide a summary of tail asymptotic
properties for the marginal distribution $\pi_{n}^{(1)}$ for all
possible cases. For this purpose, let $x_{dom}$ be the positive
dominant singularity of $\pi(x,1)$. Note that $X_{1}(1)\neq
\widetilde{x}_{1}$. The following are the all possible cases
according to which of $\widetilde{x}_{1}$, $x^{\ast }$, $x_{3}$
and $X_{1}(1)$ is $x_{dom}$.

\textbf{Case A.}  $x_{dom}=\min
\{\widetilde{x}_{1},x^{\ast},x_{3}\}<X_{1}(1)$;

\textbf{Case B.} $x_{dom}=X_{1}(1)<\min
\{\widetilde{x}_{1},x^{\ast },x_{3}\}$;

\textbf{Case C.} $x_{dom}=X_{1}(1)=x^{\ast }<\min
\{\widetilde{x}_{1},x_{3}\}$;

\textbf{Case D.} $x_{dom}=X_{1}(1)=x_{3}<x^{\ast }$;

\textbf{Case E.} $x_{dom}=X_{1}(1)=x_{3}=x^{\ast }$.
\bigskip

\begin{remark}
The cases here are different from the cases classified in the
previous section and the next section.
\end{remark}

The exact tail asymptotic properties are obtained according to the
expression of $\pi(x,1)$ and the Taubarian-like theorem.

\begin{theorem}
For the stable non-singular genus 1 random walk, the exact tail
asymptotic properties for the marginal distribution
$\pi_{n}^{(1)}$, as $n$ is large, are summarized as:

\begin{description}
\item[Case~A:]  This case includes Cases~1--4 in the previous
section. $\pi_{n}^{(1)}$ has the same types of asymptotic
properties as $\pi_{n,0}$ given in Theorem~\ref{theorem4.1-a},
respectively, with possible different expressions for the
coefficients.

\item[Case~B:] $\pi_{n}^{(1)}$ has an exact geometric decay.

\item[Case~C:] $\pi_{n}^{(1)}$ has an exact geometric decay if
$Y_{0}(x^{\ast })=1$ and a geometric decay multiplied by a factor
of $n$ if $Y_{1}(x^{\ast })=1$, respectively.

\item[Case~D:] $\pi_{n}^{(1)}$ has an exact geometric decay.

\item[Case~E:] $\pi_{n}^{(1)}$ has an exact geometric decay.

\end{description}
\end{theorem}

\section{Tail Asymptotics for Joint Probabilities}
\label{section4}

In the previous sections, we have seen how we can derive exact
tail asymptotic properties for the boundary probabilities and for
the marginal distributions based on the asymptotic property of
$\pi_1(x)$ ($\pi_2(y)$) and the kernel function. However, the
exact tail asymptotic behaviour for joint probabilities cannot be
obtained directly from them. Further tools are needed for this
purpose. Our goal is to characterize the exact tail asymptotics
for $\pi_{n,j}$ for each fixed $j$ and $\pi_{i,n}$ for each fixed
$i$. Due to the symmetry, in this section, we provide details only
for the former.

The relevant balance equations of the random walk are given by
\begin{eqnarray*}
(1-p_{0,0}^{(0)})\pi_{0,0} &=&p_{-1,0}^{(1)}\pi_{1,0}+p_{0,-1}^{(2)}\pi_{0,1}+p_{-1,-1}\pi_{1,1}, \\
(1-p_{0,0}^{(1)})\pi_{1,0} &=&p_{1,0}^{(0)}\pi_{0,0}+p_{-1,0}^{(1)}\pi_{2,0}+p_{-1,-1}\pi_{2,1}+p_{1,-1}^{(2)}\pi_{0,1}+p_{0,-1}\pi_{1,1}, \\
(1-p_{0,0}^{(1)})\pi_{i,0} &=&p_{1,0}^{(1)}\pi_{i-1,0}+p_{-1,0}^{(1)}\pi_{i+1,0}+p_{-1,-1}\pi_{i+1,1}+p_{1,-1}\pi_{i-1,1}+p_{0,-1}\pi_{i,1},\;\; i\geq 2, \\
(1-p_{0,0})\pi_{i,j} &=&p_{1,-1}\pi_{i-1,j+1}+p_{-1,-1}\pi_{i+1,j+1}+p_{0,-1}\pi_{i,j+1}+p_{1,0}\pi_{i-1,j}+p_{-1,0}\pi_{i+1,j} \\
&&+p_{1,1}\pi_{i-1,j-1}+p_{0,1}\pi_{i,j-1}+p_{-1,1}\pi_{i+1,j-1},\;\;
j\geq 2.
\end{eqnarray*}
Let
\begin{eqnarray*}
    \varphi_{j}(x) &=&\sum_{i=1}^{\infty }\pi_{i,j}x^{i-1},\;\;\; j\geq 0, \\
    \psi_{i}(y) &=&\sum_{j=1}^{\infty }\pi_{i,j}y^{i-1},\;\;\; i\geq 0.
\end{eqnarray*}
From the above definition, it is clear that $\varphi_{0}(x)=\pi_1(x)$ and $\psi_0(y)=\pi_2(y)$.
From the relevant balance equations, we obtain
\begin{eqnarray}
    c(x)\varphi_{1}(x)+b_{1}(x)\varphi_{0}(x) &=&a_{0}^{\ast}(x), \label{eqn:4.1} \\
    c(x)\varphi_{2}(x)+b(x)\varphi_{1}(x)+a_{1}(x)\varphi_{0}(x) &=&a_{1}^{\ast}(x), \label{eqn:4.2} \\
    c(x)\varphi_{j+1}(x)+b(x)\varphi_{j}(x)+a(x)\varphi_{j-1}(x) &=&a_{j}^{\ast}(x),\;\;\; j\geq 2, \label{eqn:4.3}
\end{eqnarray}
or
\begin{equation} \label{eqn:4.4}
    \varphi_{j+1}(x)=\frac{-b(x)\varphi_{j}(x)-a(x)\varphi_{j-1}(x)+a_{j}^{\ast}(x)}{c(x)},\;\;\; j \geq 0,
\end{equation}
where
\begin{eqnarray*}
    a_{0}^{\ast}(x) &=&-c_{2}(x)\pi_{0,1}-b_{0}(x)\pi_{0,0}, \\
    a_{1}^{\ast}(x) &=&-c_{2}(x)\pi_{0,2}-b_{2}(x)\pi_{0,1}-a_{0}(x)\pi_{0,0}, \\
    a_{j}^{\ast}(x) &=&-c_{2}(x)\pi_{0,j+1}-b_{2}(x)\pi_{0,j}-a_{2}(x)\pi_{0,j-1},\;\;\; j \geq 2.
\end{eqnarray*}

First, we establish the fact that a zero of $c(x)$ is not a pole of $\varphi_{j}(x)$ for all $j\geq 0$. Therefore $\varphi_{j}(x)$ has the same
singularities as $\varphi_{0}(x)$.

Let $y=Y_{0}(x)$ be in the cut plane $\widetilde{\mathbb{C}}_{x}$,
and let $y_{dom}$ and $x_{dom}$ be the positive dominating
singular points of $\psi_{0}(y)$ and $\varphi_{0}(x)$,
respectively. Let
\[
f_{k}(x)=-a_{2}(x)\sum_{j=k-1}^{\infty
}\pi_{0,j}y^{j-(k-1)}-b_{2}(x)\sum_{j=k}^{\infty
}\pi_{0,j}y^{j-k}-c_{2}(x)\sum_{j=k+1}^{\infty
}\pi_{0,j}y^{j-(k+1)},\;\;\; k\geq 1,
\]
then,
\begin{align}
    f_{1}(x) & = y f_{2}(x)-c_{2}(x)\pi _{0,2}-b_{2}(x)\pi _{0,1}, \nonumber \\
        \label{eqn:4.5}
    f_{k}(x) & =yf_{k+1}(x)+a_{k}^{\ast}(x), \quad k \geq 2.
\end{align}
According to Theorem~\ref{theorem1.1}, when $|x|<x_{dom}$, we obtain
\begin{align}
    h_{1}(x,y)\varphi_{0}(x) &= -h_{2}(x,y)\psi_{0}(y)-h_{0}(x,y)\pi_{0,0}
\nonumber \\
&= y[f_{1}(x)-a_{0}(x)\pi _{0,0}]+a_{0}^{\ast }(x) \label{eqn:4.6} \\
&= y^{2}f_{2}(x)+ya_{1}^{\ast}(x)+a_{0}^{\ast}(x) \label{eqn:4.7} \\
&= y^{3}f_{3}(x)+y^{2}a_{2}^{\ast}(x)+ya_{1}^{\ast}(x)+a_{0}^{\ast}(x). \label{eqn:4.8}
\end{align}
Let $u=\frac{y}{c(x)}=\frac{Y_{0}(x)}{c(x)}$. Since a zero of $c(x)$ is a
zero of $Y_{0}(x)$, $u$ is analytic on the cut plane $\widetilde{C}_{x}$.
Using $\frac{b(x)}{a(x)}=-Y_{1}(x)-Y_{0}(x)$ and $\frac{c(x)}{a(x)}
=Y_{1}(x)Y_{0}(x)$, we obtain
\begin{equation} \label{eqn:4.9}
1+b(x)u=-yua(x)\text{ and }\frac{1+b(x)u}{c(x)}=-a(x)u^{2}.
\end{equation}

Write $a=a(x)$, $b=b(x)$, $c=c(x)$, $a_{i}=a_{i}(x)$, $b_{i}=b_{i}(x)$,
$a_{j}^{\ast}=a_{j}^{\ast}(x)$, $f_{j}=f_{j}(x)$ and $\varphi_{j}=\varphi_{j}(x)$. We have following Lemma, which
confirms that a zero of $c(x)$ is not a pole of $\varphi_{j}(x)$ for all $j\geq 0$. Therefore,
$\varphi_{j}(x)$ has the same singularities as $\varphi_{0}(x)$.

\begin{lemma} \label{lemma4.1}
Let
\begin{equation}
    w_{-1} =-\frac{a_{1}}{au},\;\;\;
    w_{0} =-b_{1},\;\;\;  \label{eqn:4.10}
    w_{j} =buw_{j-1}+(1+bu)w_{j-2}.
\end{equation}
Then,
\begin{equation} \label{eqn:4.11}
    (-1)^{j}\left[ buw_{j-1}+(1+bu)w_{j-2}\right]+b_{1}+a_{1}y=(-1)^{j}(1+bu)w_{j-1},\;\;\; j\geq 1,
\end{equation}
and
\begin{eqnarray}
    h_{1}\varphi_{1} &=&yuf_{2}w_{0}+ug_{1}, \label{eqn:4.12} \\
    h_{1}\varphi_{j} &=&(-1)^{j+1}yuf_{j+1}w_{j-1}+u\sum_{k=0}^{j-2}(-1)^{j+1-k}a_{j-k}^{\ast}w_{j-1-k}(au)^{k}+ug_{1}(au)^{j-1},\;\;\; j\geq 2, \label{eqn:4.13}
\end{eqnarray}
where $g_{1}=a_{0}^{\ast}(x)a_{1}(x)-b_{1}(x)a_{1}^{\ast}(x)$.
\end{lemma}

\proof By applying (\ref{eqn:4.9}) and (\ref{eqn:4.10}), we easily obtain equation
(\ref{eqn:4.11}) for $j=1$. Assume that equation (\ref{eqn:4.11}) is true for $j\leq k$, we show
\begin{equation} \label{eqn:4.14}
    (-1)^{k+1}\left[ buw_{k}+(1+bu)w_{k-1}\right ] +b_{1}+a_{1}y=(-1)^{k+1}(1+bu)w_{k}.
\end{equation}
From the inductive assumption and the definition of $w_{j}$, we have
\begin{eqnarray} \label{eqn:4.15}
    b_{1}+a_{1}y &=&(-1)^{k}(1+bu)w_{k-1}-(-1)^{k}\left[ buw_{k-1}+(1+bu)w_{k-2} \right ],  \nonumber \\
    &=&(-1)^{k}(1+bu)w_{k-1}+(-1)^{k+1}w_{k},
\end{eqnarray}
which yields equation (\ref{eqn:4.14}). Equation (\ref{eqn:4.12}) is obtained by the direct
substitutions of equations (\ref{eqn:4.1}) and (\ref{eqn:4.6}).

Next, we show equation (\ref{eqn:4.13}). We use the induction again.
According to
equations (\ref{eqn:4.2}), (\ref{eqn:4.7}) and (\ref{eqn:4.12}),
\begin{eqnarray*}
    c(x)h_{1}\varphi_{2} &=&-bh_{1}\varphi_{1}-a_{1}h_{1}\varphi_{0}+h_{1}a_{1}^{\ast} \\
    &=&-b[yuf_{2}w_{0}+ug_{1}]-a_{1}[y^{2}f_{2}+ya_{1}^{\ast}+a_{0}^{\ast}]+a_{1}^{\ast}[a_{1}y+b_{1}].
\end{eqnarray*}
It follows from equations (\ref{eqn:4.9}) and (\ref{eqn:4.10}) that
\begin{eqnarray*}
    h_{1}\varphi_{2}&=&-bu^{2}f_{2}w_{0}-a_{1}yuf_{2}+u(au)g_{1}=-uf_{2} \left [buw_{0}+a_{1}\frac{1+bu}{-au} \right ]+u(au)g_{1} \\
    &=&-uf_{2}[buw_{0}+(1+bu)w_{-1}]+u(au)g_{1}=-yuf_{3}w_{1}-a_{2}^{\ast}uw_{1}+u(au)g_{1},
\end{eqnarray*}
which gives equation (\ref{eqn:4.13}) for $j=2$. Assume that equation (\ref{eqn:4.13}) is true for
$j\leq n$. We prove the result for $j=n+1$. From equations (\ref{eqn:4.9}), (\ref{eqn:4.11}), (\ref{eqn:4.5}) and the
inductive assumption, we have
\begin{eqnarray*}
  &&   c(x)h_{1}\varphi_{n+1}=-bh_{1}\varphi_{n}-ah_{1}\varphi_{n-1}+h_{1}a_{n}^{\ast} \\
&=&(-1)^{n+2}ybuf_{n+1}w_{n-1}+bu\sum_{k=0}^{n-2}(-1)^{n+2-k}a_{n-k}^{\ast}w_{n-1-k}(au)^{k}-bug_{1}(au)^{n-1} \\
&&+(-1)^{n+1}yauf_{n}w_{n-2}+au\sum_{k=0}^{n-3}(-1)^{n+1-k}a_{n-1-k}^{\ast}w_{n-2-k}(au)^{k}-g_{1}(au)^{n-1}+a_{n}^{\ast}[a_{1}y+b_{1}] \\
&=&(-1)^{n+2}yf_{n+1}[buw_{n-1}-yauw_{n-2}]+a_{n}^{\ast}\left\{
(-1)^{n+2}buw_{n-1}+(-1)^{n+1}yauw_{n-2}+a_{1}y+b_{1}\right\} \\
&&+bu\sum_{k=1}^{n-2}(-1)^{n+2-k}a_{n-k}^{\ast}w_{n-1-k}(au)^{k}+\sum_{k=0}^{n-3}(-1)^{n+1-k}a_{n-1-k}^{\ast}w_{n-2-k}(au)^{k+1}-(1+bu)g_{1}(au)^{n-1} \\
&=&(-1)^{n+2}yf_{n+1}[buw_{n-1}+(1+bu)w_{n-2}]+(-1)^{n+2}a_{n}^{\ast}(1+bu)w_{n-1} \\
&&+(1+bu)\sum_{k=1}^{n-2}(-1)^{n+2-k}a_{n-k}^{\ast}w_{n-1-k}(au)^{k}-(1+bu)g_{1}(au)^{n-1},
\end{eqnarray*}
which yields
\begin{eqnarray*}
h_{1}\varphi_{n+1} &=&(-1)^{n+2}uf_{n+1}w_{n}+(-1)^{n+1}a_{n}^{\ast}au^{2}w_{n-1}+au^{2}\sum_{k=1}^{n-2}(-1)^{n+1-k}a_{n-k}^{\ast}w_{n-1-k}(au)^{k}+ug_{1}(au)^{n} \\
&=&(-1)^{n+2}yuf_{n+1}w_{n}+u\sum_{k=0}^{n-1}(-1)^{n+2-k}a_{n+1-k}^{\ast}w_{n-k}(au)^{k}+ug_{1}(au)^{n}.
\end{eqnarray*}
This completes the proof. \QED

\begin{corollary}
$\varphi_{0}(x)$ and $\varphi_{j}(x)$, $j \geq 1$, have the same
singularities.
\end{corollary}

The following Lemma is useful in characterizing the tail asymptotics of $\pi_{n,j}$ for a fixed $j$.
\begin{lemma} \label{lemma4.2}
If $\min \{x^{\ast},\widetilde{x}_{1}\}>x_{3}$, then
\[
    \lim_{x\rightarrow x_{dom}}\sqrt{1-\frac{x}{x_{dom}}}\varphi_{j}^{\prime}(x)=c_{3,j}(x_{dom}),
\]
where $c_{3,0}(x_{dom})$ is given in Theorem~\ref{theorem3.1} and
\begin{equation} \label{eqn:4.16}
    c_{3,j+1}(x_{dom})=[A_{3}(x_{dom})+B_{3}(x_{dom})j]\left(\frac{1}{ Y_{1}(x_{dom})}\right)^{j},\;\;\; j\geq 0,
\end{equation}
with
\begin{align} \label{eqn:4.17}
    A_{3}(x_{dom})
    &=-\frac{c_{3,0}(x_{dom})b_{1}(x_{dom})}{c(x_{dom})}, \\
    B_{3}(x_{dom})
    &=\frac{-h_{1}(x_{dom},Y_{0}(x_{dom}))c_{3,0}(x_{dom})}{c(x_{dom})}. \label{eqn:4.17-b}
\end{align}
\end{lemma}

\proof  When $\min \{x^{\ast},\widetilde{x}_{1}\}>x_{3}$, we have
$x_{dom}=\pm x_{3}$. Without lose of
generality, we assume $x_{dom}=x_{3}$ in the proof.
Since $\varphi_{j}(x)$, $j\geq 0$, is continuous at $x_{3}$, $\lim_{x\rightarrow x_{3}}\sqrt{1-\frac{x}{x_{3}}}\varphi_{j}(x)=0$.
Let $j=1$. Then,
\[
\varphi_{1}^{\prime}(x)=\frac{-c^{\prime}(x)\varphi_{1}(x)-b_{1}(x)\varphi_{0}^{\prime}(x)-b_{1}(x)\varphi_{0}^{\prime}(x)+a_{0}^{\ast\prime}(x)}{c(x)}
\]
and
\[
    \lim_{x\rightarrow x_{3}}\sqrt{1-\frac{x}{x_{3}}}\varphi_{1}^{\prime}(x)=\frac{-b_{1}(x_{3})c_{3,0}(x_{3})}{c(x_{3})}=c_{3,1}(x_{3}).
\]
Assume that $\lim_{x\rightarrow x_{3}}\sqrt{1-\frac{x}{x_{3}}}\varphi_{k}^{\prime}(x)$ exists for $k \leq j$ and
\[
    \lim_{x\rightarrow x_{3}}\sqrt{1-\frac{x}{x_{3}}}\varphi_{k}^{\prime}(x)=c_{3,k}(x_{3}),
\]
we obtain
\[
    \varphi_{k+1}^{\prime}(x) =\frac{-c^{\prime}(x)\varphi_{k+1}(x)-b(x)\varphi_{k}^{\prime}(x)-b^{\prime}(x)\varphi_{k}(x)-a(x)\varphi_{k-1}^{\prime}(x)-a^{\prime}(x)\varphi_{k-1}(x)+a_{k}^{\ast\prime}(x)}{c(x)},
\]
and
\[
\lim_{x\rightarrow
x_{3}}\sqrt{1-\frac{x}{x_{3}}}\varphi_{k+1}^{\prime}(x)=
\frac{-b(x_{3})c_{3,k}(x_{3})-a(x_{3})c_{3,k-1}(x_{3})}{c(x_{3})}=c_{3,k+1}(x_{3}).
\]
Therefore, we can inductively have
\begin{eqnarray}
c_{3,1}(x_{3}) c(x_{3})+c_{3,0}(x_{3}) b_{1}(x_{3}) &=&0, \label{eqn:4.18} \\
c_{3,2}(x_{3}) c(x_{3})+c_{3,1}(x_{3}) b(x_{3})+c_{3,0}(x_{3}) a_{1}(x_{3}) &=&0, \label{eqn:4.19} \\
c_{3,j+1}(x_{3}) c(x_{3})+b(x_{3})c_{3,j}(x_{3})
+a(x_{3})c_{3,j-1}(x_{3}) &=&0, \;\;\; j\geq 2. \label{eqn:4.20}
\end{eqnarray}

It follows that $\left\{ c_{3,k}(x_{3}) \right\}$ is the solution
of the second order recursive relation determined by equations
(\ref{eqn:4.18})--(\ref{eqn:4.20}). Since $
b^{2}(x_{3})-4a(x_{3})c(x_{3})=0$, $c_{3,j}(x_{3})$ takes the form
given by equation (\ref{eqn:4.16}). $A_{3}(x_{3})$ and
$B_{3}(x_{3})$ are obtained by using the initial equations:
\begin{eqnarray*}
A_{3}(x_{3})c(x_{3})+c_{3,0}(x_{3})b_{1}(x_{3}) &=&0, \\
\frac{\lbrack A_{3}(x_{3})+B_{3}(x_{3}))c(x_{3})]}{Y_{1}(x_{3})}
+A_{3}(x_{3})b(x_{3})+c_{3,0}(x_{3})a_{1}(x_{3}) &=&0.
\end{eqnarray*}
\QED

We are now ready to prove the main theorem of this section, in
which
\begin{eqnarray}
    A_{1}(x_{dom}) &=& -B_{1}(x_{dom})+\frac{-c_{1,0}(x_{dom})b_{1}(x_{dom})%
}{c(x_{dom})}  \label{eqn:4.29} \\
&=&\left( \frac{h_{1}(x_{dom},Y_{0}(x_{dom}))}{%
a(x_{dom})[Y_{1}(x_{dom})-Y_{0}(x_{dom})]Y_{0}(x_{dom})}-\frac{b_{1}(x_{dom})%
}{c(x_{dom})}\right) c_{1,0}(x_{dom}), \nonumber
\end{eqnarray}
\begin{equation}
    A_{2}(x_{dom})=-\frac{c_{2,0}(x_{dom})b_{1}(x_{dom})}{c(x_{dom})},
\label{eqn:4.31}
\end{equation}
$A_3(x_{dom})$ is given in (\ref{eqn:4.17}),
\begin{equation} \label{eqn:A4}
    A_{4}(x_{dom})=-\frac{b_{1}(x_{dom})c_{0,4}(x_{dom})}{c(x_{dom})},
\end{equation}
\begin{equation}
    B_{1}(x_{dom})=\frac{-h_{1}(x_{dom},Y_{0}(x_{dom}))c_{1,0}(x_{dom})}{%
a(x_{dom})[Y_{1}(x_{dom})-Y_{0}(x_{dom})]Y_{0}(x_{dom})},
\label{eqn:4.28}
\end{equation}
\begin{equation}
    B_{2}(x_{dom})=\frac{-c_{2,0}(x_{dom})h_{1}(x_{dom},Y_{0}(x_{dom}))}{%
aY_0(x_{dom})^{2}},  \label{eqn:4.30}
\end{equation}
and $B_{3}(x_{dom})$ is given in (\ref{eqn:4.17-b}).

\begin{theorem} \label{theorem4.1}
Consider the stable non-singular genus 1 random walk.
Corresponding to the four case, we then have the following tail
asymptotic properties for the joint probabilities $\pi_{n,j}$ for
large $n$.
\bigskip

\textbf{1.}  If $p_{i,j}$ is not X-shaped, then there are four
types of exact tail asymptotics:

\textbf{Case 1:} (Exact geometric decay)
\begin{equation} \label{eqn:n7-1}
    \pi_{n,j} \sim \left[ A_{1}(x_{dom})\left(\frac{1}{Y_{1}(x_{dom})}\right)^{j-1} +
B_{1}(x_{dom})\left( \frac{1}{Y_{0}(x_{dom})}\right)^{j-1}\right]
\left( \frac{1}{x_{dom}}\right)^{n-1}, \;\;\;j\geq 1;
\end{equation}

\textbf{Case 2:} (Geometric decay multiplied by a factor of
$n^{-1/2}$)
\begin{equation} \label{eqn:n7-2}
    \pi_{n,j} \sim \frac{\lbrack A_{2}(x_{dom})+(j-1)B_{2}(x_{dom})]}{\sqrt{\pi}}
\left( \frac{1}{Y_{1}(x_{dom})}\right)^{j-1}n^{-1/2}\left(
\frac{1}{x_{dom}}\right)^{n-1}, \;\;\;j\geq 1;
\end{equation}

\textbf{Case 3:} (Geometric decay multiplied by a factor of
$n^{-3/2}$)

\begin{equation} \label{eqn:n7-3}
    \pi_{n,j}\sim \frac{\lbrack A_{3}(x_{dom})+(j-1)B_{3}(x_{dom})]}{\sqrt{\pi}}
\left( \frac{1}{Y_{1}(x_{dom})}\right)^{j-1}n^{-3/2}\left(
\frac{1}{x_{dom}}\right)^{n-1}, \;\;\;j\geq 1;
\end{equation}

\textbf{Case 4:} (Geometric decay multiplied by a factor of $n$)

\begin{equation} \label{eqn:n7-4}
    \pi_{n,j}\sim \left[ A_{4}(x_{dom})\left(\frac{1}{Y_{1}(x_{dom})}\right)^{j-1}\right] n
\left( \frac{1}{x_{dom}}\right)^{n-1}, \;\;\;j\geq 1.
\end{equation}
\bigskip

\textbf{2.} If $p_{i,j}$ is X-shaped, but both $p_{i,j}^{(1)}$ and $%
p_{i,j}^{(2)}$ are not X-shaped, we then have the following exact
tail asymptotic properties:

\textbf{Case 1:} (Exact geometric decay) It is given by
(\ref{eqn:n7-1});

\textbf{Case 2: }(Geometric decay multiplied by a factor of
$n^{-1/2}$) It is given by (\ref{eqn:n7-2});

\textbf{Case 3:} (Geometric decay multiplied by a factor of
$n^{-3/2}$) It is given by
\begin{equation} \label{eqn:n7-5}
    \pi_{n,j} \sim \frac{\lbrack A_{3}(x_{dom})+(-1)^{n+j}A_{3}(-x_{dom})]}{\sqrt{\pi }}
\left(\frac{1}{Y_{1}(x_{dom})}\right)^{j-1}n^{-3/2}
\left(\frac{1}{x_{dom}}\right)^{n-1}, \;\;\;j\geq 1;
\end{equation}

\textbf{Case 4:} (Geometric decay multiplied by a factor of $n$)
It is given by (\ref{eqn:n7-4}).
\bigskip

\textbf{3.} If $p_{i,j}$ and $p_{i,j}^{(1)}$ are X-shaped, but $%
p_{i,j}^{(2)} $ is not, we then have the following exact tail
asymptotic properties:

\textbf{Case 1:} (Exact geometric decay) When
$\widetilde{x}_{1}<x^{\ast }$, it is given by (\ref{eqn:n7-1});
when $\widetilde{x}_{1}=x^{\ast }=x_{3}$, it is also given by
(\ref{eqn:n7-1}); when $x^{\ast }<\widetilde{x}_{1}$, it is given
by
\begin{equation} \label{eqn:n7-6}
    \pi_{n,j} \sim \left[A_{1}(x_{dom})+(-1)^{n+j}A_{1}(-x_{dom})\right]
\left( \frac{1}{Y_{1}(x_{dom})}\right)^{j-1}
\left(\frac{1}{x_{dom}}\right)^{n-1}, \;\;\;j\geq 1;
\end{equation}

\textbf{Case 2:} (Geometric decay multiplied by a factor of
$n^{-1/2}$) When $x^{\ast }>\widetilde{x}_{1}$, it is given by
(\ref{eqn:n7-2}); when $x^{\ast }< \widetilde{x}_{1}$, it is given
by

\begin{equation} \label{eqn:n7-7}
    \pi_{n,j} \sim \frac{\lbrack A_{2}(x_{dom})+(-1)^{n+j}A_{2}(-x_{dom})]}{\sqrt{\pi }}
\left( \frac{1}{Y_{1}(x_{dom})}\right )^{j-1} n^{-1/2}
\left(\frac{1}{x_{dom}}\right)^{n-1}, \;\;\;j\geq 1;
\end{equation}

\textbf{Case 3:} (Geometric decay multiplied by a factor of
$n^{-3/2}$) It is given by
\[
    \pi_{n,j} \sim \frac{\lbrack A_{3}(x_{dom})+(-1)^{n+j}A_{3}(-x_{dom})]}{\sqrt{\pi }}
\left( \frac{1}{Y_{1}(x_{dom})}\right)^{j-1}n^{-3/2}
\left(\frac{1}{x_{dom}}\right)^{n-1}, \;\;\;j\geq 1;
\]

\textbf{Case 4:} (Geometric decay multiplied by a factor of $n$)
It is given by (\ref{eqn:n7-4}).
\bigskip

\textbf{4.} If $p_{i,j}$ and $p_{i,j}^{(2)}$ are X-shaped, but $%
p_{i,j}^{(1)} $ is not, then it is the symmetric case to
\textbf{3}. All
expression in \textbf{3} are valid after switching $x^{\ast }$ and $%
\widetilde{x}_{1}$.
\bigskip

\textbf{5.} If all $p_{i,j}$, $p_{i,j}^{(1)}$ and $p_{i,j}^{(2)}$
are X-shaped, we then have the following exact tail asymptotic
properties:

\textbf{Case 1:} (Exact geometric decay) When $x^{\ast }\leq
\widetilde{x}_{1}$, it is given by (\ref{eqn:n7-6}); when
$x^{\ast}>\widetilde{x}_{1}$, it is also given by (\ref{eqn:n7-6})
by replacing the dominant singularity $x^{\ast}$ by
$\widetilde{x}_{1}$;

\textbf{Case 2:} (Geometric decay multiplied by a factor of
$n^{-1/2}$) When $x^{\ast}<\widetilde{x}_{1}$, it is given by
(\ref{eqn:n7-7}); when $x^{\ast }> \widetilde{x}_{1}$, it is also
given by (\ref{eqn:n7-7}) by replacing the dominant singularity
$x^{\ast}$ by $\widetilde{x}_{1}$;

\textbf{Case 3:} (Geometric decay multiplied by a factor of
$n^{-3/2}$) It is given by (\ref{eqn:n7-5});

\textbf{Case 4:} (Geometric decay multiplied by a factor of $n$)
It is given by
\[
    \pi_{n,j} \sim \left[A_{4}(x_{dom})+(-1)^{n+j}A_{4}(-x_{dom})\right]
\left( \frac{1}{Y_{1}(x_{dom})}\right)^{j-1} n \left(
\frac{1}{x_{dom}}\right)^{n-1}, \;\;\;j\geq 1.
\]
\end{theorem}

\proof \textbf{1.}

\textbf{Case 1:} It follows from Section~\ref{section3} that
$\lim_{x\rightarrow x_{dom}}\left( 1-\frac{x}{x_{dom}}\right)
\varphi_{0}(x)=c_{0,1}(x_{dom})$.  By the induction and equations
(\ref{eqn:4.1})--(\ref{eqn:4.3}), $\lim_{x\rightarrow
x_{dom}}\left( 1-\frac{x}{x_{dom}}\right)
\varphi_{j}(x)=c_{1,j}(x_{dom})$
 with
\begin{eqnarray*}
c_{1,1}(x_{dom})c(x_{dom})+c_{1,0}(x_{dom})b_{1}(x_{dom}) &=&0, \\
c_{1,2}(x_{dom})c(x_{dom})+c_{1,1}(x_{dom})b(x_{dom})+c_{1,0}(x_{dom})a_{1}(x_{dom}) &=&0, \\
c_{1,j+1}(x_{dom})c(x_{dom})+c_{1,j}(x_{dom})b(x_{dom})+c_{1,j-1}(x_{dom})a(x_{dom})
&=&0,\;\;\; j\geq 2.
\end{eqnarray*}
Since $c_{1,j}(x_{dom})$, $j\geq 0$, satisfies the second order
recursive relation above, it takes the form of
\[
    c_{1,j+1}(x_{dom})=A_{1}(x_{dom}) \left( \frac{1}{Y_{1}(x_{dom})}\right )^{j}
+B_{1}(x_{dom}) \left( \frac{ 1}{Y_{0}(x_{dom})}\right)^{j},
\;\;\;j\geq 0.
\]
To determine $A_1=A_{1}(x_{dom})$ and $B_1=B_{1}(x_{dom})$, we use
the initial equations:
\begin{align}
    (A_{1}+B_{1})c(x_{dom})+c_{1,0}(x_{dom})b_{1}(x_{dom}) &=0, \label{eqn:4.26} \\
    \left[ A_{1}\left( \frac{1}{Y_{1}(x_{dom})}\right) +B_{1}\left(\frac{1}{Y_0(x_{dom})} \right) \right]
c(x_{dom})+(A_{1}+B_{1})b(x_{dom})+c_{1,0}(x_{dom})a_{1}(x_{dom})
&=0. \label{eqn:4.27}
\end{align}
Multiplying both sides of equation (\ref{eqn:4.27}) by
$Y_{0}(x_{dom})$, adding the resulting one to (\ref{eqn:4.26}),
and taking into account
$a(x_{dom})Y_{0}^{2}(x_{dom})+b(x_{dom})Y_{0}(x_{dom})+c(x_{dom})=0$,
$h_{1}(x_{dom},Y_{0}(x_{dom}))=a_{1}(x_{dom})Y_{0}(x_{dom})+b_{1}(x_{dom})$
and $c(x_{dom})=Y_{0}(x_{dom})Y_{1}(x_{dom})a(x_{dom})$ yield:
\begin{eqnarray*}
(A_{1}+B_{1})c(x_{dom})+c_{1,0}(x_{dom})b_{1}(x_{dom}) &=&0, \\
A_{1}\frac{Y_{0}(x_{dom})}{Y_{1}(x_{dom})}%
c(x_{dom})+B_{1}c(x_{dom})+(A_{1}+B_{1})b(x_{dom})Y_{0}(x_{dom})+c_{1,0}(x_{dom})a_{1}(x_{dom})Y_{0}(x_{dom})
&=&0,
\end{eqnarray*}%
which gives (\ref{eqn:4.28}) and (\ref{eqn:4.29}).
So, $B_{1}(x_{dom})=0$ if $x_{dom}=x^{\ast }$ and $B_{1}(x_{dom})\neq 0$ if $%
x_{dom}=\widetilde{x}_{1}$. By the Tauberian-like theorem, we
obtain (\ref{eqn:n7-1}).

\textbf{Case 2:} Similar to that for 1-Case 1. From the proof, we
have (\ref{eqn:4.30}) and (\ref{eqn:4.31}).
%Let $\lim_{x\rightarrow x_{3}}\sqrt{1-\frac{x}{x3}} \varphi_{0}(x)=c_{2,0}(x_{dom})=c_{0,2}(x_{dom})$.
%By the induction and equations (\ref{eqn:4.1})--(\ref{eqn:4.3}),
%$\lim_{x\rightarrow
%x_{3}}\sqrt{1-\frac{x}{x3}}\varphi_{j}(x)=c_{2,j}(x_{dom})$.
%Similarly, as in Case~1, we have
%\[
%    c_{2,j+1}(x_{dom})=[A_{2}+jB_{2}]\left( \frac{1}{Y_{0}(x_{3})}\right)^{j},\;\;\; j\geq 0,
%\]
%since $Y_{0}(x_{3})=Y_{1}(x_{3})$. To determine $A_{2}$ and $B_{2}$, we use
%\begin{eqnarray*}
%    A_{2}c+c_{2,0}(x_{dom})b_{1} &=&0, \\
%    (A_{2}+B_{2})\frac{c}{Y_{0}}+A_{2}b+c_{2,0}(x_{dom})a_{1} &=&0.
%\end{eqnarray*}
%Multiplying the second equation by $Y_{0}$ and then adding the
%resulting one to the first one lead to
%\[
%    (A_{2}+B_{2})c+A_{2}[bY_{0}+c]+c_{2,0}(x_{dom})h_{1}(x_{3},Y_{0})=0.
%\]
%Using $aY_{0}^{2}+bY_{0}+c=0$,
%$h_{1}(x_{dom},Y_{0})=a_{1}Y_{0}+b_{1}$ and $
%c=Y_{0}Y_{1}a=aY_{0}^{2}$, we obtain
%\[
%    (A_{2}+B_{2})aY_{0}^{2}-aA_{2}Y_{0}^{2}+c_{2,0}(x_{dom})h_{1}(x_{3},Y_{0})=0,
%\]
%which gives
%\begin{equation} \label{eqn:4.30}
%B_{2}=\frac{-c_{2,0}(x_{dom})h_{1}(x_{3},Y_{0})}{aY_{0}^{2}}
%\end{equation}
%and
%\begin{equation} \label{eqn:4.31}
%A_{2}=-\frac{c_{2,0}(x_{dom})b_{1}}{c}.
%\end{equation}
%So, $B_{2}=0$ if $x_{dom}=x^{\ast}$ and $B_{2}\neq 0$ if
%$x_{dom}=\widetilde{x}_{1}$. By the Tauberian-like theorem, we
%obtain (\ref{eqn:4.22}).

\textbf{Case 3:} Write
\[
    \varphi_{j}^{\prime}(x)=\sum_{n=0}^{\infty }(n+1)\pi_{n+2,j}x^{n}=\sum_{n=0}^{\infty }(n+1)x_{3}^{n}\pi_{n+2,j}\left( \frac{x}{
x_{3}}\right)^{n}.
\]
According Lemma~\ref{lemma4.2} and the Tauberian-like theorem, we
have
\[
    (n+1)x_{3}^{n}\pi_{n+2,j}\sim \frac{c_{3,j}(x_{3})}{\sqrt{\pi}}n^{-1/2},
\]
which is equivalent to (\ref{eqn:n7-3}).

\textbf{Case 4:} The results can be proved in the same fashion as
in Case~1 and Case~2.

The proofs of the other cases are omitted due to the similarity to
\textbf{1} and Theorem~\ref{theorem4.1-a}. \QED

%\begin{remark}
%For the random walk in 4,
%$h_{2}(-x^{\ast},-Y_{0})=-h_{2}(x^{\ast},Y_{0})$. So, if we choose
%$p_{i,j}^{(0)}$ such that
%$h_{0}(x^{\ast},Y_{0})+h_{0}(-x^{\ast},-Y_{0})=0$, then when
%$n+j+1$ is even (for example, $n$ odd, $j$ even),
%\begin{eqnarray*}
%&&c_{1,0}(x^{\ast})+c_{1,0}(-x^{\ast}) \\
%&=&[\pi_2(Y_{0})-\pi_2(-Y_{0})]\times \frac{
%h_{2}(x^{\ast},Y_{0})h_{1}(x^{\ast},Y_{1})a(x^{\ast})}{x^{\ast}f^{\prime}(x^{\ast})},
%\end{eqnarray*}
%which is equal to zero if $\pi_2(y)$ is an even function.
%\end{remark}

\section{Examples and Concluding Remarks}

In this paper, for a non-singular genus 1 random walk, we proposed
a kernel method to study the exact tail asymptotic behaviour of
the joint stationary probabilities along a coordinate direction,
when the value of the other coordinate is fixed,
 and also the exact tail asymptotic behaviour for the two marginal
 distributions.
A total of four different types of exact tail asymptotics exists.
The fourth one, a geometric decay multiplied by a factor $n$, was
not reported before for this discrete-time model (the same type
was reported recently for a continuous-time random walk model by
Dai and Miyazawa~\cite{Dai-Miyazawa:10}). In this study, we also
revealed a new periodic phenomena for all four types of exact tail
asymptotics when there are two dominant singularities for the
unknown generating function, say $\pi_1(x)$, with the same
asymptotic property at them.

The key idea of this kernel method is simple and the use of the
Tauberian-like theorem greatly simplifies the analysis, which,
unlike in the situation when a standard Tauberian theorem is used,
is also rigorous. Under the assumption that there is only one
dominant singularity, this method provides a straightforward
routine analysis for the exact tail asymptotic behaviour. However,
without this assumption, the analysis is not simple, at least to
our best effort, for telling how many dominant singularities and
when a pole is simple. It is also challenging to characterize the
exact tail asymptotic along a coordinate direction when the value
of the other coordinate is not zero, since it is not a direct
consequence of the kernel method.

This kernel method can also be used for characterizing the exact
tail asymptotics for the non-singular genus 0 case and the
singular random walks (see Li, Tavakoli and
Zhao~\cite{Li-Tavakoli-Zhao:11}). With the detailed analysis
provided in this paper, we expect further research in applying
this kernel method to more general models.

The complete characterization of the exact tail asymptotic
behaviour provided in this paper does not necessarily imply that
for any specific model, a characterization explicitly in terms of
the system parameters exists. However, we are confident that for
any specific model, if using a different method could lead to a
such characterization, in terms of system parameters, then it can
be done using the kernel method.
 Finally, we mention two
examples, which have been analyzed by using the proposed kernel
method.

\textbf{Example 1.} A generalized two-demand model was considered
in Li and Zhao~\cite{Li-Zhao:10b} using the same idea proposed in
this paper. For this model, let $\lambda$ and $\lambda_k$
($k=1,2$) be the Poisson arrival rate with two demands and the
arrival rate of the two dedicated Poisson arrivals, respectively.
Furthermore, let $\mu_k$ ($k=1,2$) be the exponential service
rates of the two independent parallel servers. For a detailed
description of the model, one may refer to \cite{Li-Zhao:10b}. For
this model, the three regions, on which the joint probabilities
along a coordinate direction, say queue 1, have an exact geometric
decay, a geometric decay multiplied by a factor $n^{-1/2}$ and a
geometric decay multiplied by a factor $n^{-3/2}$ are extremely
simple, which are: (a) $\frac{\mu_1}{\lambda+\lambda_1} <
\frac{\mu_2-\lambda_2}{\lambda}$; (b)
$\frac{\mu_1}{\lambda+\lambda_1} =
\frac{\mu_2-\lambda_2}{\lambda}$; and (c)
$\frac{\mu_1}{\lambda+\lambda_1} >
\frac{\mu_2-\lambda_2}{\lambda}$, respectively.
\medskip

\textbf{Example 2.} Consider the simple random walk, or a random
walk for which $p_{i,j}$ and both $p^{(k)}_{i,j}$ ($k=1,2$) are
cross-shaped. We then can follow the general results obtained in
this paper to have refined properties. For example, consider the
case of $M_y>0$ and $M_x<0$ and assume that the system is stable.
Then, along the $x$-direction, $\pi_{n,j}$ has three types exact
asymptotics in the following respective regions:

\textbf{1. Exact geometric:}
\[
    \frac{x_{3}}{x_{3}-1}\left[\sqrt{\frac{p_{0,-1}}{p_{0,1}}}-1\right]
p_{0,1}^{(1)}+p_{1,0}^{(1)}x_{3}>p_{-1,0}^{(1)};
\]

\textbf{2. Geometric with a factor $n^{-1/2}$:}
\[
    \frac{x_{3}}{x_{3}-1}\left[\sqrt{\frac{p_{0,-1}}{p_{0,1}}}-1\right]
p_{0,1}^{(1)}+p_{1,0}^{(1)}x_{3}=p_{-1,0}^{(1)};
\]

 \textbf{3. Geometric with a factor $n^{-3/2}$:}
\[
    \frac{x_{3}}{x_{3}-1}\left[\sqrt{\frac{p_{0,-1}}{p_{0,1}}}-1\right]
p_{0,1}^{(1)}+p_{1,0}^{(1)}x_{3}<p_{-1,0}^{(1)}.
\]

When $M_y<0$ and $M_x<0$, this example also reveals the fourth
type of exact tail asymptotic property, or a geometric decay
multiplied by the factor $n$ along the $x$-coordinate direction in
the region defined by the following conditions:
\begin{equation} \label{eqn:9.1}
    \frac{x_{3}}{x_{3}-1}\left[ \sqrt{\frac{p_{0,-1}}{p_{0,1}}}-1\right]
p_{0,1}^{(1)}+p_{1,0}^{(1)}x_{3}\geq p_{-1,0}^{(1)},
\end{equation}

\begin{equation} \label{eqn:9.2}
    \frac{y_{3}}{y_{3}-1}\left[ \sqrt{\frac{p_{-1,0}}{p_{1,0}}}-1\right]
p_{1,0}^{(2)}+p_{0,1}^{(2)}y_{3}\geq p_{0,-1}^{(2)},
\end{equation}
\begin{eqnarray}
    h_{1}(x^{\ast },\widetilde{y}_0) &=&0,  \\
    \frac{p_{-1,0}}{p_{1,0}} &<&\frac{p_{-1,0}^{(1)}}{p_{1,0}^{(1)}},
\end{eqnarray}
and
\begin{equation}
    \frac{(x^{\ast}-1)p_{0,-1}^{(2)}p_{1,0}+p_{1,0}^{(2)}p_{0,-1}}{(x^{\ast}-1)p_{0,1}^{(2)}p_{1,0}+p_{1,0}^{(2)}p_{0,1}}
=1+\frac{(x^{\ast}-1)[p_{-1,0}^{(1)}-p_{1,0}^{(1)}x^{\ast}]}{p_{0,1}^{(1)}x^{\ast}}.
\end{equation}
Here, $x^{\ast} \in (1,x_{3}]$ and $y^* \in (1,y_{3}]$
are the zero $h_{1}(x,Y_{0}(x))$ and $h_{2}(X_{0}(y),y)$,
respectively, whose existence is guaranteed by
Lemma~\ref{lemma2.7} under conditions (\ref{eqn:9.1}) and
(\ref{eqn:9.2}); $\widetilde{y}_0=Y_{0}(x^{\ast})$ and in this case we
have $\widetilde{y}_0=y^{\ast}$; and $\widetilde{x}_0=X_{0}(Y_{0}(x^{\ast}))$.

It is not very difficult to see this is not an empty region. The
last thing which we need to check is the coefficient
\begin{equation} \label{eqn:9.6}
    c_{0,4}(x_{dom})=\frac{h_{2}(x_{dom},y^{\ast })[h_{1}(\widetilde{x}_0,
y^{\ast })\pi (\widetilde{x}_0) + h_{0}(\widetilde{x}_0,y^{\ast})]\pi_{0,0}]}{x^{\ast2}h_{1}^{\prime}(x_{dom},y^{\ast})Y_{0}^{\prime}(x_{dom})
h_{2}^{\prime}(X_{0}(y^{\ast}),y^{\ast})}\neq 0,
\end{equation}
or
\[
    h_{1}(\widetilde{x}_0,y^{\ast })\pi (\widetilde{x}_0)+h_{0}(\widetilde{x}_0,y^{\ast })\pi _{0,0}\neq 0,
\]
which is true since $h_{2}(x_{dom},y^{\ast })=h_{2}(X_{1}(y^{\ast
}),y^{\ast})>h_{2}(X_{0}(y^{\ast }),y^{\ast })=0$.

 \vspace*{1cm}

\noindent \textbf{Acknowledgements:} The authors thank the
anonymous referee for the valuable comments and suggestions, which
significantly improved the quality of the paper, and the late
Dr.~P. Flajolet of INRIA for the discussion of the Tauberian-like
theorem. This work was supported in part by Discovery Grants from
NSERC of Canada.

\appendix

\section{Proof to Lemmas~\ref{lemma2.2}--\ref{lemma2.6} and Propositions~\ref{theorem2.2}--\ref{theorem2.3}}
\label{appendix1}

\underline{\proof of Lemma~\ref{lemma2.2}.} \textbf{1.} From $h(x,y)=0$, we have
\begin{equation*}
    y^{\prime}=-\frac{a^{\prime}(x)y^{2}+b^{\prime}(x)y+c^{\prime}(x)}{2a(x)y+b(x)}.
\end{equation*}
Using $a(1)+b(1)+c(1)=0$, the property in (\ref{eqn:MxMy}) and the
expression for $Y_k(1)$ in Lemma~\ref{lemma1.1}, we obtain (a) and
(b). (c) is obvious.

\textbf{2.} There are two possible cases: $M_{y}<0$ and $M_{y}>0$. If
$M_{y}<0$, according to the ergodicity condition in Theorem~\ref{theoremergodicity}, $M_{y}^{(1)}M_{x}-M_{y}M_{x}^{(1)}<0$ must hold, which yields
\begin{eqnarray*}
    f^{\prime}(1) &=&a(1)h_{1}^{\prime}(1,Y_{0}(1))h_{1}(1,Y_{1}(1)) \\
    &=&a(1)\left[ a_{1}^{\prime}(1)+a_{1}(1)Y_{0}^{\prime}(1)+b_{1}^{\prime}(1)\right] h_{1}(1,Y_{1}(1)) \\
    &=&\frac{a(1)h_{1}(1,Y_{1}(1))}{-M_{y}}\left[M_{y}^{(1)}M_{x}-M_{y}M_{x}^{(1)}\right] <0.
\end{eqnarray*}
From equation (\ref{eqn:2.1}), $f(x_{3})\geq 0$, it follows that $f(x)=0$ has a root in
$(1,x_{3}]$ since $f(1)=0$ and $f^{\prime}(1)<0$.

If $M_{y}>0$, we have
\begin{eqnarray*}
    f^{\prime}(1) &=&a(1)h_{1}(1,Y_{0}(1))h_{1}^{\prime}(1,Y_{1}(1)) \\
    &=&\frac{-a(1)h_{1}(1,Y_{0}(1))[M_{x}M_{y}^{(1)}-M_{y}M_{x}^{(1)}]}{M_{y}}.
\end{eqnarray*}
If $M_{x}M_{y}^{(1)}-M_{y}M_{x}^{(1)}<0$, from $f(x_{2})\geq 0$,
$f(1)=0$ and $f^{\prime}(1)>0$, $f(x)=0$ has a root in
$[x_{2},1)$. Similarly, if $M_{x}M_{y}^{(1)}-M_{y}M_{x}^{(1)}>0$,
we have $f^{\prime}(1)<0$, which implies that $f(x)=0$ has a root
in $(1,x_{3}]$. Also, 1 is not a repeated root of $f(x)=0$ since
$f^{\prime}(1)\neq 0$ when $M_{y}\neq 0$. \QED
\bigskip

\underline{\proof of Lemma~\ref{lemma2.3}.} \textbf{1.} Suppose
$f(z)=0$. From equation (\ref{eqn:2.1}), we have $F(z)=0$. So we
can write $F(x)=(x-z)G(z)$. Similarly, since $D_{1}(z)=0$ (Recall
$D_{1}(x)=b^{2}(x)-4a(x)c(x)$), we can write
$D_{1}(x)=(x-z)D^{\ast}(x)$, where $D^{\ast}(x)$ is a polynomial.
It follows that $f(x)=(x-z)T(x)$, where
\[
    T(x)=a(x)[a_{1}(x)]^{2}\left\{ (x-z)[G(z)]^{2}-\frac{D^{\ast}(x)}{4a^{2}(x)}\right\}.
\]
 Since the random walk has genus 1, $z$ is not a repeated
root of $D_{1}(x)=0$, which implies $a(x)[a_{1}(x)]^{2}\frac{D^{\ast}(z)}{4a^{2}(x)}\neq 0$ (note that $a(z)\neq 0$ since $D_{1}(z)=0$ and $b(z)>0$
when $z<0$). It follows that $T(z)\neq 0$, that is, $z$ is not a repeated
root of $f(x)=0$.

\textbf{2.} This is a direct consequence of equation
(\ref{eqn:2.1}).

\textbf{3.} Suppose $x^{\prime}$ is a common root. If
$a(x^{\prime})\neq 0$, it is easy to obtain that $x^{\prime}$ is a
branch point. Assume $a(x^{\prime})=0$. Clearly, $x^{\prime}$
cannot be a positive number. Since
$\widetilde{f}_{1}(x^{\prime})=a_{1}(x^{\prime})[-2b(x^{\prime})]=0$,
$f_{0}(x^{\prime})=\frac{a_{1}(x^{\prime})c(x^{\prime})}{-b(x^{\prime})}+b_{1}(x^{\prime})=0$
and $b(x^{\prime})\neq 0$, we obtain $a_{1}(x^{\prime})=0$ and
$b_{1}(x^{\prime})=0$, which implies that $x^{\prime}=0$ since
$b_{1}(x)$ has only nonnegative zeros.

\textbf{4.} Let $-|z|$ be a negative root of $f_{0}(x)=0$ in
$[-x_{3},-1)$. From the definition of $f_{0}(x)$,  we have
$\sum_{i\geq -1,j\geq 0}p_{i,j}^{(1)}[-|z|]^{i}Y_{0}^{j}(-|z|)=1$,
which implies $f_{0}(|z|)>0$ since $Y_{0}(|z|)>|Y_{0}(-|z|)|$.
According to $f_{0}(1)\leq 0$ and Lemma~\ref{lemma2.2}-1,
$f_{0}(x)=0$ has a root, say $z^{\prime }$ in $(1,|z|)$. Again,
from $Y_{0}(|z^{\prime }|)>|Y_{0}(-|z^{\prime }|)|$,
$f_{0}(-|z^{\prime}|)<0$, which implies $f_{0}(x)=0$ has a root in
$(-|z^{\prime }|,-1)$ since $f_{0}(-1)>0$. Clearly, this root is
greater than $-|z|$.

\textbf{5.} Let $|x|\in (1,x_{3}]$. Since $-b(-|x|)<0$,
$Y_{0}(-|x|)=\frac{-b(-|x|)+\sqrt{D_{1}(x)}}{a(-|x|)}
=\frac{2c(-|x|)}{-b(-|x|)-\sqrt{D_{1}(-|x|)}}$. From $b(-|x|)\geq
b(|x|)$,  $\sqrt{D_{1}(-|x|)}> \sqrt{D_{1}(|x|)}$ and
$|c(-|x|)|\leq c(|x|)$, we obtain $|Y_{0}(-|x|)|<
\frac{2c(|x|)}{b(|x|)+\sqrt{D_{1}(|x|)}}=Y_{0}(|x|)$. \QED

\underline{\proof of Lemma~\ref{lemma2.4}.} Assume $|z|=1$ and
$z\neq 1$ or $-1$. From Lemma~\ref{lemma1.1-b}-1, $|Y_{0}(z)|<1$.
Since
\[
h_{1}(x,y)=a_{1}(x)y+b_{1}(x)=x\left( \sum_{i\geq -1,j\geq
0}p_{i,j}x^{i}y^{j}-1\right)
\]
and when $|z|=1$, $\left\vert \frac{1}{z}\right\vert =1$ as well,
we obtain
\[
\left\vert \sum_{i\geq -1,j\geq
0}p_{i,j}z^{i}Y_{0}(z)^{j}\right\vert \leq \sum_{i\geq -1,j\geq
0}p_{i,j}|z^{i}||Y_{0}(z)|^{j}<1,
\]
which yields $f_{0}(x)=h_{1}(z,Y_{0}(z))\neq 0$.

For $z=-1$, $|Y_{0}(-1)|<1$ if $p_{i,j}$ is not X-shaped, and
$Y_{0}(-1)=-1$ if $p_{i,j}$ is X-shaped and $p_{i,j}^{(1)}$ is not
X-shaped. It follows that $f_{0}(-1)>0$ in both cases since
$b_{1}(-1)\geq | a_{1}(-1)|$ and $b_{1}(-1)>0$  in the first case
and $b_{1}(-1)>|a_{1}(-1)|$ in the second case. \QED
\bigskip

For special ransom walk 1, we have
\begin{equation}
    a_{1}(x)=p_{0,1}^{(1)}x+p_{1,1}^{(1)}x^{2}\; \text{ and }\;
    b_{1}(x)=-x+p_{1,0}^{(1)}x^{2},
\end{equation}
\begin{equation} \label{eqn:2.4}
    a(x)=p_{0,1}x, b(x)=p_{-1,0}-x+p_{1,0}x^{2}\; \text{ and }\;
    c(x)=p_{0,-1}x.
\end{equation}
In this case, $f(x)$ becomes
\begin{equation} \label{eqn:2.5}
    f(x) = x^{2}f^{\ast}(x),
\end{equation}
where
\begin{eqnarray*}
    f^{\ast}(x) &=&d_{4}^{\ast}x^{4}+d_{3}^{\ast}x^{3}+d_{2}^{\ast}x^{2}+d_{1}^{\ast}x+d_{0}^{\ast} \\
    &=&p_{0,1}x[1-p_{1,0}^{(1)}x]^{2}+[p_{-1,0}-x+p_{1,0}x^{2}][1-p_{1,0}^{(1)}x][p_{0,1}^{(1)}+p_{1,1}^{(1)}x]+p_{0,-1}x[p_{0,1}^{(1)}+p_{1,1}^{(1)}x]^{2}
\end{eqnarray*}
with
\begin{equation*}
    d_{0}^{\ast}=p_{-1,0}p_{0,1}^{(1)}\text{ and }d_{4}^{\ast}=-p_{1,0}p_{1,0}^{(1)}p_{1,1}^{(1)}.
\end{equation*}

\underline{\proof of Proposition~\ref{theorem2.2}.} \textbf{1.}
Obviously, From equation (\ref{eqn:2.5}) and Lemma~\ref{lemma2.2},
$f(x)=0$ has at least four real roots with two in $[x_{2},x_{3}]$
and two equal to zero. The facts that $f(x_{1})\geq 0$,
$f(x_{4})\geq 0$ and $f(\pm \infty )=-\infty$ yield one root in
$(-\infty ,x_{1}]$ and another root in $ [x_{4},+\infty )$.

\textbf{2.} It is a direct result of
Proposition~\ref{theorem2.2}-1
 and Lemma~\ref{lemma2.3}-4.
 \QED
\bigskip

For the random walk considered in Theorem~\ref{theorem2.1}-2 (or both $p_{i,j}$ and $p^{(1)}_{i,j}$
are X-shaped, we have
\begin{equation}
    a_{1}(x)=p_{-1,1}^{(1)}+p_{1,1}^{(1)}x^{2},\;\;\; b_{1}(x)=-x,
\end{equation}
\begin{equation}
    a(x)=p_{-1,1}+p_{1,1}x^{2}, \;\;\; b(x)=-x,\;\;\;c(x)=p_{-1,-1}+p_{1,-1}x^{2}.
\end{equation}
Therefore, $f(x)$ becomes
\begin{eqnarray}
    f(x) &=& a(x)b_{1}^{2}(x)-b(x)b_{1}(x)a_{1}(x)+c(x)a_{1}^{2}(x)  \notag \\
    &=&x^{2}[p_{-1,1}+p_{1,1}x^{2}]-x^{2}[p_{-1,1}^{(1)}+p_{1,1}^{(1)}x^{2}]+[p_{-1,-1}+p_{1,-1}x^{2}][p_{-1,1}^{(1)}+p_{1,1}^{(1)}x^{2}]^{2} \notag \\
    &=&d_{6}x^{6}+d_{4}x^{4}+d_{2}x^{2}+d_{0}, \label{eqn:2.8}
\end{eqnarray}
where
\begin{eqnarray*}
    d_{6} &=&\left[ p_{1,1}^{(1)}\right] ^{2}p_{1,-1}, \\
    d_{4} &=&p_{1,1}-p_{1,1}^{(1)}+2p_{1,-1}p_{-1,1}^{(1)}p_{1,1}^{(1)}+p_{-1,-1}\left[ p_{1,1}^{(1)}\right] ^{2}, \\
    d_{2} &=&p_{-1,1}-p_{-1,1}^{(1)}+2p_{-1,-1}p_{-1,1}^{(1)}p_{1,1}^{(1)}+p_{1,-1}\left[ p_{-1,1}^{(1)}\right] ^{2}, \\
    d_{0} &=&\left[ p_{-1,1}^{(1)}\right] ^{2}p_{-1,-1}.
\end{eqnarray*}

\underline{\proof of Lemma~\ref{lemma2.6}} $f(1)=f(-1)=0$ follows
from $Y_{i}(1)=-Y_{i}(-1)$, $a_{1}(1)=a_{1}(-1)$ and
$b_{1}(1)=-b_{1}(-1)$. From Lemma~\ref{lemma2.3}, there exists an
$x_0 \neq 1$, $x_0 \in [x_{2},x_{3}]$ such that $f(x_0)=0$. We
provide details for the case of $f_{1}(z)=0$ and a similar proof
can be found for the other case. Since $x_{2}< x_0 \leq x_{3}$,
$-b(z)=-b(-z)>0$. It follows that
\begin{equation}
    Y_{1}(x_0)=\frac{-b(x_0)}{2a(x_0)}+\frac{\sqrt{b^{2}(x_0)-4a(x_0)c(x_0)}}{2a(x_0)}
\end{equation}
and
\begin{equation}
    f_{1}(x_0)=a_{1}(x_0)Y_{1}(x_0)+b_{1}(x_0)=0.
\end{equation}
On the other hand, from $-x_{3}\leq -x_0<-1$ we have $-b(-x_0)<0$, which yields
\begin{equation}
    Y_{1}(-x_0)=\frac{b(x_0)}{2a(x_0)}-\frac{\sqrt{b^{2}(x_0)-4a(x_0)c(x_0)}}{2a(x_0)} = -Y_{1}(x_0)
\end{equation}
and
\begin{equation*}
    f_{1}(-x_0)=a_{1}(-x_0)Y_{1}(-x_0)+b_{1}(-x_0)=-f_{1}(x_0)=0.
\end{equation*}
It follows from equation (\ref{eqn:2.8}) that $f(x)$ can be written as
\begin{equation*}
    f(x)=d_{6}(x^{2}-1)(x^{2}-x_0^{2})(x^{2}+\eta ).
\end{equation*}
Since $\frac{d_{0}}{d_{6}}>0$, we have $\eta >0$, which indicates that $f(x)=0$
has two complex roots.
\QED
\bigskip

\underline{\proof of Proposition~\ref{theorem2.3}.} Suppose that
one of the two complex roots is a root of $f_{0}(x)=0$. First
assume $\frac{d_{0}}{d_{6}}\leq 1$. Then, $z^{2}\eta =
\frac{d_{0}}{d_{6}}\leq 1$ implies $|\eta |<1$. In the case of
$\frac{d_{0}}{d_{6}}>1$, we choose a path $\ell$ to connect the
random walk here to the one with $\frac{d_{0}}{d_{6}}\leq 1$. Then
on $\ell $, the two complex roots of $f(x)=0$ have to pass through
the unit circle, which is impossible according to
Remark~\ref{remark2.1} and Lemma~\ref{lemma2.4}.

%\red{(Hui: Lemma~\ref{lemma2.3} in your revision list, but I thought it should be Lemma~\ref{lemma2.4}. Please check!)} \QED

%\section{Proof of Theorem~\ref{theorem3.1}}
%\label{appendix2}

%\section{Proof to Theorem~\ref{theorem4.1}}
%\label{appendix3}


\begin{thebibliography}{99}

\bibitem{Abate-Whitt:97} Abate, J. and Whitt, W. (1997)
Asymptotics for $M/G/1$ low-priority waiting-time tail
probabilities, \textit{Queueing Systems}, \textbf{25}, 173--233.

\bibitem{Adan-Foley-McDonald:09} Adan, I., Foley, R.D. and McDonald,
D.R. (2009) Exact asymptotics for the stationary distribution of a
Markov chain: a production model, \textit{Queueing Systems},
\textbf{62}, 311--344.

\bibitem{B-BM-D-F-G-GB:02}
Banderier, C., Bousquet-M\'{e}lou, M,  Denise, A, Flajolet, P.,
Gardy, D. and Gouyou-Beauchamps, D. (2002) Generating functions of
generating trees, \textit{Discrete Math.}, \textbf{246}, 29--55.


\bibitem{Bender:1974} Bender, E. (1974)
Asymptotic methods in enumeration, \textit{SIAM Review},
\textbf{16}, 485--513.

\bibitem{Borovkov-M:01}
Borovkov, A.A. and Mogul'skii, A.A. (2001) Large deviations for
Markov chains in the positive quadrant, \textit{Russian Math.
Surveys}, \textbf{56}, 803--916.

\bibitem{Bousquet-Melou:05} Bousquet-M\'{e}lou, M. (2005) Walks in the quarter plane: Kreweras' algebraic model, \textit{Annals
of Applied Probability}, \textbf{15}, 1451-–1491.

\bibitem{cohen-boxma:83} Cohen, J. W. and Boxma, O. J. (1983) \textit{Boundary Value
Problems in Queueing System Analysis}, North-Holland, Amsterdam.

\bibitem{Dai-Miyazawa:10} Dai, J. and Miyazawa, M. (2010)
Reflecting Brownian motion in two dimensions: Exact asymptotics
for the stationary distribution, submitted.

\bibitem{FI:1979} Fayolle, G. and Iasnogorodski, R.  (1979)
Two coupled processors: the reduction to a
Riemann-Hilbert problem, \textit{Z. Wahrscheinlichkeitsth}, \textbf{47}, 325--351.

\bibitem{FKM:82}
Fayolle, G., King, P.J.B. and  Mitrani, I. (1982)
The solution of certain two-dimensional Markov models,
 \textit{Adv. Appl. Prob.}, \textbf{14}, 295--308.

\bibitem{FIM:1999} Fayolle, G., Iasnogorodski, R. and Malyshev, V. (1999)
\textit{Random Walks in the Quarter-Plane}, Springer, New York.

\bibitem{FO:1990} Flajolet, P. and Odlyzko, A. (1990)
Singularity analysis of generating functions, \textit{SIAM J.
Disc. Math.}, \textbf{3}, 216--240.

\bibitem{Flajolet-Sedgewick:09} Flajolet, F. and Sedgewick, R.
(2009) \textit{Analytic Combinatorics}, Cambridge University
Press.


\bibitem{Flatto-McKeqn:77}
Flatto, L. and McKean, H.P. (1977)
Two queues in parallel, \textit{Comm. Pure Appl. Math.}, \textbf{30}, 255--263.

\bibitem{Flatto-Hahn:84}
Flatto, L. and Hahn, S. (1984)
Two parallel queues created by arrivals with two demands I,
\textit{SIAM J. Appl. Math.}, \textbf{44}, 1041--1053.

\bibitem{Flatto:85}
Flatto, L. (1985)
Two parallel queues created by arrivals with two demands II,
\textit{SIAM J. Appl. Math.}, \textbf{45}, 861--878.

\bibitem{Foley-McDonald:01} Foley, R.D. and McDonald, D.R. (2001)
Join the shortest queue: stability and exact asymptotics, \textit{Annals of Applied Probability}, \textbf{11}, 569--607.

\bibitem{Foley-McDonald:05a} Foley, R.D. and McDonald, R.D. (2005) Large deviations of a modified Jackson network:
stability and rough asymptotics, \textit{Annals of Applied Probability}, \textbf{15}, 519-–541.

\bibitem{Foley-McDonald:05b} Foley, R.D. and McDonald, R.D. (2005)
Bridges and networks: exact asymptotics, \textit{Annals of Applied Probability}, \textbf{15}, 542–-586.

\bibitem{Guillemin-Leeuwaarden:09} Guillemin, F. and Leeuwarden, J. (2009)
Rare event asymptotics for a random walk in the quarter plane,
submitted.

\bibitem{Haque:03} Haque, L. (2003)
\textit{Tail Behaviour for Stationary Distributions for Two-Dimensional Stochastic Models}, Ph.D. Thesis, Carleton
University, Ottawa, ON, Canada.

\bibitem{Haque-Liu-Zhao:05}
Haque, L., Liu, L. and Zhao, Y.Q. (2005) Sufficient conditions for a
geometric tail in a QBD process with countably many levels and
phases, \textit{Stochastic Models}, \textbf{21}(1), 77--99.

\bibitem{He-Li-Zhao:08} He, Q., Li, H. and Zhao, Y.Q. (2009)
Light-tailed behaviour in QBD process with countably many phases,
\textit{Stochastic Models}, \textbf{25}, 50--75.

%\bibitem{Hou-Mansour:08} Hou, Q.-H. and Mansour, T. (2008)
%Kernel method and linear recurrence system,
%\textit{Journal of Computational and Applied Mathematics}, \textbf{216}, 227-–242.

\bibitem{Khanchi:08} Khanchi, Aziz (2008) \textit{State of a network when one node overloads},
Ph.D. Thesis, University of Ottawa.

\bibitem{Khanchi:09} Khanchi, Aziz (2009) Asymptotic hitting distribution for a reflected random walk in the
positive quardrant, \textit{Stochastic Models}, \textbf{27}, 169--201.

\bibitem{Kobayashi-Miyazawa:2011} Kobayashi, M. and  Miyazawa, M. (2011)
Tail asymptotics of the stationary distribution of a two dimensional reflecting random walk with unbounded upward jumps,
submitted.

\bibitem{Kobayashi-Miyazawa-Zhao:10} Kobayashi, M.,  Miyazawa, M. and Zhao, Y.Q. (2010)
Tail asymptotics of the occupation measure for a Markov additive
process with an M/G/1-type background process, accepted by
\textit{Stochastic Models}.

\bibitem{Knuth:69} Knuth, D.E. (1969)
\textit{The Art of Computer Programming, Fundamental Algorithms},
vol. 1, second ed., Addison-Wesley.

\bibitem{KST:04} Kroese, D.P., Scheinhardt, W.R.W. and
Taylor, P.G. (2004) Spectral properties of the tandem Jackson
network, seen as a quasi-birth-and-death process, \textit{Annals
of Applied Probability},  \textbf{14}(4), 2057--2089.

\bibitem{Kurkova-Suhov:03}
Kurkova, I.A. and  Suhov, Y.M. (2003)
Malyshev's theory and JS-queues. Asymptotics of stationary probabilities,
\textit{The Annals of Applied Probability},
\textbf{13}, 1313-–1354.


\bibitem{Li-Miyazawa-Zhao:07} Li, L., Miyazawa, M. and Zhao, Y. (2007) Geometric decay in a QBD process with countable background
states with applications to a join-the-shortest-queue model,
\textit{Stochastic Models}, \textbf{23}, 413–-438.

\bibitem{Li-Tavakoli-Zhao:11} Li, H., Tavakoli, J. and Zhao, Y.Q. (2012)
Analysis of exact tail asymptotics for singular random walks in
the quarter plane, accepte by \textit{Queueing Systems}.


\bibitem{Li-Zhao:05}
Li, H. and Zhao, Y.Q. (2005) A retrial queue with a constant
retrial rate, server break downs and impatient customers,
\textit{Stochastic Models}, \textbf{21}, 531--550.

\bibitem{Li-Zhao:09} Li, H. and Zhao, Y.Q. (2009)
Exact tail asymptotics in a priority queue---characterizations of
the preemptive model, \textit{Queueing Systems}, \textbf{63},
355--381.

\bibitem{Li-Zhao:10} Li, H. and Zhao, Y.Q. (2011)
Exact tail asymptotics in a priority queue---characterizations of
the non-preemptive model, \textit{Queueing Systems}, \textbf{68}, 165--192.

\bibitem{Li-Zhao:10b} Li, H. and Zhao, Y.Q. (2011)
Tail asymptotics for a generalized two demand queueing model
--- A kernel method, \textit{Queueing Systems}, \textbf{69}, 77--100.


\bibitem{LM:2008} Lieshout, P. and Mandjes, M. (2008)
Asymptotic analysis of L\'{e}vy-driven tandem queues,
\textit{Queueing Systems}, \textbf{60}, 203--226.

\bibitem{Liu-Miyazawa-Zhao:08} Liu, L., Miyazawa, M. and Zhao, Y.Q. (2008)
Geometric decay in level-expanding QBD models, \textit{Annals of
Operations Research}, \textbf{160}, 83--98.

%\bibitem{Mandjes:07} Mandjes, M. (2007)
%\textit{Large Deviations for Gaussian Queues: Modelling Communication Networks}, Wiley, Hoboken, NJ.

%\bibitem{MWB:07}  Maertens, T.,  Walraevens, J. and  Bruneel H. (2007)
%Priority queueing systems: from probability generating functions
%to tail probabilities, \textit{Queueing Systems}, \textbf{55}, 27--39.

\bibitem{Malyshev:72}
Malyshev, V.A. (1972)
An analytical method in the theory of two-dimensional positive random walks, \textit{Siberian Math. Journal}, \textbf{13}, 1314--1329.

\bibitem{Malyshev:73}
Malyshev, V.A. (1973)
Asymptotic behaviour of stationary probabilities for two dimensional positive random walks, \textit{Siberian Math. Journal}, \textbf{14}, 156--169.

\bibitem{McDonald:99} McDonald, D.R. (1999) Asymptotics of first
passage times for random walk in an orthant, \textit{Annals of
Applied Probability}, \textbf{9}, 110--145.


\bibitem{Mishna:06} Mishna, M. (2009)
Classifying lattice walks restricted to the quarter plane,
\textit{Journal of Combinatorial Theory, Series A}, \textbf{116},
460--477.

\bibitem{Miyazawa:04} Miyazawa, M. (2004) The Markov renewal approach to
$M/G/1$ type queues with countably many background states,
\textit{Queueing Systems}, \textbf{46}, 177--196.

\bibitem{Miyazawa:07} Miyazawa, M. (2007) Doubly QBD process and a solution to the tail decay rate problem,
in \textit{Proceedings of the Second Asia-Pacific Symposium on Queueing Theory and
Network Applications}, Kobe, Japan.

\bibitem{Miyazawa:09} Miyazawa, M. (2009) Two sided DQBD process and solutions to the tail decay rate
problem and their applications to the generalized join shortest queue,
in \textit{Advances in Queueing Theory and Network Applications}, edited by W. Yue, Y. Takahashi and H. Takaki, 3--33,
Springer, New York.

\bibitem{Miyazawa:08}
Miyazawa, M. (2009) Tail decay rates in double QBD processes and
related reflected random walks, \textit{Math. OR}, \textbf{34}, 547--575.

\bibitem{Miyazawa:11} Miyazawa, M. (2011) Light tail asymptotics in multidimensional reflecting processes
for queueing networks, \textit{TOP}, online first at
http://www.springerlink.com/content/120409/?Content+Status=Accepted.

\bibitem{Miyazawa-Rolski:09} Miyazawa, M. and Rolski, T. (2009)
Tail asymptotics for a L\'{e}vey-driven tandem queue with an
intermediate input, \textit{Queueing Systems}, \textbf{63},
323--353.

\bibitem{Miyazawa-Zhao:04}
Miyazawa, M. and Zhao, Y.Q. (2004) The stationary tail asymptotics
in the $GI/G/1$ type queue with countably many background states,
\textit{Adv. in Appl. Probab.},  \textbf{36}(4), 1231--1251.

\bibitem{Morrison:07} Morrison, J.A. (2007)
Processor sharing for two queues with vastly different rates,
\textit{Queueing Systems}, \textbf{57}, 19-–28.

\bibitem{Motyer-Taylor:06} Motyer, Allan J. and Taylor, Peter G. (2006)
Decay rates for quasi-birth-and-death process with countably many
phases and tri-diagonal block generators, \textit{Advances in
Applied Probability}, \textbf{38}, 522--544.

\bibitem{Raschel:10} Raschel, K. (2010)
Green functions and Martin compactification for killed random
walks related to SU(3), \textit{Elect. Comm. in Probab.},
\textbf{15}, 176--190.

\bibitem{TFM:01} Takahashi, Y., Fujimoto, K. and Makimoto, N. (2001)
Geometric decay of the steady-state probabilities in a
quasi-birth-and-death process with a countable number of phases,
\textit{Stochastic Models}, \textbf{17}(1), 1--24.

\bibitem{Tang-Zhao:08} Tang, J. and Zhao, Y.Q. (2008)
Stationary tail asymptotics of a tandem queue with feedback,
\textit{Annals of Operations Research}, \textbf{160}, 173-189.

%\bibitem{vanUitert:03} van Uitert, M.J.G. (2003)
%\textit{Generalized Processor Sharing Queues}, Ph.D. thesis,
%Eindhoven University of Technology, Eindhoven, The Netherlands.

%\bibitem{Wischik:01} Wischik, D. (2001)
%Sample path large deviations for queues with many inputs,
%\textit{Annals of Applied Probability}, \textbf{11}, 379-–404.

\bibitem{Wright:92} Wright, P. (1992) Two parallel processors with coupled inputs, \textit{Adv. Appl. Prob.}, \textbf{24}, 986--1007.
\end{thebibliography}
\end{document}